\documentclass[12pt,a4paper]{article}

\usepackage[margin=1in]{geometry}
\usepackage{graphicx}
\usepackage{float}
\usepackage{booktabs}
\usepackage{amsmath}
\usepackage{caption}
\usepackage[utf8]{inputenc}
\usepackage{amssymb}
\usepackage{hyperref}
\usepackage{amsmath}
\usepackage{times}
\usepackage{xcolor}
\usepackage{authblk}
\usepackage{tabularx}
\usepackage{booktabs}
\usepackage{subcaption}
\usepackage{multirow}
\usepackage{placeins}

\hypersetup{
    colorlinks=true,
    linkcolor=blue,
    urlcolor=blue,
    citecolor=blue
}

\providecommand{\keywords}[1]
{
  \small	
  \textbf{\textit{Keywords---}} #1
}

\title{{Learning constitutive laws under explicit strain limits:
An interpretable strain limiting elasticity--Kolmogorov Arnold neural network framework}}
\author[1]{Chandana Pati}
\author[1,$\dagger$]{S. M. Mallikarjunaiah}
\affil[1]{Department of Mathematics \& Statistics, Texas A\&M University-Corpus Christi, Texas- 78412, USA}
\affil[$\dagger$]{Corresponding author}
\affil[ ]{\textit{E-mail addresses:} \texttt{cpati@islander.tamucc.edu} (Chandana Pati), \texttt{M.Muddamallappa@tamucc.edu} (S.M. Mallikarjunaiah)}
\date{}

\usepackage{algorithm}
\usepackage{algpseudocode}
\usepackage[utf8]{inputenc}

\begin{document}
 
\maketitle

\begin{abstract}
A physically consistent framework for modeling materials in which deformation saturates under increasing stress, such as elastomers, biological tissues, and soft polymeric solids, is provided by strain-limiting elasticity. By the enforcement of bounded strain and vanishing tangent stiffness at large stresses, fundamental limitations of classical elasticity are addressed; however, significant challenges for data-driven learning and experimental calibration are posed by the strongly nonlinear and asymptotic character of these laws. In this work, an interpretable hybrid constitutive modeling framework integrating strain-limiting elasticity (SLE) with Kolmogorov--Arnold Networks (KANs) is proposed to balance mechanical admissibility with controlled data-driven flexibility. The dominant nonlinear elastic response is captured by the SLE model serving as a physics-based backbone, while smooth, low-amplitude residual corrections in stress space are learned exclusively via a KAN. Essential mechanical principles including odd symmetry, monotonic response, bounded strain, and vanishing tangent modulus are embedded directly into the model structure, by which physically admissible behavior is ensured across all stress regimes. The framework is initially assessed on synthetic benchmarks spanning moderate and strong strain-limiting regimes, where near-exact recovery is achieved in smooth regimes and physical consistency is retained under sharply nonlinear transitions. Subsequently, the proposed approach is applied to the classical Treloar rubber elasticity experiments under uniaxial, biaxial, and planar deformation modes. In this experimental setting, agreement with measured stress--stretch responses is systematically improved by the hybrid SLE--KAN formulation, while explicit strain-limiting constraints are preserved. A transparent trade-off between data fidelity and mechanical admissibility is revealed by a regime-based analysis, demonstrating that deviations observed under stronger strain limits are caused by deliberately imposed physical restrictions rather than unconstrained model expressivity. Overall, a transparent and physics-consistent alternative to black-box neural networks for constitutive modeling is provided by the SLE--KAN framework, which is of particular relevance to experimental solid mechanics applications where interpretability, robustness, and strict adherence to physical limits are essential.
\end{abstract}

\noindent\keywords{Strain-limiting elasticity; Kolmogorov--Arnold Networks; Interpretable machine learning; Constitutive modeling; Rubber elasticity; Physics-informed learning}

\section{Introduction}

Accurate constitutive modeling is central to the predictive simulation of material
behavior in solid mechanics, underpinning applications ranging from structural
engineering and soft robotics to biomechanics and polymer science
\cite{marsden1994mathematical,truesdell1952mechanical}. Classical elastic and
hyperelastic theories have been widely adopted due to their analytical
tractability and computational efficiency
\cite{rivlin1948large,arruda1993three}. However, these formulations typically
assume that strain may grow without bound as the applied stress increases. While
this assumption is adequate for small to moderate deformations, it becomes
physically questionable for materials such as elastomers, polymers, biological
tissues, and soft composite solids, which are routinely subjected to large strains
in both natural and engineered settings
\cite{treloar1975physics,boyce2000constitutive}.

A substantial body of experimental evidence demonstrates that, beyond a certain
stress level, incremental increases in stress produce progressively smaller
increments of strain. This phenomenon, commonly referred to as
\emph{strain saturation} or \emph{strain-limiting behavior}, reflects fundamental
microstructural mechanisms, including molecular chain alignment, network
entanglement, chain locking, and geometric constraints within the material
architecture
\cite{treloar1944stress,jones1975properties,deam1976theory}. Once these mechanisms
become dominant, the material’s ability to deform is intrinsically limited.
Constitutive models that neglect such effects may predict unrealistically large
deformations, experience loss of ellipticity, or suffer from numerical
instabilities when extrapolated beyond their calibration regime
\cite{gent1996new,horgan2004constitutive}.

Strain-limiting elasticity (SLE) was introduced to address these deficiencies by
explicitly enforcing bounded strain responses together with a vanishing tangent
modulus at large stresses
\cite{rajagopal2011non,bulivcek2014elastic,rajagopal2011modeling,gou2015modeling}. By construction, SLE models ensure
physical admissibility across all stress levels and provide improved mathematical
well-posedness for boundary value problems involving extreme loading
\cite{rajagopal2003implicit,bulivcek2015analysis}. These properties make
strain-limiting formulations particularly attractive for modeling soft materials
and biological tissues. Nevertheless, the strong nonlinearity and asymptotic
character of strain-limiting constitutive laws pose practical challenges.
Analytical treatment becomes difficult, and purely data-driven approximation is
complicated by the need to simultaneously capture near-linear behavior at small
stresses and saturation at large stresses, especially in the presence of
experimental noise and material variability
\cite{bulivcek2015existence,itou2017nonlinear,mallikarjunaiah2015direct,yoon2021quasi,yoon2022finite}.

In parallel, machine learning techniques particularly neural networks have gained
increasing attention as flexible tools for constitutive modeling
\cite{shen2004neural,abdolazizi2024viscoelastic}. While conventional black-box
neural networks can approximate complex stress--strain relationships, they often
struggle to respect essential physical constraints such as symmetry,
monotonicity, bounded deformation, and asymptotic stiffness degradation. Even when
trained to high accuracy within a limited data range, such models may exhibit
nonphysical behavior under extrapolation, thereby limiting their reliability in
mechanics applications where predictive robustness is critical
\cite{lee2025constitutive}.

Motivated by these limitations, we propose an interpretable hybrid constitutive
modeling framework that integrates strain-limiting elasticity with
Kolmogorov--Arnold Networks (KANs)
\cite{liu2024kan}. The central idea is to employ the SLE model as a
physics-based backbone that captures the dominant nonlinear elastic response,
while using a KAN exclusively to learn smooth, low-dimensional residual
corrections from data. This separation of roles ensures that the primary
constitutive behavior remains governed by physically grounded principles, while
data-driven components provide controlled flexibility to account for deviations
arising from model idealization or experimental uncertainty
\cite{abdolazizi2025constitutive}.

The proposed framework embeds essential mechanical structure including odd
symmetry, monotonic response, bounded strain, and vanishing tangent modulus directly
into the constitutive representation, rather than relying on implicit learning or
purely penalty-based regularization
\cite{rajagopal2003implicit,liu2024kan}. Moreover, the use of KANs enables a
transparent internal representation in which learned parameters admit direct
physical interpretation, in contrast to the distributed and opaque
parameterizations typical of conventional dense neural networks
\cite{ji2024comprehensive,somvanshi2024survey}.

The framework is evaluated through a combination of synthetic benchmarks and
classical experimental data. Synthetic studies spanning moderate and strong
strain-limiting regimes are used to systematically assess approximation accuracy,
interpretability, and representational limitations. The approach is then applied
to the Treloar rubber elasticity experiments under uniaxial, biaxial, and planar
deformation
\cite{treloar1944stress,treloar1975physics}, providing a stringent test in the
presence of experimental noise and material heterogeneity. A regime-based
analysis highlights how the hybrid formulation balances data fidelity with
mechanical admissibility across varying levels of strain limitation.

By integrating physics-based constitutive modeling with interpretable machine
learning, the proposed SLE--KAN framework provides a transparent, robust, and
physically consistent alternative to black-box neural network approaches
\cite{abdolazizi2025constitutive,panahi2025data,pati2025neural}. The methodology is particularly
relevant for experimental solid mechanics applications, where predictive accuracy
must be accompanied by interpretability and strict adherence to fundamental
mechanical principles.

\section{Strain-Limiting Elasticity Model}

\subsection{Motivation and Physical Background}

Having established the experimental prevalence of strain saturation in soft
materials, we now consider its implications for constitutive modeling. Classical
elastic and hyperelastic formulations lack an intrinsic mechanism to restrict
strain growth and therefore permit unbounded deformation under increasing stress
\cite{rivlin1948large,arruda1993three}. While such behavior may be acceptable within
limited loading regimes, it becomes problematic when constitutive laws are
extrapolated to large stresses, where physical admissibility and numerical
robustness are essential
\cite{gent1996new,horgan2004constitutive}.

From a constitutive perspective, strain-limiting behavior imposes two fundamental
requirements. First, the strain response must remain bounded for all admissible
stress states. Second, the incremental stiffness must progressively decrease as
the deformation limit is approached, ultimately leading to a vanishing tangent
modulus. These requirements are not modeling artifacts; rather, they reflect
underlying microstructural constraints such as molecular chain alignment, network
entanglement, and geometric locking that govern large-deformation behavior in
elastomers, polymers, and biological tissues
\cite{treloar1975physics,deam1976theory}.

Strain-limiting elasticity (SLE) models were introduced precisely to encode these
features directly at the constitutive level
\cite{rajagopal2011non,bulivcek2014elastic}. By construction, SLE formulations
enforce a finite strain bound together with a smooth transition from near-linear
elasticity at small stresses to saturation at large stresses. This structure
ensures physical admissibility across the entire stress domain and provides
improved mathematical well-posedness for boundary value problems involving extreme
loading
\cite{rajagopal2003implicit,bulivcek2015analysis}.

Despite these advantages, the strongly nonlinear and asymptotic nature of
strain-limiting constitutive laws introduces practical challenges. Analytical
manipulation becomes less tractable, and numerical implementation requires careful
treatment to avoid instability as the strain limit is approached
\cite{bulivcek2015existence,itou2017nonlinear}. In experimental settings, model
calibration further complicates this picture, as fidelity to measured responses
must be balanced against strict enforcement of bounded deformation
\cite{treloar1944stress,jones1975properties}.

These considerations motivate the development of modeling frameworks that preserve
the essential structure of strain-limiting elasticity while remaining amenable to
systematic calibration and data-driven refinement. In the following sections, we
formalize the strain-limiting constitutive law adopted in this work and introduce
an interpretable learning framework designed to augment the physics-based
structure without compromising its fundamental mechanical properties
\cite{liu2024kan,abdolazizi2025constitutive}.

\subsection{Mathematical Formulation}
\label{sec:mathematical_formulation}

We restrict attention to a one-dimensional stress--strain mapping, treated
mode-wise for each deformation path, as a phenomenological constitutive surrogate
rather than a full three-dimensional hyperelastic potential. This setting is
sufficient for the present synthetic benchmarks and experimental Treloar studies,
which are analyzed independently for each loading mode.

To adopt a single constitutive framework applicable to both the synthetic
benchmarks and the experimental data, we introduce a strain-limiting elasticity
(SLE) law with an explicit small-strain modulus
\cite{rajagopal2011non,bulivcek2014elastic}. In one dimension, the constitutive
relation is written as
\begin{equation}
\varepsilon(\tau)
=
\frac{\tau/E}{\left(1 + (\beta |\tau|)^{\alpha}\right)^{1/\alpha}},
\label{eq:strain_limiting}
\end{equation}
where $\tau$ denotes the applied stress, $\varepsilon(\tau)$ the corresponding
strain, $E>0$ the Young's modulus (small-strain tangent modulus), and
$\alpha>0$, $\beta>0$ are strain-limiting parameters.

In the small-stress regime, $(\beta|\tau|)^{\alpha}\ll 1$, and
Eq.~\eqref{eq:strain_limiting} reduces to linear elasticity,
\begin{equation}
\varepsilon(\tau) \approx \frac{\tau}{E},
\label{eq:linear_limit}
\end{equation}
so that $E$ governs the initial slope of the stress--strain response, consistent
with classical linearized elasticity
\cite{truesdell1952mechanical,marsden1994mathematical}. In the synthetic studies,
we work in normalized units and set $E=1$ without loss of generality, whereas in
the experimental setting $E$ is calibrated in physical units and directly reflects
the small-strain stiffness of the material.

As $|\tau|$ increases, the nonlinear term in
Eq.~\eqref{eq:strain_limiting} becomes dominant and the strain approaches a finite
bound,
\begin{equation}
\lim_{|\tau|\to\infty} |\varepsilon(\tau)|
=
\frac{1}{E\beta},
\label{eq:strain_limit}
\end{equation}
which is a defining feature of strain-limiting elasticity
\cite{rajagopal2003implicit,bulivcek2015analysis}. For notational convenience, we
introduce the compound parameter
\begin{equation}
\gamma = E\beta,
\label{eq:gamma_def}
\end{equation}
so that the strain limit can be expressed compactly as
$\varepsilon_{\max}=1/\gamma$.

The parameters $\alpha$ and $\beta$ play distinct and complementary roles in
shaping the constitutive response. The exponent $\alpha$ controls the sharpness
of the transition from the near-linear regime in
Eq.~\eqref{eq:linear_limit} to saturation, while $\beta$ (or equivalently
$\gamma$) determines the saturation level, in agreement with general
strain-limiting formulations
\cite{bulivcek2014elastic,rajagopal2014nonlinear}.

Equation~\eqref{eq:strain_limiting} is odd symmetric,
\begin{equation}
\varepsilon(-\tau) = -\varepsilon(\tau),
\label{eq:odd_symmetry}
\end{equation}
implying identical constitutive behavior in tension and compression, consistent
with isotropic elastic response
\cite{rivlin1948large}.

The tangent modulus,
\begin{equation}
E_t(\tau) = \frac{d\varepsilon}{d\tau},
\label{eq:tangent_modulus}
\end{equation}
decreases monotonically with increasing $|\tau|$ and vanishes in the saturation
limit,
\begin{equation}
\lim_{|\tau|\to\infty} E_t(\tau) = 0,
\label{eq:vanishing_tangent}
\end{equation}
a hallmark of strain-limiting constitutive behavior
\cite{rajagopal2011non,bulivcek2015analysis}.

Together, Eqs.~\eqref{eq:linear_limit}--\eqref{eq:vanishing_tangent} demonstrate
that the SLE law provides a smooth interpolation between linear elasticity at
small stresses and bounded deformation with vanishing incremental stiffness at
large stresses. The explicit separation of small-strain stiffness ($E$) from
strain saturation ($\beta$ or $\gamma$) enables a consistent transition from
normalized synthetic studies to experimental calibration in physical units, as
required for meaningful comparison with classical rubber elasticity experiments
\cite{treloar1944stress,treloar1975physics}.

\subsection{Interpretation of Model Parameters}

The strain-limiting constitutive law introduced in
Eq.~\eqref{eq:strain_limiting} involves a small number of parameters, each of
which admits a clear and physically meaningful interpretation
\cite{rajagopal2011non,bulivcek2014elastic}. For reference,
Table~\ref{tab:parameters} summarizes the role of these quantities within the
constitutive formulation.

\begin{table}[h!]
\centering
\caption{Interpretation of strain-limiting model parameters}
\label{tab:parameters}
\begin{tabular}{ll}
\toprule
\textbf{Parameter} & \textbf{Physical Interpretation} \\
\midrule
$\tau$ & Applied stress \\
$\varepsilon(\tau)$ & Resulting strain response \\
$E$ & Small-strain (Young's) modulus governing initial stiffness \\
$\alpha$ & Controls the sharpness of the transition to strain saturation \\
$\beta$ & Strain-limiting parameter \\
$\gamma = E\beta$ & Effective strain-limiting strength (see Sec.~\ref{sec:mathematical_formulation}) \\
\bottomrule
\end{tabular}
\end{table}

The Young's modulus $E$ governs the small-stress response of the material.
As shown by the linear limit in Eq.~\eqref{eq:linear_limit}, $E$ determines the
initial slope of the stress--strain curve and therefore represents the
small-strain stiffness, consistent with classical linearized elasticity
\cite{truesdell1952mechanical,marsden1994mathematical}. In the synthetic benchmarks
considered later, $E$ is set to unity through normalization, whereas in the
experimental studies $E$ is identified from data in physical units and directly
reflects the measured small-strain stiffness
\cite{treloar1944stress,treloar1975physics}.

The parameters $\alpha$ and $\beta$ jointly control the nonlinear,
strain-limiting behavior. Their combined influence is conveniently expressed
through the compound parameter $\gamma = E\beta$, which characterizes the
strength of strain limitation and governs the asymptotic deformation bound
introduced in Sec.~\ref{sec:mathematical_formulation}
\cite{rajagopal2003implicit,bulivcek2015analysis}. Smaller values of $\gamma$
correspond to more compliant responses with larger admissible strains, whereas
larger values enforce stronger strain limitation and earlier saturation, in
accordance with strain-limiting elastic theories
\cite{rajagopal2014nonlinear}.

The exponent $\alpha$ governs the manner in which the material transitions from
near-linear elasticity to the saturation regime. Smaller values of $\alpha$
produce a smooth and gradual transition characterized by a slowly decaying
tangent modulus, while larger values yield a sharper change in curvature and a
more abrupt approach to the strain limit
\cite{bulivcek2014elastic,bulivcek2015analysis}.

Taken together, the parameter set $(E,\alpha,\gamma)$ provides a compact and
physically interpretable description of strain-limiting elasticity. The explicit
separation of small-strain stiffness ($E$), saturation strength ($\gamma$), and
transition sharpness ($\alpha$) is particularly advantageous when integrating
physics-based constitutive modeling with data-driven refinement, as each
parameter retains a clear mechanical meaning across both synthetic and
experimental settings
\cite{rajagopal2011non,abdolazizi2025constitutive}.

\subsection{Key Physical Properties}

The strain-limiting constitutive law defined in
Eq.~\eqref{eq:strain_limiting} satisfies several fundamental physical
properties that are essential for constitutive admissibility and numerical
robustness. These properties arise directly from the mathematical structure
of the formulation and play a central role in both analytical modeling and
data-driven approximation
\cite{rajagopal2011non,bulivcek2014elastic,bulivcek2015analysis}.

\subsubsection{Odd Symmetry}

The constitutive response is odd symmetric with respect to the applied stress,
as expressed by Eq.~\eqref{eq:odd_symmetry}. This symmetry guarantees identical
mechanical behavior in tension and compression, with the strain changing sign
upon reversal of the applied stress. Odd symmetry is a fundamental requirement
for one-dimensional isotropic elastic materials and ensures consistency with
basic principles of continuum mechanics
\cite{truesdell1952mechanical,marsden1994mathematical}. 

From a modeling perspective, exact enforcement of odd symmetry eliminates the
need to treat tensile and compressive responses separately and provides a
natural structural constraint for learning-based constitutive
representations. Embedding this property at the model level reduces redundancy
and improves numerical behavior near the stress-free configuration
\cite{rajagopal2003implicit}.

\subsubsection{Bounded Strain}

A defining feature of the strain-limiting formulation is the boundedness of the
strain response. For all admissible stress values, the strain magnitude remains
uniformly bounded by the finite strain limit introduced in
Sec.~\ref{sec:mathematical_formulation}. This property ensures that deformation
remains finite even under arbitrarily large applied stresses, reflecting the
experimentally observed saturation of deformation in polymeric and soft
materials and preventing the nonphysical prediction of unbounded strain
\cite{rajagopal2011non,bulivcek2014elastic}.

From a numerical standpoint, the presence of a finite strain bound enhances
stability and robustness, particularly in simulations involving extreme loading
conditions or extrapolation beyond the calibration range. In such regimes,
bounded strain responses mitigate the risk of loss of ellipticity and
pathological mesh distortion
\cite{bulivcek2015existence,bulivcek2015analysis}.

\subsubsection{Vanishing Tangent Modulus}

The tangent modulus, defined in Eq.~\eqref{eq:tangent_modulus}, characterizes
the incremental stiffness of the material. Differentiation of
Eq.~\eqref{eq:strain_limiting} reveals that the tangent modulus decreases
monotonically with increasing stress magnitude and asymptotically vanishes, as
stated in Eq.~\eqref{eq:vanishing_tangent}
\cite{rajagopal2014nonlinear,bulivcek2014elastic}. As the strain limit is
approached, additional increases in stress produce progressively smaller
increments of strain, indicating a gradual loss of incremental deformability.

This vanishing stiffness behavior is a hallmark of strain-limiting elasticity
and is consistent with the physical interpretation of microstructural
constraints becoming dominant at large deformations
\cite{horgan2002constitutive,horgan2004constitutive}. From a modeling and
learning perspective, explicit representation of stiffness decay is critical
for ensuring physically meaningful asymptotic behavior and stable numerical
performance in high-stress regimes.

\subsubsection{Implications for Constitutive Modeling and Learning}

Taken together, odd symmetry, bounded strain, and vanishing tangent modulus define
a constitutively admissible response that is both physically realistic and
mathematically well behaved. These properties ensure stable material behavior
under extreme loading conditions and play a critical role in preventing
nonphysical predictions such as unbounded deformation or artificial stiffness
recovery at large stresses
\cite{rajagopal2011conspectus,bulivcek2014elastic}.

From the perspective of data-driven constitutive modeling, preservation of these
properties is particularly important. Machine learning models are frequently
trained on finite datasets spanning limited stress or strain ranges, yet are
often deployed in regimes that require extrapolation. In the absence of explicit
structural constraints, even highly accurate black-box models may violate basic
mechanical principles when evaluated outside the training domain, leading to
instabilities or loss of predictive reliability.

These considerations motivate the use of structured and interpretable learning
frameworks that encode essential physical properties directly at the model
level, rather than relying on their implicit emergence from data. Explicit
enforcement of symmetry, boundedness, and asymptotic stiffness decay provides a
robust foundation for learning-based constitutive representations and enables
controlled extrapolation consistent with mechanical admissibility
\cite{liu2024kan,abdolazizi2025constitutive}.

This viewpoint underpins the Kolmogorov--Arnold Network--based constitutive
formulation introduced in the following section, in which the architectural
structure of the learning model is aligned with the fundamental features of
strain-limiting elasticity. By embedding these properties into the representation
itself, the resulting framework balances expressivity with interpretability and
ensures physically consistent behavior across both data-rich and data-sparse
regimes.

\subsection{Behavior Across Stress Regimes}

The strain-limiting constitutive response defined in
Eq.~\eqref{eq:strain_limiting} exhibits qualitatively distinct behavior depending
on the magnitude of the applied stress. From a physical and modeling standpoint,
the stress domain may be partitioned into three regimes that reflect the
progressive activation of strain-limiting mechanisms
\cite{rajagopal2011non,bulivcek2014elastic}.

\begin{enumerate}
    \item \textbf{Low-stress regime.}  
    For sufficiently small stress magnitudes, the nonlinear contribution in the
    denominator of Eq.~\eqref{eq:strain_limiting} is negligible. The constitutive
    response therefore reduces to its linear elastic approximation given by
    Eq.~\eqref{eq:linear_limit}. In this regime, the stress--strain relationship
    is nearly linear, and the tangent modulus remains approximately constant,
    governed by the small-strain modulus $E$
    \cite{truesdell1952mechanical,marsden1994mathematical}.

    \item \textbf{Transition regime.}  
    As the applied stress increases, nonlinear effects become progressively
    significant. The stress--strain response departs smoothly from linearity,
    and the tangent modulus defined in Eq.~\eqref{eq:tangent_modulus} begins to
    decrease. This regime marks the gradual activation of strain-limiting
    mechanisms, in which additional stress increments produce diminishing strain
    increments while deformation continues to evolve
    \cite{bulivcek2014elastic,bulivcek2015analysis}.

    \item \textbf{High-stress regime.}  
    At sufficiently large stress magnitudes, the strain asymptotically approaches
    the finite bound introduced in Sec.~\ref{sec:mathematical_formulation}. In
    this regime, the tangent modulus vanishes asymptotically, as stated in
    Eq.~\eqref{eq:vanishing_tangent}, indicating that further increases in stress
    result in negligible additional deformation. This behavior reflects the
    saturation of strain associated with microstructural constraints and is a
    defining characteristic of strain-limiting elasticity
    \cite{rajagopal2014nonlinear,horgan2004constitutive}.
\end{enumerate}

This regime-based interpretation provides a concise and physically transparent
description of the constitutive response across the full stress domain. It also
highlights a central challenge for data-driven constitutive modeling: the need to
accurately represent near-linear behavior at small stresses, smooth nonlinear
transition behavior, and asymptotic strain saturation within a single unified
model. Capturing all three regimes simultaneously, while preserving physical
admissibility and numerical stability, motivates the use of structured learning
frameworks aligned with strain-limiting mechanics
\cite{liu2024kan,abdolazizi2025constitutive}.

\subsection{Algorithmic Evaluation of the Constitutive Law}

\begin{algorithm}[h!]
\caption{Algorithmic Evaluation of the Strain-Limiting Constitutive Law}
\label{alg:strain_eval}
\begin{algorithmic}[1]
\Require Applied stress $\tau$, material parameters $\alpha$, $\beta$, Young's modulus $E$
\Ensure Strain $\varepsilon$
\State Compute stress magnitude: $s \gets |\tau|$
\State Evaluate nonlinear denominator:
$d \gets \left(1 + (\beta s)^{\alpha}\right)^{1/\alpha}$
\State Compute strain magnitude:
$\varepsilon_s \gets \dfrac{s}{E\, d}$
\State Recover signed strain:
$\varepsilon \gets \mathrm{sign}(\tau)\,\varepsilon_s$
\State \Return $\varepsilon$
\end{algorithmic}
\end{algorithm}

Algorithm~\ref{alg:strain_eval} provides an explicit and transparent procedure for
evaluating the strain response associated with a prescribed applied stress using
the strain-limiting constitutive law defined in
Eq.~\eqref{eq:strain_limiting}
\cite{rajagopal2011non,bulivcek2014elastic}. Importantly, the algorithm is not
merely a numerical implementation; it reveals an intrinsic structural
decomposition of the constitutive mapping that underlies both its physical
interpretation and its suitability for learning-based augmentation.

\paragraph{Stress magnitude extraction.}
The evaluation begins by isolating the stress magnitude,
$s = |\tau|$, exploiting the odd symmetry of the constitutive response. This
operation reduces the constitutive mapping to a nonlinear transformation defined
on a nonnegative scalar domain, while the sign of the stress is treated
separately. Such a decomposition is consistent with the symmetry requirements of
one-dimensional isotropic elasticity and simplifies both analysis and numerical
treatment
\cite{truesdell1952mechanical,rajagopal2003implicit}.

\paragraph{Nonlinear magnitude transformation.}
The nonlinear denominator governs the evolution from linear elasticity to strain
saturation. For small stress magnitudes, the denominator remains close to unity,
recovering an approximately linear response. As $s$ increases, the denominator
grows monotonically, progressively attenuating the incremental strain generated
by additional stress. This mechanism encodes the gradual activation of
strain-limiting effects and ensures a smooth transition toward deformation
saturation
\cite{bulivcek2014elastic,rajagopal2014nonlinear}.

\paragraph{Bounded strain computation.}
The strain magnitude is computed as $\varepsilon_s = s/(E d)$, guaranteeing that
the strain remains uniformly bounded for all admissible stress values. This step
enforces the finite strain limit intrinsic to strain-limiting elasticity and
ensures vanishing incremental stiffness at large stresses
\cite{bulivcek2015analysis}. In the synthetic studies, the Young’s modulus is set
to $E=1$ through normalization, whereas in the experimental studies $E$ is
calibrated in physical units and retained explicitly
\cite{treloar1944stress,treloar1975physics}.

\paragraph{Sign-preserving reconstruction.}
Finally, the signed strain response is recovered by reintroducing the sign of the
applied stress. This deterministic reconstruction enforces odd symmetry exactly
and guarantees identical constitutive behavior in tension and compression, without
relying on approximate or data-driven symmetry enforcement
\cite{truesdell1952mechanical,marsden1994mathematical}.

\medskip
Taken together, Algorithm~\ref{alg:strain_eval} exposes the constitutive mapping as
a composition of two interpretable stages: (i) a nonlinear transformation acting
solely on the stress magnitude, followed by (ii) an exact sign-preserving
reconstruction. This separable structure is intrinsic to the strain-limiting
constitutive law and is independent of any learning strategy. Crucially, it
provides the structural blueprint for the Kolmogorov--Arnold Network formulation
introduced in the next section, where each stage of the algorithm is mirrored by a
corresponding, physically interpretable network component
\cite{liu2024kan,abdolazizi2025constitutive}.

\subsection{Visualization of Strain-Limiting Behavior}

\begin{figure}[h!]
\centering
\includegraphics[width=0.65\linewidth]{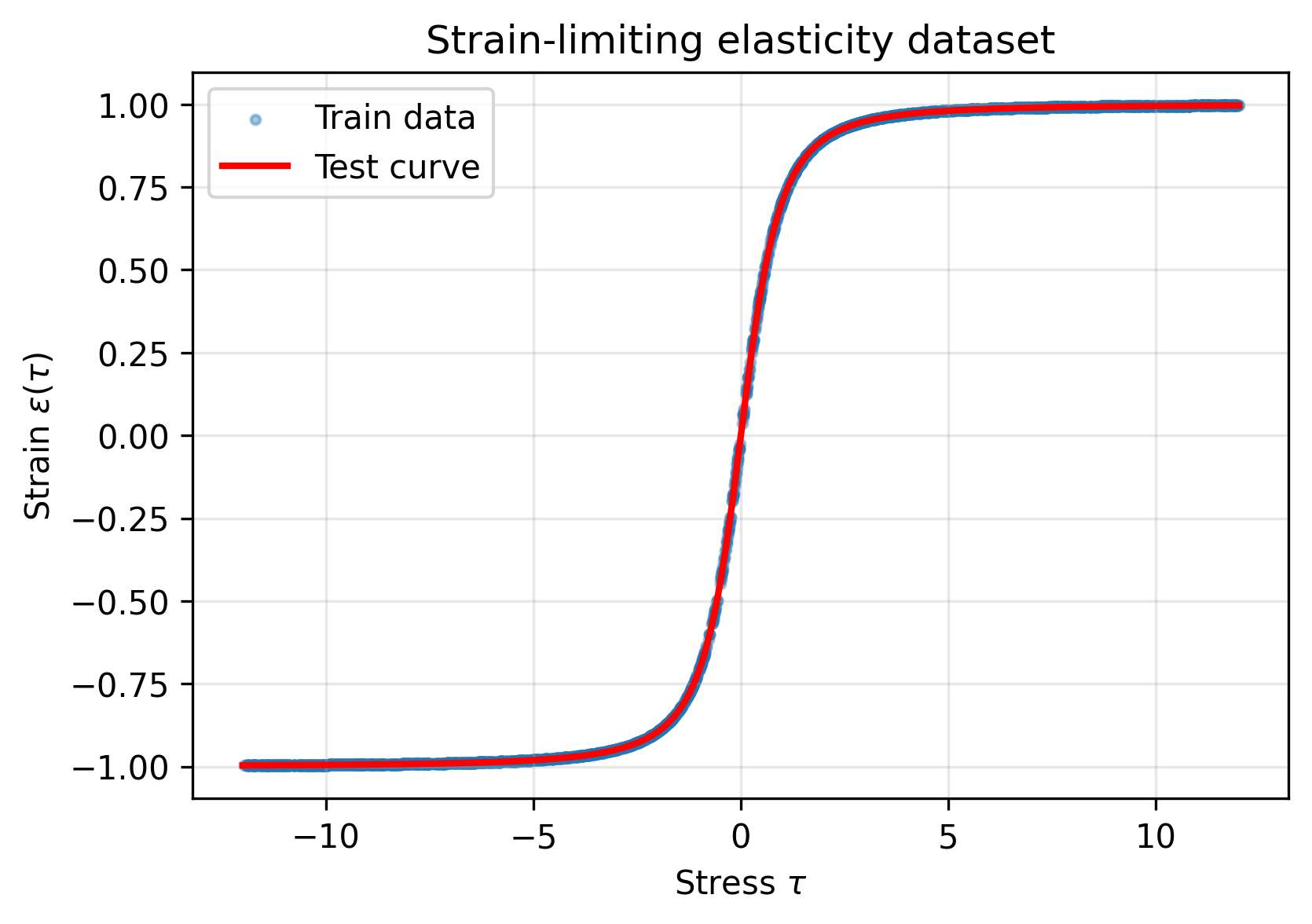}
\caption{Representative stress--strain response exhibiting strain-limiting
behavior. The response is approximately linear at small stresses and
asymptotically approaches a finite strain bound as the applied stress increases.}
\label{fig:stress_strain_example}
\end{figure}

To provide an intuitive and diagnostic illustration of the characteristic
features of strain-limiting elasticity, Fig.~\ref{fig:stress_strain_example}
shows a representative stress--strain response generated using the constitutive
law in Eq.~\eqref{eq:strain_limiting}
\cite{rajagopal2011non,bulivcek2014elastic}. The figure encapsulates, in a single
curve, the essential mechanical signatures of strain-limiting behavior and
highlights the continuous transition between distinct deformation regimes.

At small stress magnitudes, the constitutive response is nearly linear, in
agreement with the classical elastic limit described by
Eq.~\eqref{eq:linear_limit} and established theories of linearized elasticity
\cite{truesdell1952mechanical,marsden1994mathematical}. In this regime, the
tangent modulus remains approximately constant and equal to the small-strain
modulus $E$, reflecting conventional elastic behavior. As the applied stress increases, nonlinear effects become progressively
dominant. The slope of the stress--strain curve decreases monotonically,
indicating a continuous reduction of the tangent modulus. This gradual stiffness
degradation reflects the activation of microstructural constraints—such as chain
alignment and network locking—and constitutes a defining feature of
strain-limiting elastic models
\cite{rajagopal2014nonlinear,bulivcek2015analysis}.

In the high-stress regime, the strain asymptotically approaches a finite upper
bound, demonstrating deformation saturation in accordance with
Eq.~\eqref{eq:strain_limit}
\cite{bulivcek2014elastic,horgan2004constitutive}. In this limit, further
increases in stress produce negligible incremental strain, corresponding to a
vanishing tangent modulus. Importantly, this saturation is achieved smoothly,
without abrupt changes in curvature or loss of regularity. This visualization reinforces the regime-based interpretation discussed earlier
and highlights the central challenge for data-driven constitutive modeling:
capturing both the near-linear response at small stresses and the asymptotic
flattening associated with bounded strain within a single unified
representation. The smooth yet highly nonlinear transition evident in
Fig.~\ref{fig:stress_strain_example} underscores why unconstrained black-box
approximations often struggle to extrapolate reliably, and motivates the use of
structured, physics-aligned learning frameworks that preserve boundedness,
monotonicity, and asymptotic behavior by construction
\cite{rajagopal2011conspectus,bulivcek2014elastic,liu2024kan}.

\subsection{Relevance for Data-Driven Modeling}

Although the strain-limiting constitutive law in
Eq.~\eqref{eq:strain_limiting} is physically well founded, its strongly nonlinear
and asymptotic structure presents substantial challenges for data-driven
approximation
\cite{rajagopal2011non,bulivcek2014elastic}. Any reliable surrogate model must
simultaneously reproduce the near-linear elastic response at small stresses and
the progressive saturation of strain at large stresses, while strictly preserving
boundedness, odd symmetry, and monotonic degradation of incremental stiffness
across the entire stress domain.

Conventional black-box learning models, including fully connected neural
networks, can often achieve low training error over restricted stress ranges.
However, their unconstrained expressivity frequently leads to poor extrapolation
behavior. In the context of strain-limiting elasticity, such models may violate
essential mechanical requirements—predicting unbounded strain, loss of symmetry,
or spurious stiffening or softening at large stresses—resulting in nonphysical
responses outside the calibration regime
\cite{shen2004neural,lee2025constitutive}. These deficiencies are particularly
problematic in solid mechanics, where extrapolation beyond available experimental
data is often unavoidable.

These limitations motivate the use of structured learning approaches that embed
constitutive principles directly into the model representation rather than relying
on implicit inference from data alone. By enforcing physical admissibility at the
architectural level, such frameworks promote robustness, interpretability, and
predictive reliability, especially in regimes where experimental observations are
sparse or unavailable
\cite{rajagopal2003implicit,abdolazizi2024viscoelastic}.

Kolmogorov--Arnold Networks (KANs) provide a natural foundation for this purpose by
representing nonlinear mappings as structured compositions of low-dimensional,
interpretable functions
\cite{liu2024kan,ji2024comprehensive}. This architecture aligns closely with the
separable structure of strain-limiting constitutive laws, enabling symmetry
preservation, bounded response, and asymptotic behavior to be enforced by
construction rather than approximated through penalty-based regularization. As a
result, KANs offer a principled and transparent framework for learning
strain-limiting elasticity models that maintain physical consistency while
exhibiting robust extrapolation behavior beyond the training domain
\cite{abdolazizi2025constitutive,panahi2025data}.

\section{KAN-Based Constitutive Representation}

The constitutive decomposition introduced in
Eq.~\eqref{eq:kan_decomposition} is directly motivated by the intrinsic structure
of strain-limiting elasticity (SLE) models
\cite{rajagopal2011non,bulivcek2014elastic}. From a mechanical standpoint, the
stress--strain response exhibits two fundamental characteristics: \emph{odd
symmetry}, corresponding to identical material behavior in tension and
compression, and \emph{strain saturation}, whereby deformation depends primarily
on the magnitude of the applied stress rather than its sign.

These properties naturally motivate a sign--magnitude decomposition of the
constitutive response,
\begin{equation}
\varepsilon(\tau) = \mathrm{sign}(\tau)\, g(|\tau|),
\label{eq:kan_decomposition}
\end{equation}
where $g(|\tau|)$ is a nonnegative scalar function mapping the stress magnitude to
the corresponding strain magnitude. Such decompositions are standard in isotropic
elasticity and implicit constitutive theories, where symmetry and monotonicity
are enforced structurally rather than inferred empirically
\cite{truesdell1952mechanical,rajagopal2003implicit}.

By isolating the sign of the applied stress, the learning task is reduced to the
approximation of a scalar nonlinear function defined on a nonnegative domain.
Consequently, odd symmetry is enforced exactly by construction, ensuring
constitutive admissibility for all stress levels, including regimes well beyond
the training data
\cite{rajagopal2011conspectus}. This structural enforcement distinguishes the
proposed formulation from unconstrained black-box regressions, where symmetry
must be learned implicitly and may be violated under extrapolation.

From a numerical perspective, restricting learning to the magnitude domain
$|\tau|\ge0$ significantly improves conditioning during optimization. In
conventional neural networks, the sign change at $\tau=0$ often introduces sharp
gradients and spurious oscillations near the origin, degrading convergence and
stability. By treating the sign deterministically and learning only the smooth
magnitude response, the present formulation yields a better-conditioned
optimization problem and more reliable training behavior
\cite{shen2004neural,lee2025constitutive}.

The function $g(|\tau|)$ admits a direct physical interpretation as the
\emph{strain magnitude envelope} associated with increasing stress. Its
near-linear behavior at small stress magnitudes reflects compliant elastic
response, while its asymptotic saturation captures the strain-limiting behavior
observed experimentally in rubber-like and soft materials
\cite{treloar1944stress,treloar1975physics,horgan2004constitutive}. When
$g(|\tau|)$ is represented using bounded and smooth basis functions, the resulting
constitutive law automatically respects the prescribed strain limit without
requiring penalty terms or post-processing corrections
\cite{bulivcek2015analysis}. This decomposition aligns naturally with the Kolmogorov--Arnold representation
theorem, which expresses nonlinear mappings as compositions of low-dimensional,
interpretable functions
\cite{liu2024kan,ji2024comprehensive}. In the present one-dimensional constitutive
setting, the KAN architecture reduces to a single nonlinear mapping that
represents $g(|\tau|)$, followed by a deterministic sign reconstruction. Each
component of the network therefore corresponds directly to a physically
meaningful operation, enabling transparent interpretation and validation of the
learned constitutive behavior
\cite{abdolazizi2025constitutive}.

Overall, the KAN-based constitutive representation provides a principled balance
between modeling flexibility and mechanical consistency. By embedding symmetry,
boundedness, and smoothness directly into the model architecture, the proposed
approach enables accurate approximation of strain-limiting behavior while
maintaining interpretability, numerical robustness, and reliable extrapolation
under large stresses
\cite{panahi2025data}.

\subsection{Conceptual Overview of the KAN Architecture}

\begin{figure}[h!]
\centering
\includegraphics[width=0.65\linewidth]{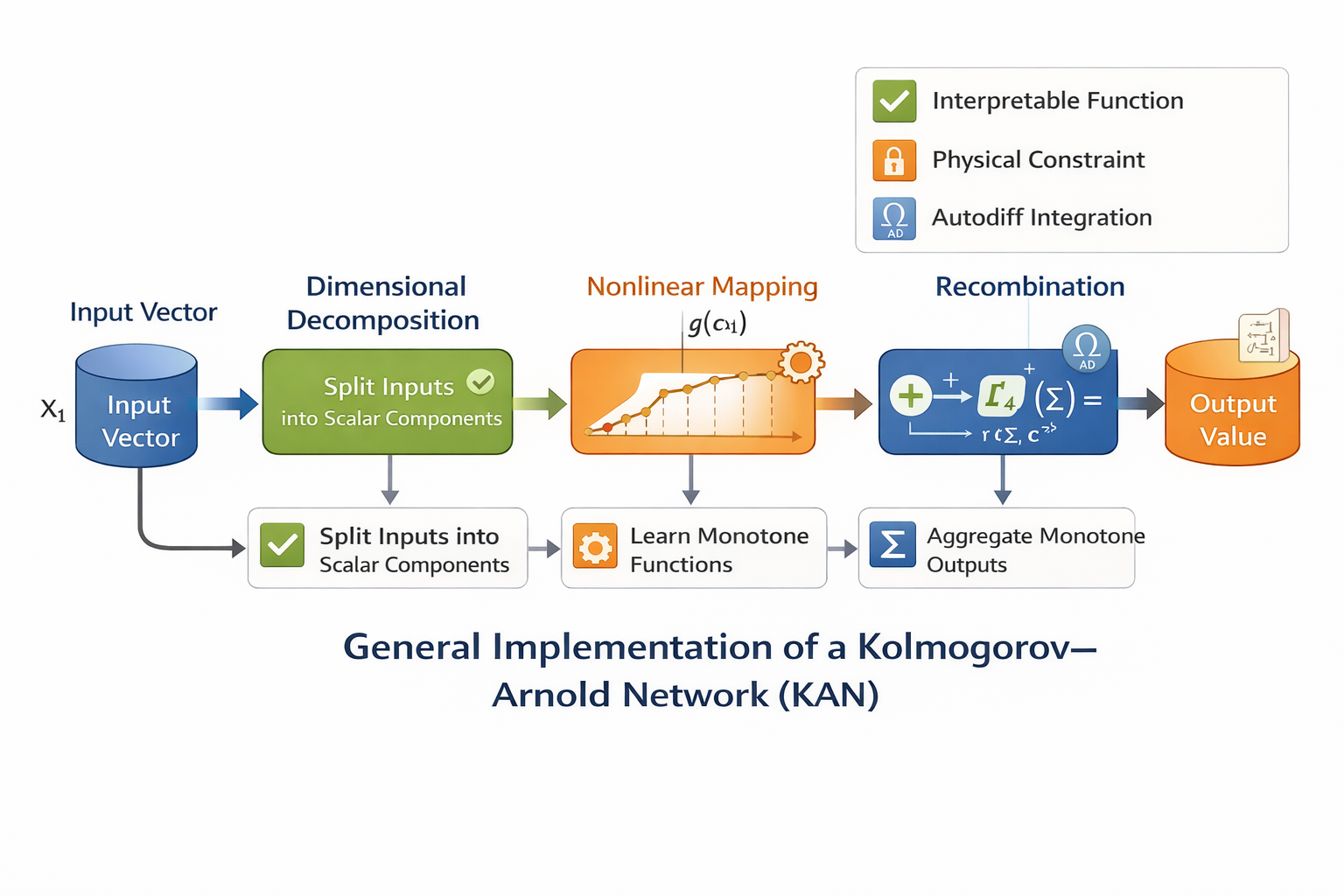}
\caption{General structure of a Kolmogorov--Arnold Network (KAN). The input is
decomposed into scalar components, transformed through interpretable nonlinear
functions, and recombined to form the output.}
\label{fig:kan_overview}
\end{figure}

Figure~\ref{fig:kan_overview} illustrates the general structure of a
Kolmogorov--Arnold Network (KAN). The defining principle of the KAN architecture
is the representation of nonlinear mappings as structured compositions of
low-dimensional, interpretable functions, rather than as densely connected
black-box transformations
\cite{liu2024kan,ji2024comprehensive}. This design is rooted in the
Kolmogorov--Arnold representation theorem and is particularly well suited for
physics-informed modeling, where transparency, identifiability, and consistency
with underlying physical principles are essential.

As depicted in Fig.~\ref{fig:kan_overview}, the KAN architecture can be viewed as
a structured three-stage mapping. In the first stage, the input is decomposed
explicitly into its scalar components. This decomposition isolates the influence
of individual input variables and avoids the entangled feature representations
commonly encountered in fully connected neural networks
\cite{carneros2024comparison}.

In the second stage, each scalar component is processed independently through a
univariate nonlinear transformation. These transformations are represented
using structured function approximators, such as spline-based mappings. Unlike
standard neural network activation functions, these univariate mappings are
explicitly parameterized and directly interpretable. Their shape, slope, and
saturation characteristics can be examined and assessed against known physical
behavior, providing insight into how nonlinear response mechanisms are encoded
within the model
\cite{liu2024kan,essahraui2025kolmogorov}.

In the final stage, the transformed scalar outputs are recombined through a
simple aggregation operation to produce the network output. This aggregation is
typically linear or weakly nonlinear, ensuring that representational complexity
is concentrated within the interpretable univariate functions rather than
distributed across opaque dense layers. Automatic differentiation can be
applied consistently across all stages, enabling efficient gradient-based
training while preserving architectural transparency
\cite{panahi2025data}.

In the present one-dimensional constitutive modeling context, this general KAN
framework simplifies naturally. The applied stress $\tau$ is decomposed into its
magnitude $|\tau|$ and sign. The nonlinear mapping $g(|\tau|)$ governing the
strain magnitude is learned explicitly using a spline-based representation,
while the sign of the stress is reintroduced deterministically to recover the
signed strain response according to Eq.~\eqref{eq:kan_decomposition}. This
construction is fully consistent with the structure of strain-limiting
elasticity models
\cite{rajagopal2011non,bulivcek2014elastic}.

This specialization highlights a central advantage of the KAN approach for
constitutive modeling: the network architecture mirrors the intrinsic structure
of the underlying constitutive law. Each component of the model corresponds to
a physically meaningful operation, ensuring exact symmetry preservation,
improved numerical behavior near the origin, and enhanced interpretability of
the learned response. As a result, the KAN-based architecture provides a
transparent and physically aligned alternative to conventional black-box neural
networks for constitutive modeling
\cite{abdolazizi2025constitutive,lee2025constitutive}.

\subsection{Spline-Based Representation of the Constitutive Function}

Following the architectural principles outlined in the previous subsection,
the nonlinear function $g(|\tau|)$ appearing in the KAN-based decomposition
in Eq.~\eqref{eq:kan_decomposition} is represented using a piecewise-linear
spline defined over a fixed grid of stress magnitudes
\cite{liu2024kan,essahraui2025kolmogorov}. Let
$\{\tau_i\}_{i=1}^{N}$ denote a set of uniformly spaced knot locations spanning
the interval $[0,\tau_{\max}]$, where $\tau_{\max}$ denotes the maximum stress
magnitude encountered during training. The constitutive function is approximated as
\begin{equation}
g(|\tau|) \approx \sum_{i=1}^{N} c_i \, \phi_i(|\tau|),
\label{eq:spline_representation}
\end{equation}
where $\phi_i(|\tau|)$ are piecewise-linear basis functions associated with the
knot locations and $\{c_i\}$ are trainable spline coefficients.

\begin{figure}[h!]
\centering
\includegraphics[width=0.65\linewidth]{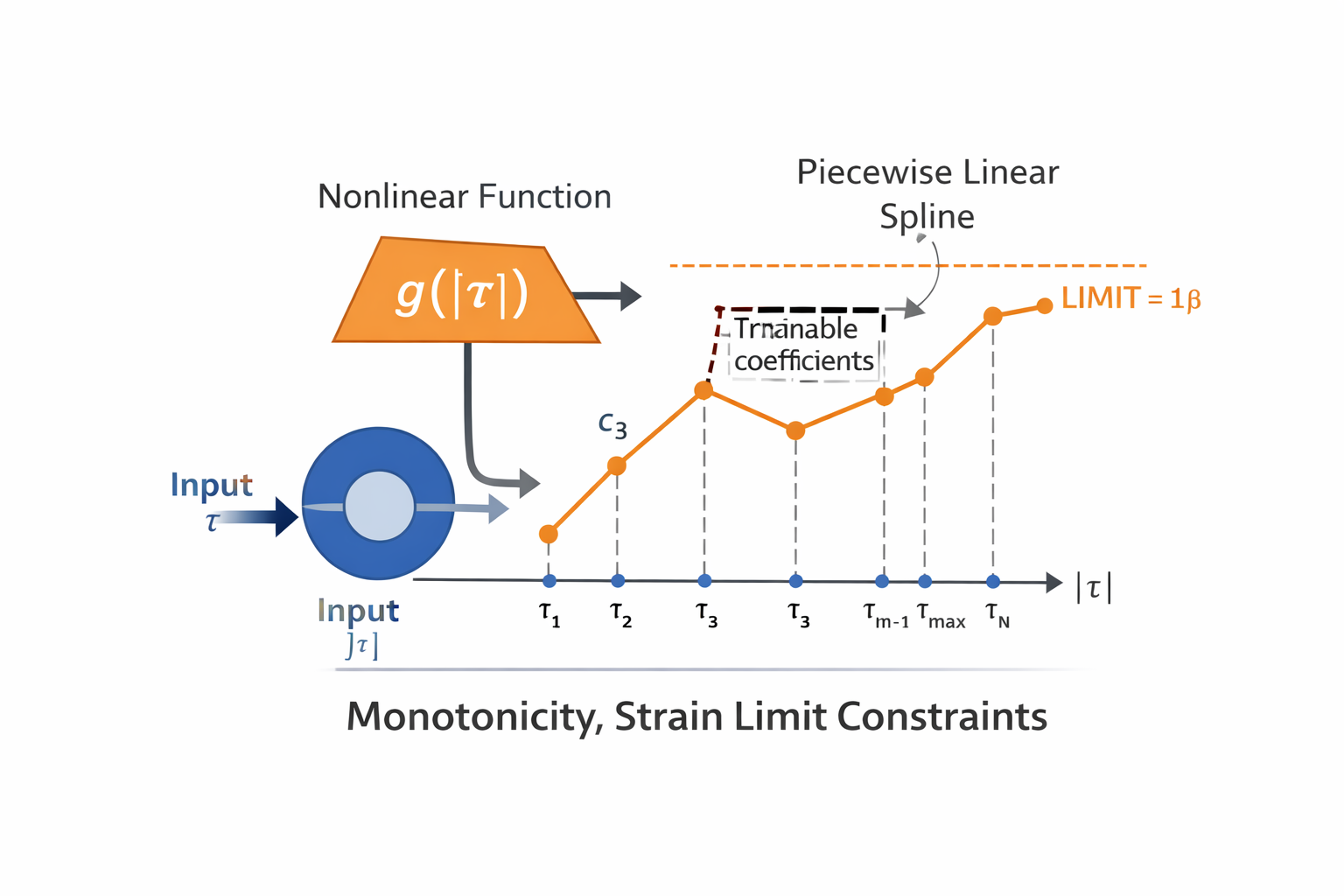}
\caption{Spline-based representation of the nonlinear constitutive function
$g(|\tau|)$. Trainable coefficients control the strain magnitude at discrete
stress knots, while the piecewise-linear structure ensures continuity and
controlled local behavior.}
\label{fig:spline_representation}
\end{figure}

Figure~\ref{fig:spline_representation} illustrates this spline-based
representation. Each coefficient $c_i$ specifies the strain magnitude attained
at the corresponding stress knot $\tau_i$, while linear interpolation between
adjacent knots defines the response at intermediate stress values. This
construction ensures continuity of the learned constitutive function and
enforces locally linear behavior within each stress interval, consistent with
structured KAN representations
\cite{liu2024kan,ji2024comprehensive}.

This representation offers several advantages for constitutive modeling.
First, the spline coefficients admit a direct and physically meaningful
interpretation: each parameter corresponds to the strain magnitude at a
prescribed stress level. In contrast to conventional neural network weights,
which are distributed across multiple layers and lack explicit mechanical
meaning, the spline parameters can be directly inspected, visualized, and
interpreted in terms of material response
\cite{abdolazizi2025constitutive,carneros2024comparison}.

Second, the derivative of the spline function is piecewise constant within each
interval between adjacent knots. Consequently, the local slope of the spline
provides an explicit representation of the tangent modulus over each stress
interval. This feature enables direct analysis of stiffness evolution as stress
increases and aligns naturally with strain-limiting elasticity, where progressive
stiffness degradation and asymptotic saturation are defining physical
characteristics
\cite{rajagopal2014nonlinear,bulivcek2015analysis}.

Finally, the spline-based formulation enables essential physical constraints to
be imposed directly at the parameter level. Monotonicity of the constitutive
response is enforced by constraining spline slopes to remain nonnegative, while
strain saturation is ensured by bounding the spline coefficients according to
the theoretical strain limit derived in
Sec.~\ref{sec:mathematical_formulation}
\cite{rajagopal2011non,bulivcek2014elastic}. Because these constraints act on
interpretable parameters, physical admissibility is preserved by construction,
without reliance on ad hoc penalty terms or post hoc correction strategies.

Overall, the spline-based representation provides a transparent, physically
interpretable, and numerically stable mechanism for approximating the nonlinear
constitutive function $g(|\tau|)$. When embedded within the KAN framework, this
approach yields a constitutive representation that balances expressive power with
explicit enforcement of mechanical structure, making it particularly well suited
for data-driven modeling of strain-limiting elasticity
\cite{panahi2025data}.

\subsection{Enforcement of Physical Constraints}

A central advantage of the spline-based KAN formulation is its ability to enforce
essential physical constraints directly at the level of model parameters, rather
than relying on implicit inference from data or auxiliary penalty-based
regularization
\cite{liu2024kan,abdolazizi2025constitutive}. By embedding admissibility conditions
into the model architecture, the learned constitutive response remains physically
consistent across the entire stress domain, including regimes not explicitly
sampled during training.

Monotonicity of the constitutive response is enforced by constraining the slopes
of the piecewise-linear spline to be nonnegative. Because each spline segment is
defined by adjacent coefficients, this constraint guarantees that the function
$g(|\tau|)$ is non-decreasing with respect to stress magnitude. Physically, this
reflects the fundamental requirement that strain magnitude must not decrease as
applied stress increases, thereby preserving mechanical consistency and
suppressing nonphysical oscillations in the learned response
\cite{rajagopal2011non,bulivcek2014elastic}.

Strain saturation is enforced by bounding the spline coefficients according to
the theoretical strain limit prescribed by the strain-limiting elasticity model.
By restricting the maximum admissible values of the spline coefficients, the
learned constitutive function is guaranteed to satisfy the bounded strain
condition for all stress levels, including those beyond the range of the training
data
\cite{bulivcek2015analysis,rajagopal2014nonlinear}. This architectural constraint
precludes unbounded extrapolation, a well-documented failure mode of unconstrained
black-box neural networks in constitutive modeling
\cite{shen2004neural,lee2025constitutive}.

\begin{figure}[h!]
\centering
\includegraphics[width=0.65\linewidth]{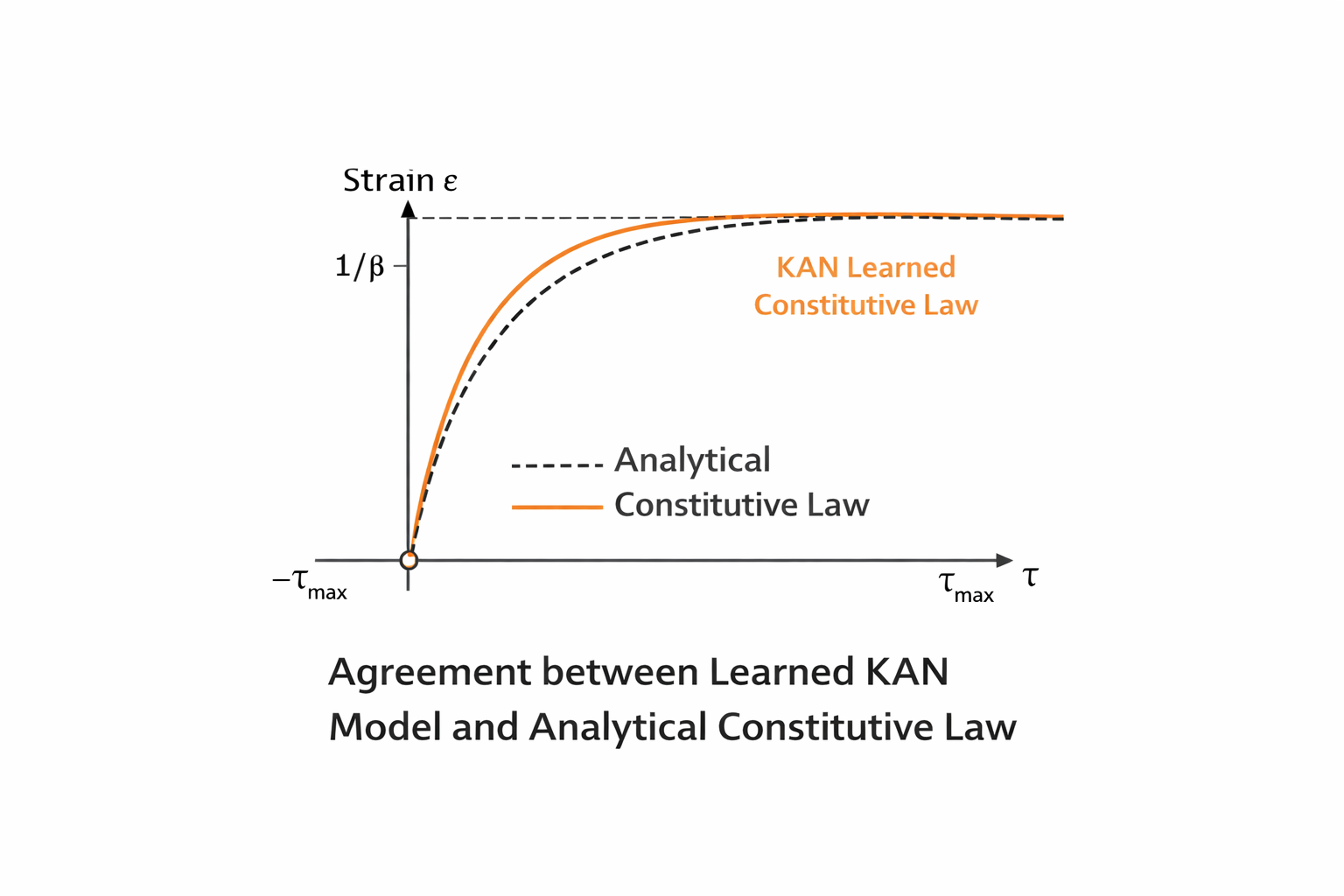}
\caption{Comparison between the learned KAN constitutive response and the
analytical strain-limiting model. The KAN reproduces the linear regime, transition
region, and asymptotic saturation behavior while remaining consistent with the
imposed physical constraints.}
\label{fig:strain_limit}
\end{figure}

Figure~\ref{fig:strain_limit} illustrates the resulting agreement between the
learned KAN constitutive response and the analytical strain-limiting law. The KAN
accurately reproduces the near-linear response at small stresses, captures the
smooth transition regime, and recovers the correct asymptotic saturation behavior
at large stresses. Importantly, saturation emerges directly from the constrained
spline representation, without the need for post-processing, explicit clipping,
or corrective enforcement
\cite{liu2024kan,panahi2025data}.

In contrast to conventional neural networks, which typically require extensive
data augmentation or penalty-based regularization to approximate properties such
as monotonicity and boundedness, the proposed KAN formulation embeds these
features directly within its representation. As a result, the learned model
exhibits improved numerical stability, particularly near the origin where the
stress changes sign, and demonstrates reliable extrapolation behavior under large
stresses
\cite{abdolazizi2024viscoelastic,carneros2024comparison}.

Overall, the explicit incorporation of physical constraints within the
spline-based KAN framework enables accurate approximation of strain-limiting
constitutive behavior while preserving physical realism, interpretability, and
numerical robustness. This constraint-aware design is central to the suitability
of the proposed approach for both synthetic benchmarks and experimental
constitutive modeling
\cite{rajagopal2011conspectus}.

\subsection{Comparison with Classical Neural Networks}

\begin{figure}[h!]
\centering
\includegraphics[width=0.65\linewidth]{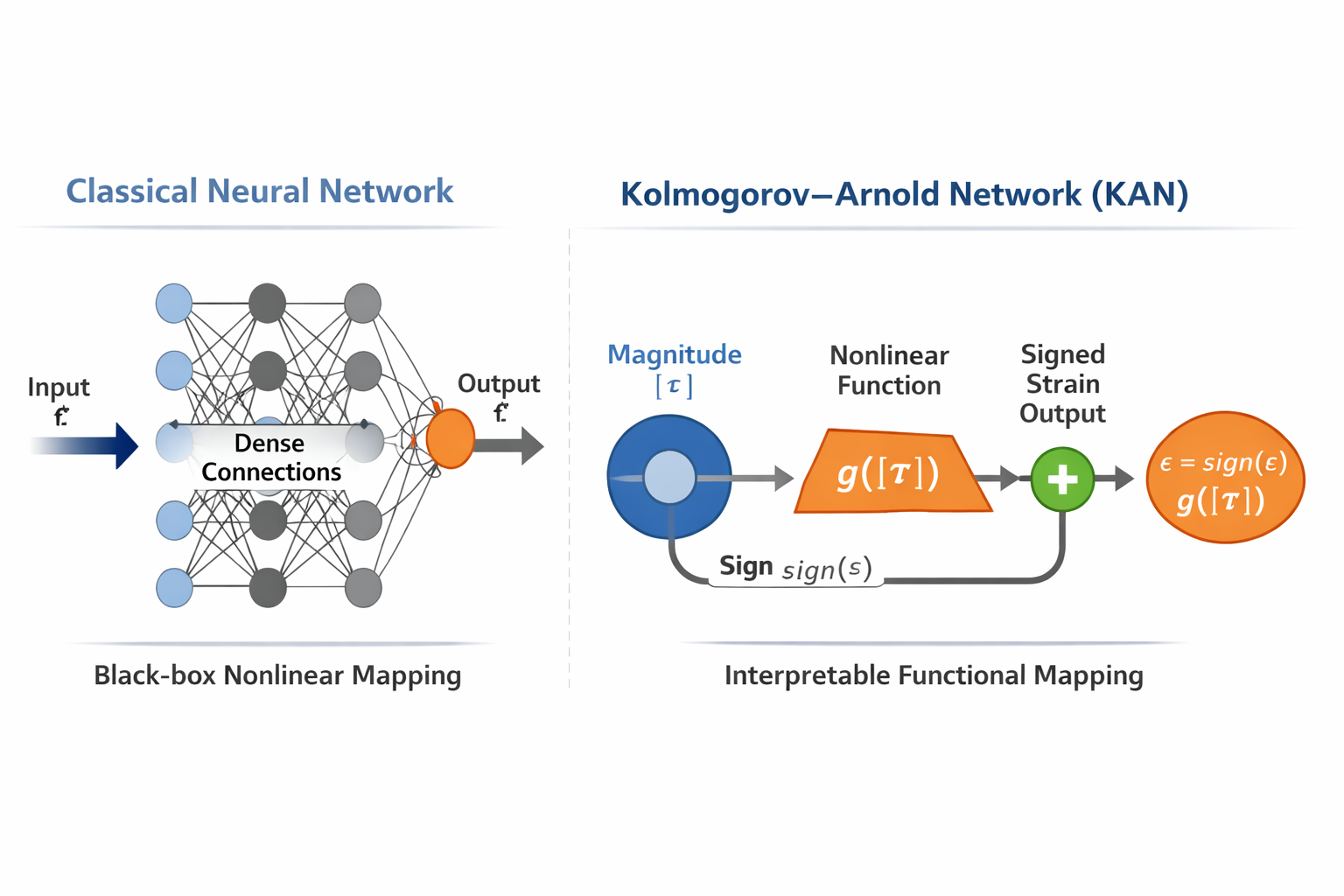}
\caption{Conceptual comparison between classical neural networks and
Kolmogorov--Arnold Networks. Classical networks rely on distributed nonlinear
representations, whereas KANs employ structured and interpretable functional
decompositions aligned with constitutive principles.}
\label{fig:cnn_vs_kan}
\end{figure}

Figure~\ref{fig:cnn_vs_kan} provides a conceptual comparison between classical
fully connected neural networks and the proposed Kolmogorov--Arnold Network
(KAN) formulation for constitutive modeling.

Classical neural networks approximate stress--strain relationships using multiple
layers of densely connected neurons combined with nonlinear activation functions.
Such architectures are highly expressive and have been shown to fit complex
constitutive datasets with high accuracy
\cite{shen2004neural,lee2025constitutive}. However, the resulting representation
is distributed across a large number of parameters, making it difficult to
associate individual network components with identifiable physical mechanisms.
Consequently, while predictive performance may be high, physical interpretability
and diagnostic insight are often limited.

Moreover, essential mechanical properties such as symmetry, monotonicity, and
bounded deformation are not enforced explicitly in conventional neural networks.
Instead, these properties must be inferred implicitly from data, often requiring
extensive data augmentation, carefully tuned regularization strategies, or large
training datasets. Even under such conditions, there is no guarantee that physical
admissibility will be preserved outside the training domain, which is a
significant limitation in constitutive modeling where extrapolation is frequently
unavoidable
\cite{rajagopal2003implicit,abdolazizi2024viscoelastic}.

In contrast, the KAN-based formulation embeds the structure of the
strain-limiting constitutive law directly into the network architecture. By
decomposing the stress--strain mapping into a magnitude-dependent nonlinear
function and a deterministic sign-preserving reconstruction, the KAN mirrors the
intrinsic physics of strain-limiting elasticity. As a result, key properties such
as odd symmetry, monotonicity, and bounded strain response are enforced at the
representation level, rather than being approximated through data-driven
optimization
\cite{rajagopal2011non,bulivcek2014elastic}.

A further distinction lies in interpretability. The spline-based functions used
within the KAN admit parameters with direct physical meaning: spline coefficients
correspond to strain magnitudes at prescribed stress levels, while their local
slopes represent tangent moduli over specific stress intervals. This structure
enables direct inspection and validation of the learned constitutive response
against analytical models or experimental observations, an analysis that is
generally not feasible with conventional densely connected neural networks
\cite{liu2024kan,abdolazizi2025constitutive}.

Overall, while classical neural networks emphasize representational flexibility,
the KAN formulation provides a structured alternative that balances expressive
power with physical consistency. By aligning the network architecture with
constitutive principles, the KAN-based approach promotes improved numerical
stability, enhanced interpretability, and more reliable extrapolation behavior,
making it particularly well suited for data-driven modeling of strain-limiting
elasticity
\cite{carneros2024comparison,panahi2025data}.

\subsection{Relation to Kolmogorov--Arnold Networks}

Kolmogorov--Arnold Networks (KANs) are motivated by the
Kolmogorov--Arnold representation theorem, which establishes that any continuous
multivariate function can be expressed as a finite composition of univariate
nonlinear functions and linear combinations
\cite{liu2024kan,ji2024comprehensive}. In contrast to classical neural networks,
which distribute nonlinear effects across densely connected layers, KANs
localize nonlinearity within explicit, low-dimensional functional components.
This architectural organization promotes transparency of the learned mapping and
facilitates the incorporation of problem-specific structural constraints
\cite{essahraui2025kolmogorov}.

In the present one-dimensional constitutive modeling setting, the general KAN
framework simplifies naturally. Because the stress--strain relationship depends
on a single scalar input, the multivariate Kolmogorov--Arnold decomposition
reduces to a single nonlinear univariate mapping followed by a deterministic
reconstruction step. The spline-based function $g(|\tau|)$ captures the nonlinear
dependence of strain magnitude on stress magnitude, while the sign-preserving
operation restores the correct deformation direction, as defined in
Eq.~\eqref{eq:kan_decomposition}. This structure is fully consistent with the
intrinsic form of strain-limiting elasticity
\cite{rajagopal2011non,bulivcek2014elastic}.

Importantly, this architectural decomposition mirrors the algorithmic evaluation
procedure of the strain-limiting constitutive law described in
Algorithm~\ref{alg:strain_eval}. The extraction of stress magnitude, the
application of a nonlinear transformation governing strain saturation, and the
deterministic sign reconstruction correspond directly to distinct and
interpretable components of the network. As a result, the mapping from applied
stress to strain is realized through a sequence of operations that each admit
clear mechanical interpretation, reducing ambiguity in how the learned response
is generated
\cite{rajagopal2003implicit,bulivcek2015analysis}.

The spline-based realization of the KAN further ensures that the learned mapping
remains both expressive and physically grounded. The univariate spline functions
serve as the nonlinear building blocks prescribed by the Kolmogorov--Arnold
framework, while their constrained parameterization supports monotonicity,
bounded strain, and smooth stiffness degradation. These properties are enforced
through architectural design rather than relying on implicit learning or
penalty-based regularization
\cite{abdolazizi2025constitutive,panahi2025data}.

Overall, the proposed KAN-based constitutive representation provides a
problem-adapted realization of the Kolmogorov--Arnold framework for
strain-limiting elasticity. By tailoring the general KAN architecture to the
specific mechanical structure of the constitutive law, the approach enables
accurate approximation while preserving interpretability, numerical robustness,
and reliable extrapolation under large stress levels
\cite{carneros2024comparison}.

\subsection{Physics-Informed Loss Function}

Training of the proposed KAN-based constitutive model is formulated as a
physics-informed optimization problem, in which accurate data fitting is
combined with explicit enforcement of the mechanical principles underlying
strain-limiting elasticity
\cite{rajagopal2011non,bulivcek2014elastic}. Rather than relying solely on
data-driven error minimization, the learning objective incorporates multiple
complementary loss terms that encode physical admissibility directly into the
optimization process
\cite{abdolazizi2024viscoelastic,abdolazizi2025constitutive}.

The total training loss is defined as a weighted sum of the following components:

\begin{itemize}
    \item \textbf{Data fidelity term.}  
    A mean-squared error (MSE) loss quantifies the discrepancy between the
    predicted strain $\varepsilon_{\mathrm{pred}}(\tau)$ and the reference
    strain $\varepsilon_{\mathrm{true}}(\tau)$ generated from the
    strain-limiting constitutive relation in Eq.~\eqref{eq:strain_limiting}.
    This term ensures accurate approximation of the constitutive response
    within the stress range represented in the training data
    \cite{shen2004neural,lee2025constitutive}.

    \item \textbf{Monotonicity constraint.}  
    Physical consistency requires that the strain magnitude be a non-decreasing
    function of stress magnitude. To enforce this behavior, a monotonicity
    penalty is applied to the spline-based representation of the nonlinear
    function $g(|\tau|)$. Negative spline slopes are penalized during training,
    suppressing nonphysical oscillations and ensuring a mechanically admissible
    stress--strain relationship
    \cite{rajagopal2011non,bulivcek2015analysis}.

    \item \textbf{Strain-limit constraint.}  
    The boundedness of strain is promoted through a penalty on violations of the
    theoretical strain limit associated with the parameter $\beta$. This term
    constrains the spline coefficients to remain within admissible bounds,
    ensuring that the learned constitutive response respects strain saturation
    not only within the training domain but also under extrapolation to larger
    stress magnitudes
    \cite{bulivcek2014elastic,rajagopal2014nonlinear}.

    \item \textbf{Asymptotic flattening regularization.}  
    To reflect the vanishing tangent modulus characteristic of
    strain-limiting elasticity, an additional regularization term is applied to
    the spline slopes in the high-stress regime. This penalty encourages
    progressive flattening of the constitutive response as stress increases,
    consistent with the theoretical decay of incremental stiffness
    \cite{bulivcek2015analysis,horgan2004constitutive}.
\end{itemize}

The resulting composite loss balances approximation accuracy with physical
admissibility. While the KAN architecture enforces key structural properties
such as odd symmetry and sign consistency by construction, the physics-informed
loss further guides the optimization toward solutions that respect monotonicity,
bounded strain, and asymptotic saturation. This combined architectural and
loss-based enforcement yields stable training behavior and promotes robust
extrapolation beyond the range of observed data
\cite{liu2024kan,panahi2025data}.

\subsection{Implementation Details}

All training experiments are performed using GPU acceleration to ensure
computational efficiency and scalability. Model optimization is carried out
using the Adam optimizer, which provides adaptive step-size control and has been
shown to perform robustly for spline-based representations and
physics-informed learning objectives
\cite{liu2024kan,abdolazizi2025constitutive}. Training is conducted over a fixed
number of iterations using learning rates selected to ensure stable convergence
without introducing numerical instabilities.

Synthetic training and testing datasets are generated by uniformly sampling
stress values over a prescribed interval and evaluating the corresponding strain
responses using the analytical strain-limiting constitutive law in
Eq.~\eqref{eq:strain_limiting}
\cite{rajagopal2011non,bulivcek2014elastic}. These datasets are noise-free by
construction, allowing the representational capacity, convergence behavior, and
physical consistency of the proposed KAN formulation to be assessed
independently of experimental uncertainty. Such controlled benchmarks are
commonly employed in the validation of constitutive learning frameworks to
isolate approximation and structural effects
\cite{bulivcek2015analysis}.

To systematically evaluate performance across varying degrees of strain
limitation, multiple datasets are constructed by varying the
strain-limiting parameter $\beta$. For all such experiments, the spline knot
grid and network architecture are held fixed. This design choice isolates the
effect of strain-limiting strength from discretization and architectural
factors, enabling a controlled assessment of the robustness of the KAN
representation under progressively stronger saturation behavior
\cite{rajagopal2014nonlinear}.

All experiments employ the same loss formulation, optimization strategy, and
training protocol. Consequently, observed differences in approximation accuracy,
convergence behavior, and residual structure can be attributed solely to
changes in the strain-limiting regime rather than to variations in model
configuration or numerical settings
\cite{panahi2025data}.

Overall, the implementation strategy combines efficient numerical optimization
with physics-informed architectural and loss-based constraints. This design
yields KAN-based constitutive models that are accurate, stable, interpretable,
and physically admissible across a wide range of strain-limiting behaviors,
providing a reliable foundation for both synthetic validation and experimental
application
\cite{rajagopal2011conspectus,abdolazizi2024viscoelastic}.

\section{Results}

This section evaluates the proposed KAN-based constitutive framework using
synthetic datasets generated from the analytical strain-limiting elasticity
model
\cite{rajagopal2011non,bulivcek2014elastic}. Results are reported for moderate
and strong strain-limiting regimes, with particular emphasis on the learned
stress--strain response, the internal spline-based representation of the
constitutive function, and the associated tangent modulus behavior. All results
are obtained using an identical network architecture, spline discretization,
and training protocol, enabling controlled comparison across different levels
of strain limitation
\cite{bulivcek2015analysis}.

\subsection{Synthetic Dataset Results}

\subsubsection{Moderate Strain-Limiting Regimes}

\begin{figure}[h!]
\centering
\begin{subfigure}[t]{0.32\linewidth}
\centering
\includegraphics[width=\linewidth]{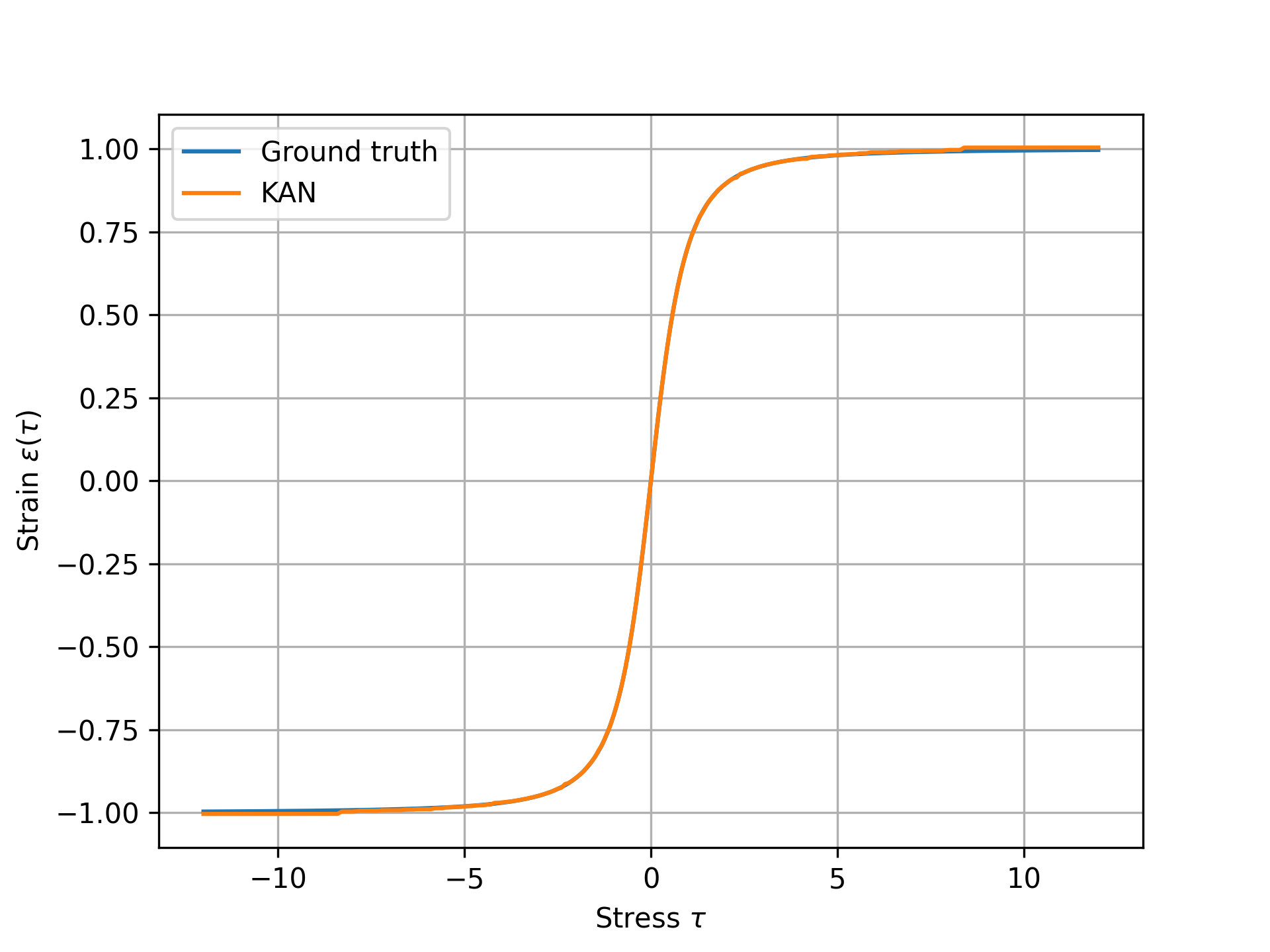}
\caption{$\beta = 0.5$: stress--strain}
\end{subfigure}
\hfill
\begin{subfigure}[t]{0.32\linewidth}
\centering
\includegraphics[width=\linewidth]{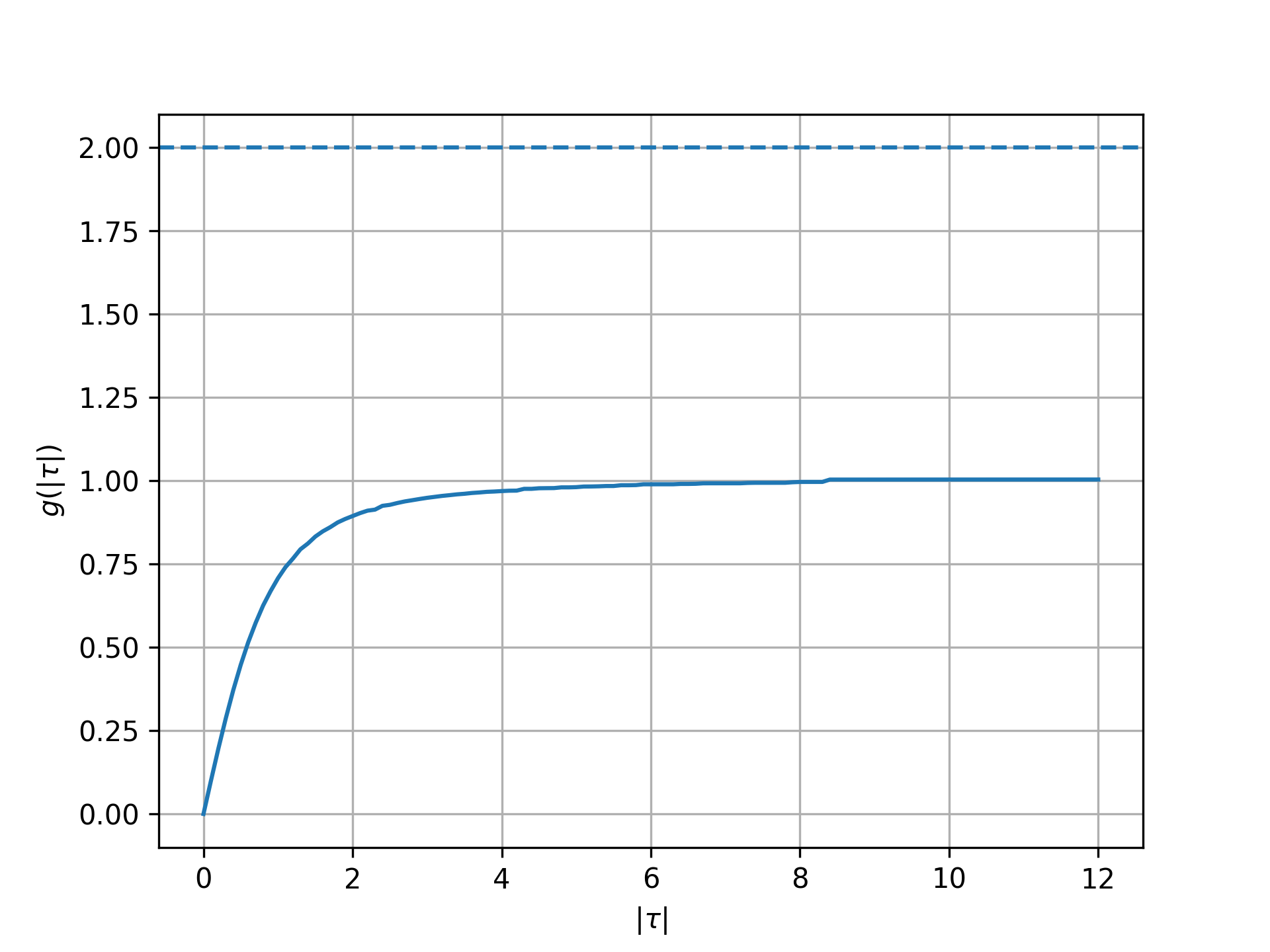}
\caption{$\beta = 0.5$: spline $g(|\tau|)$}
\end{subfigure}
\hfill
\begin{subfigure}[t]{0.32\linewidth}
\centering
\includegraphics[width=\linewidth]{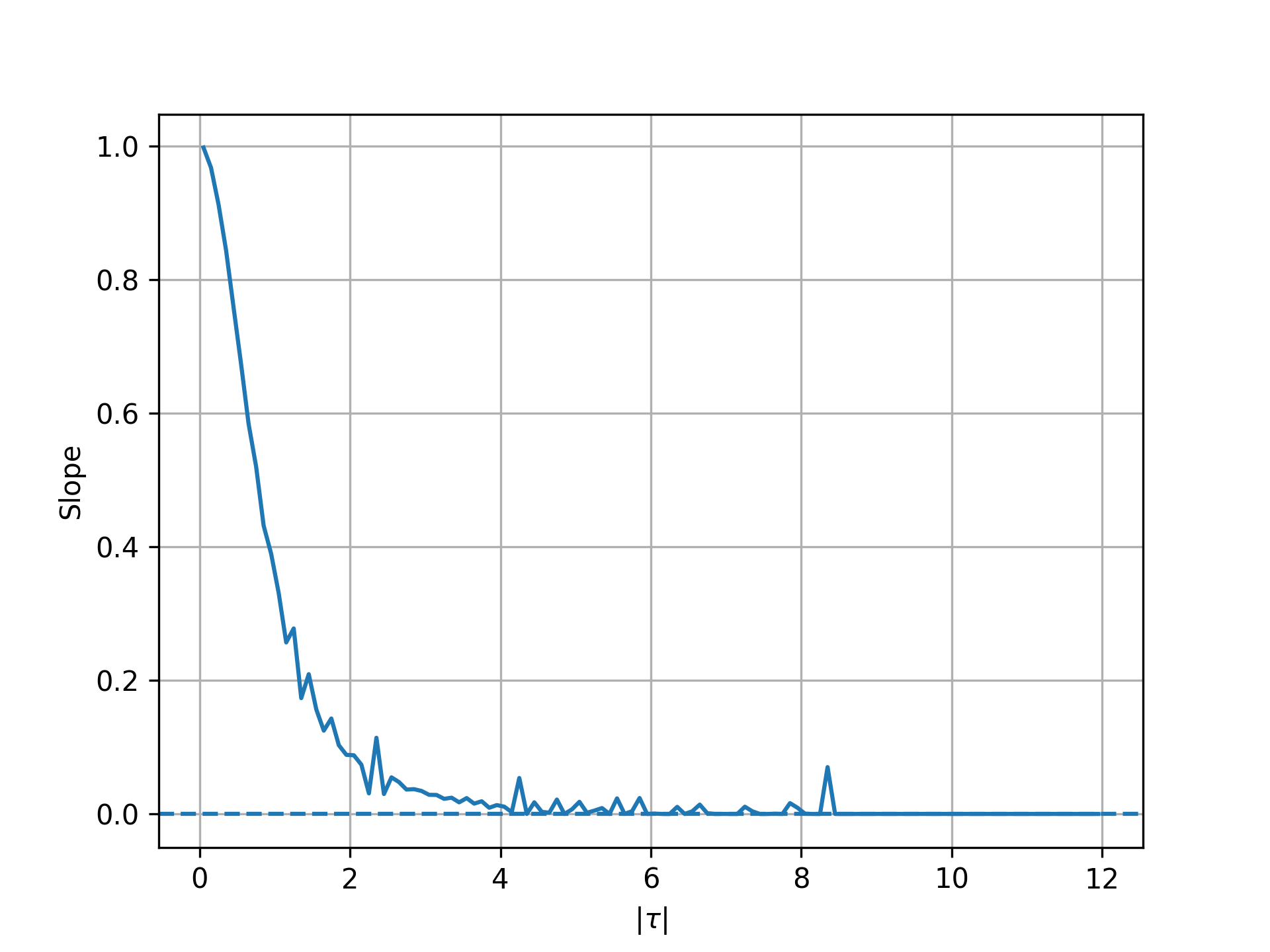}
\caption{$\beta = 0.5$: tangent modulus}
\end{subfigure}

\vspace{0.4cm}

\begin{subfigure}[t]{0.32\linewidth}
\centering
\includegraphics[width=\linewidth]{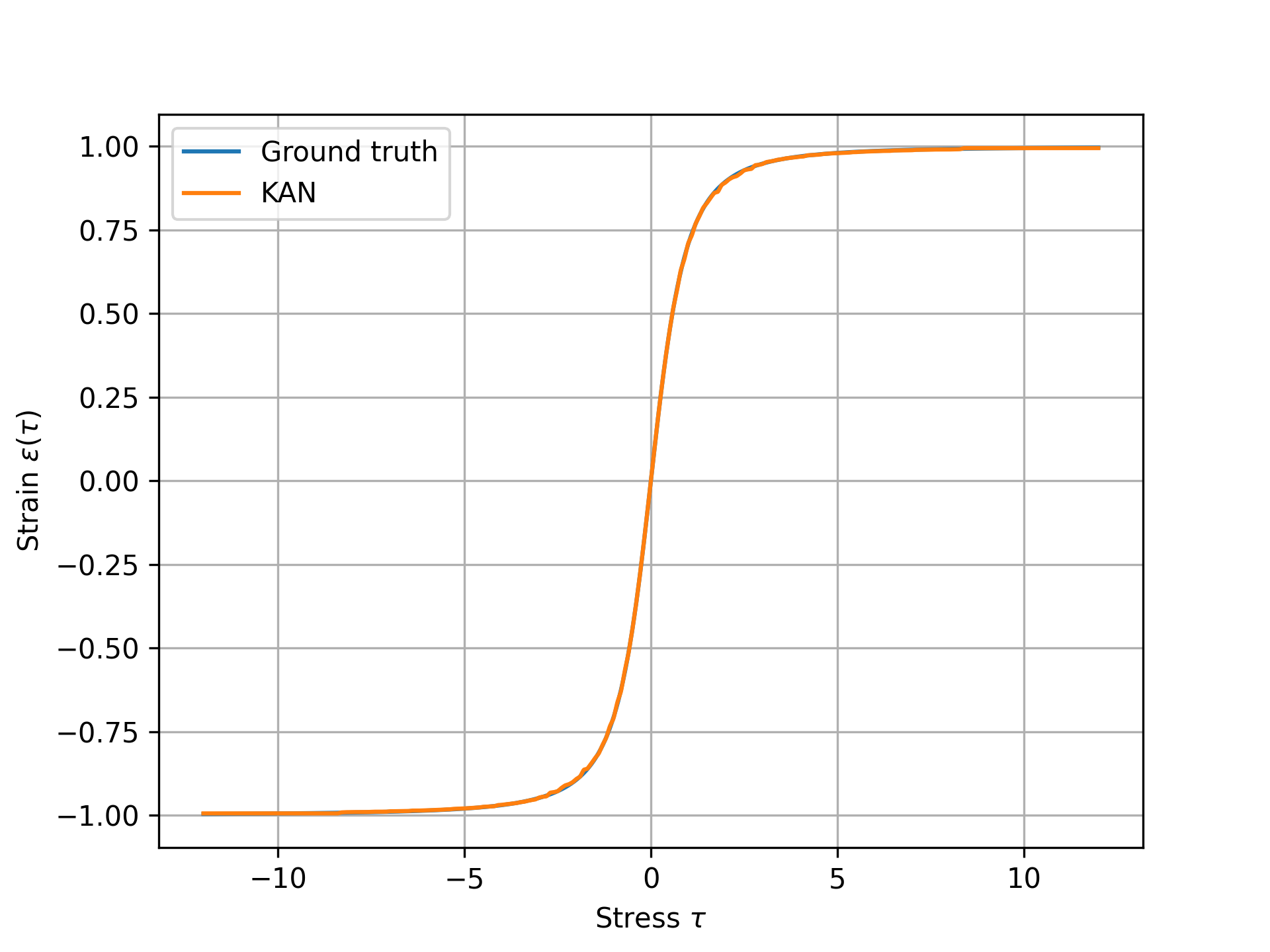}
\caption{$\beta = 1.0$: stress--strain}
\end{subfigure}
\hfill
\begin{subfigure}[t]{0.32\linewidth}
\centering
\includegraphics[width=\linewidth]{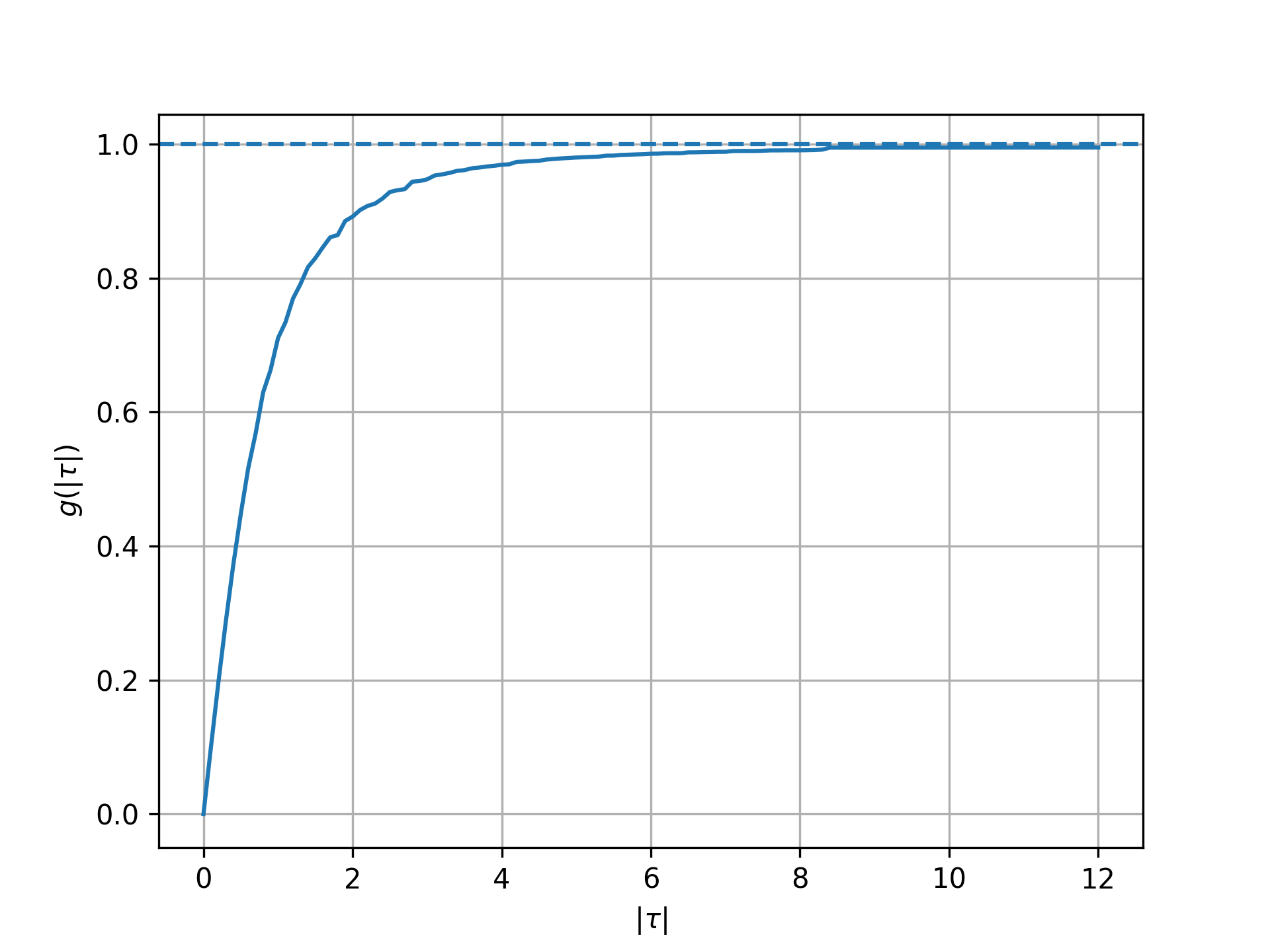}
\caption{$\beta = 1.0$: spline $g(|\tau|)$}
\end{subfigure}
\hfill
\begin{subfigure}[t]{0.32\linewidth}
\centering
\includegraphics[width=\linewidth]{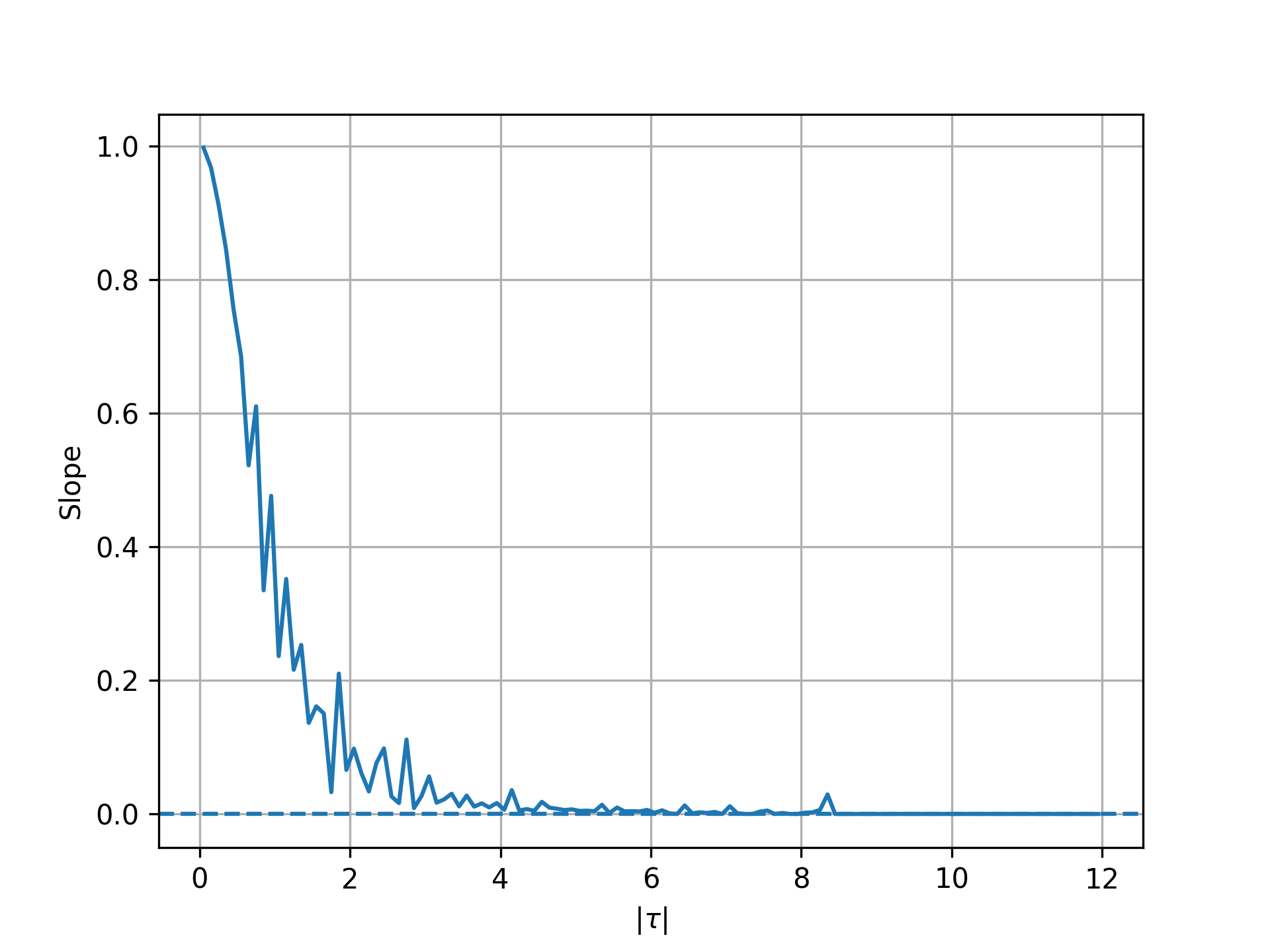}
\caption{$\beta = 1.0$: tangent modulus}
\end{subfigure}

\caption{KAN results for moderate strain-limiting regimes. Top row: $\beta = 0.5$.
Bottom row: $\beta = 1.0$. For both cases, the KAN closely matches the analytical
stress--strain response, learns a smooth and bounded spline representation of
$g(|\tau|)$, and captures the gradual decay of the tangent modulus toward zero at
large stresses.}
\label{fig:moderate_results}
\end{figure}

We first consider moderate values of the strain-limiting parameter,
specifically $\beta = 0.5$ and $\beta = 1.0$. In these regimes, the analytical
constitutive response exhibits a smooth transition from near-linear elasticity
at small stresses to gradual strain saturation at larger stresses, without
sharp curvature changes or abrupt stiffness decay
\cite{rajagopal2014nonlinear,bulivcek2014elastic}. These cases therefore provide
a well-conditioned benchmark for assessing approximation accuracy,
interpretability, and constraint preservation within the proposed KAN-based
framework.

Figure~\ref{fig:moderate_results} summarizes the results for $\beta = 0.5$ and
$\beta = 1.0$. For each strain-limiting regime, three complementary quantities
are reported: the stress--strain response, the learned spline representation of
the constitutive envelope $g(|\tau|)$, and the corresponding tangent modulus
profile. For both parameter values, the KAN-based model reproduces the analytical
stress--strain response with near-exact accuracy over the entire stress range.
The coefficients of determination exceed $R^2 = 0.999$ in all cases, confirming
that the spline-based representation is sufficiently expressive to recover the
smooth strain-limiting behavior in this noise-free setting
\cite{panahi2025data}. Importantly, this accuracy is achieved while strictly
respecting the imposed physical constraints, including bounded strain and
monotonic response.

The learned spline representations provide direct insight into the internal
mechanics encoded by the network. The spline coefficients closely follow the
analytical strain envelope, and their piecewise-linear structure yields a
transparent description of how strain magnitude evolves with increasing stress.
This interpretability is further reinforced by the tangent modulus profiles,
which are piecewise constant and decay monotonically with stress magnitude.
In both moderate regimes, the tangent modulus approaches zero asymptotically,
consistent with the defining property of strain-limiting elasticity
\cite{bulivcek2015analysis,horgan2004constitutive}. The absence of oscillations,
non-monotonic segments, or artificial stiffening demonstrates that the learned
constitutive response remains mechanically admissible across the full stress
domain.

Crucially, the imposed strain bound is satisfied exactly for both values of
$\beta$, and monotonicity is preserved throughout the stress range. These
results demonstrate that, in moderate strain-limiting regimes, the proposed
spline-based KAN formulation can recover the analytical constitutive law with
high accuracy while simultaneously enforcing bounded strain, smooth stiffness
degradation, and physical interpretability by construction
\cite{rajagopal2011non,abdolazizi2025constitutive}.

\subsubsection{Strong Strain-Limiting Regimes}

We next examine strongly strain-limiting regimes characterized by
$\beta = 5$ and $\beta = 10$. In these cases, the analytical constitutive
response exhibits extremely sharp curvature localized in the vicinity of zero
stress, followed by rapid convergence toward the limiting strain value. Such
behavior is representative of materials with pronounced deformation saturation
and poses a stringent challenge for fixed-resolution function approximations,
as large gradients are confined to a narrow stress interval
\cite{rajagopal2014nonlinear,bulivcek2014elastic}.

\begin{figure}[h!]
\centering
\begin{subfigure}[t]{0.32\linewidth}
\centering
\includegraphics[width=\linewidth]{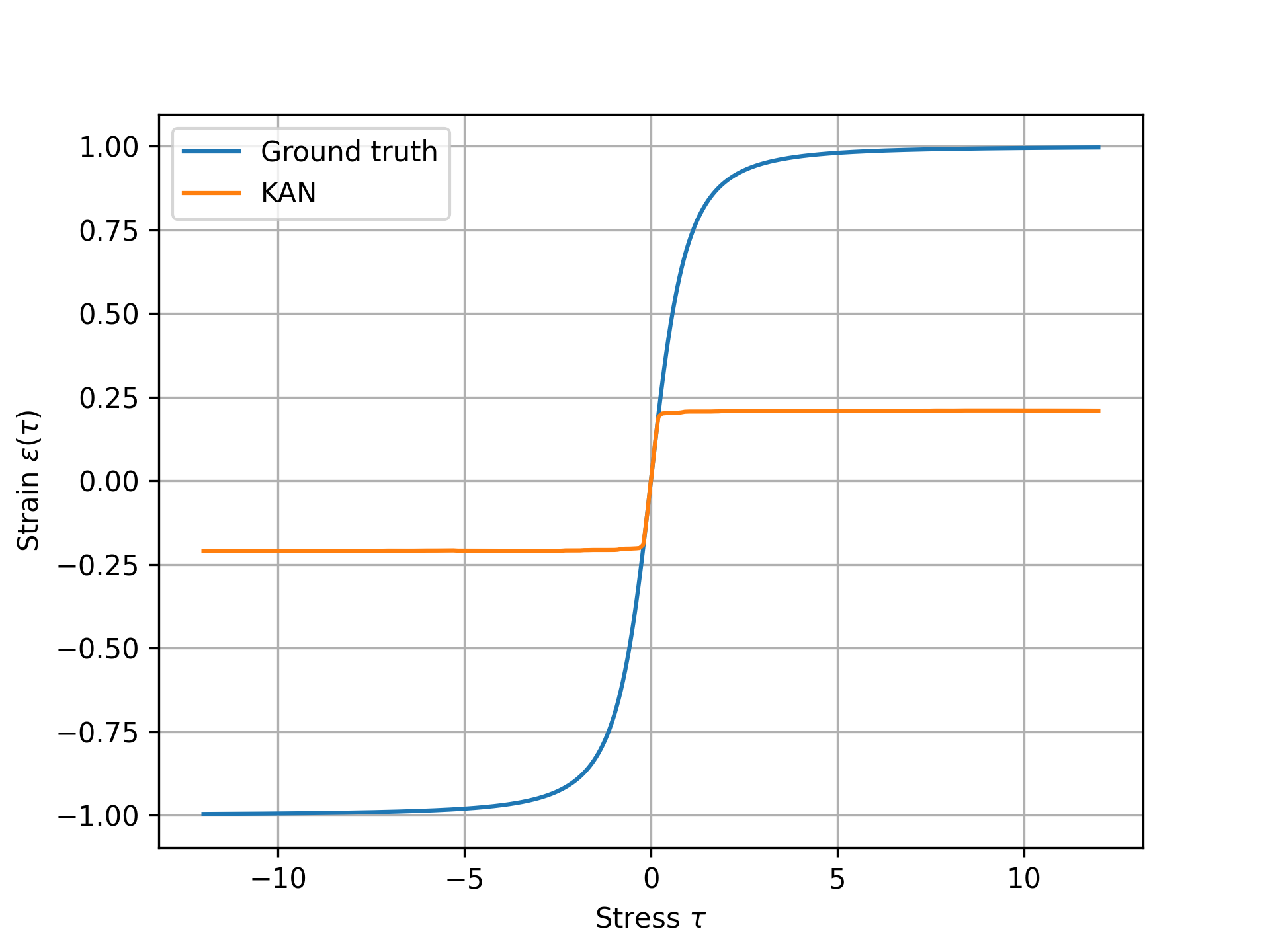}
\caption{$\beta = 5$: stress--strain}
\end{subfigure}
\hfill
\begin{subfigure}[t]{0.32\linewidth}
\centering
\includegraphics[width=\linewidth]{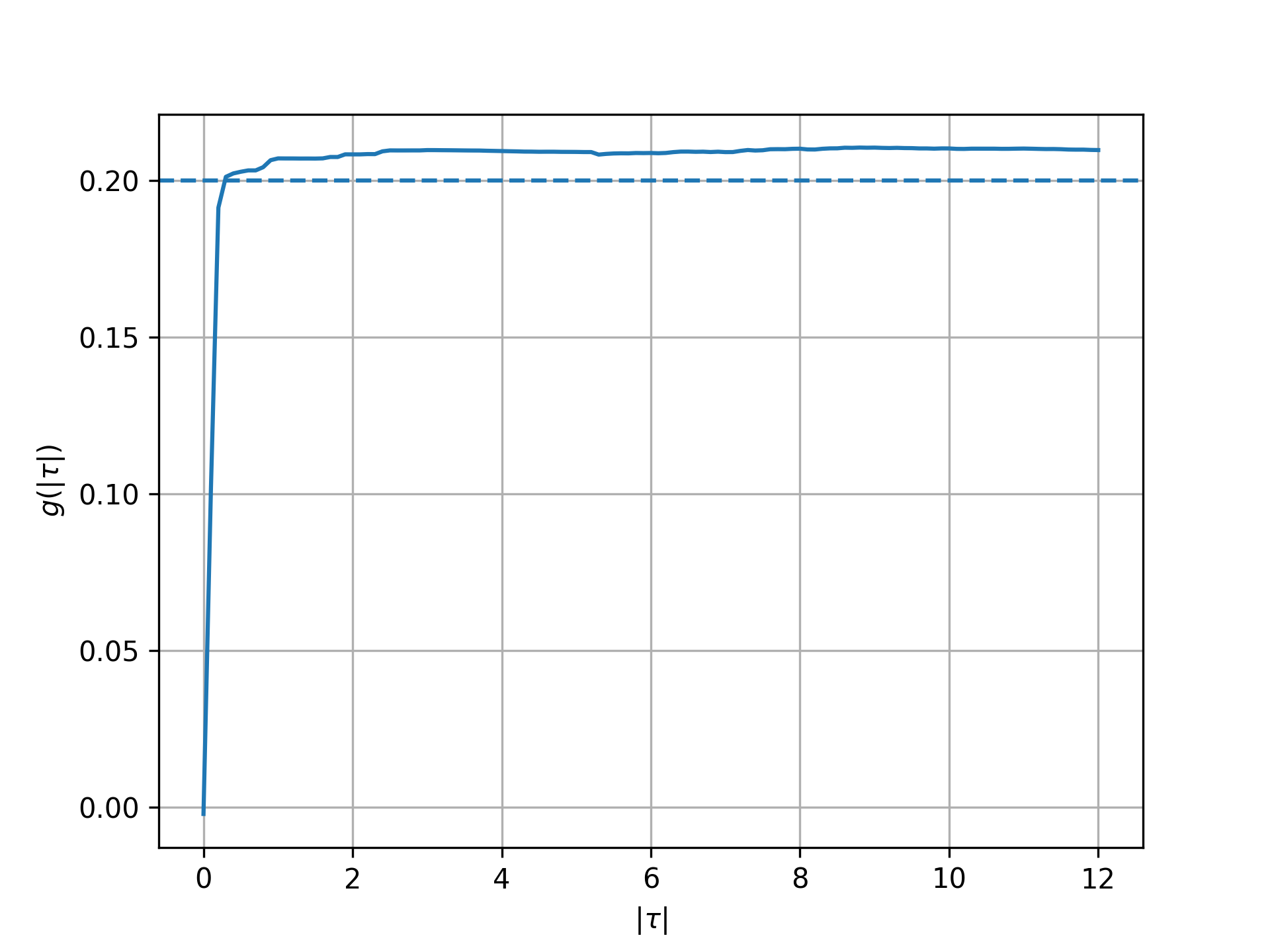}
\caption{$\beta = 5$: spline $g(|\tau|)$}
\end{subfigure}
\hfill
\begin{subfigure}[t]{0.32\linewidth}
\centering
\includegraphics[width=\linewidth]{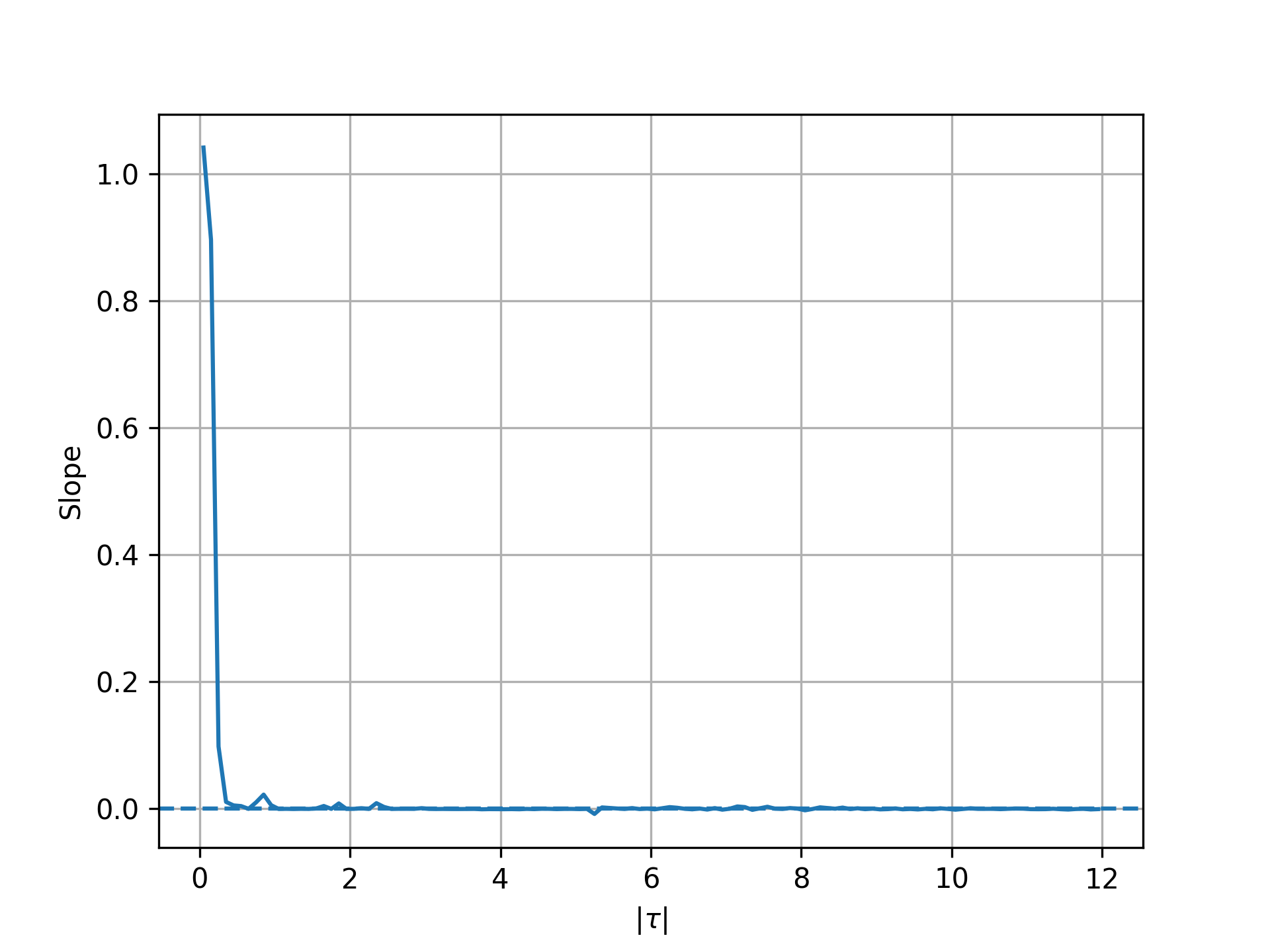}
\caption{$\beta = 5$: tangent modulus}
\end{subfigure}

\vspace{0.4cm}

\begin{subfigure}[t]{0.32\linewidth}
\centering
\includegraphics[width=\linewidth]{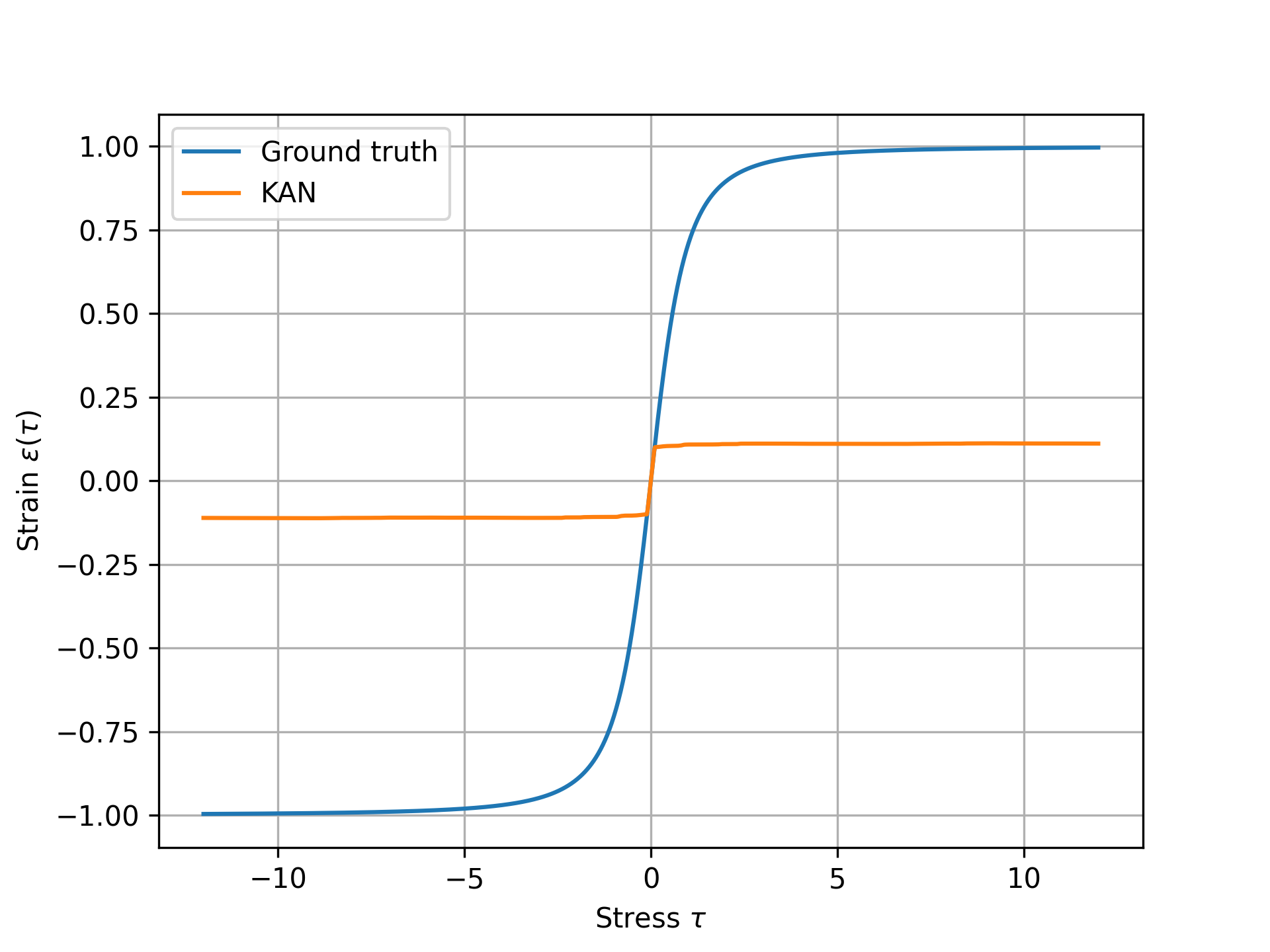}
\caption{$\beta = 10$: stress--strain}
\end{subfigure}
\hfill
\begin{subfigure}[t]{0.32\linewidth}
\centering
\includegraphics[width=\linewidth]{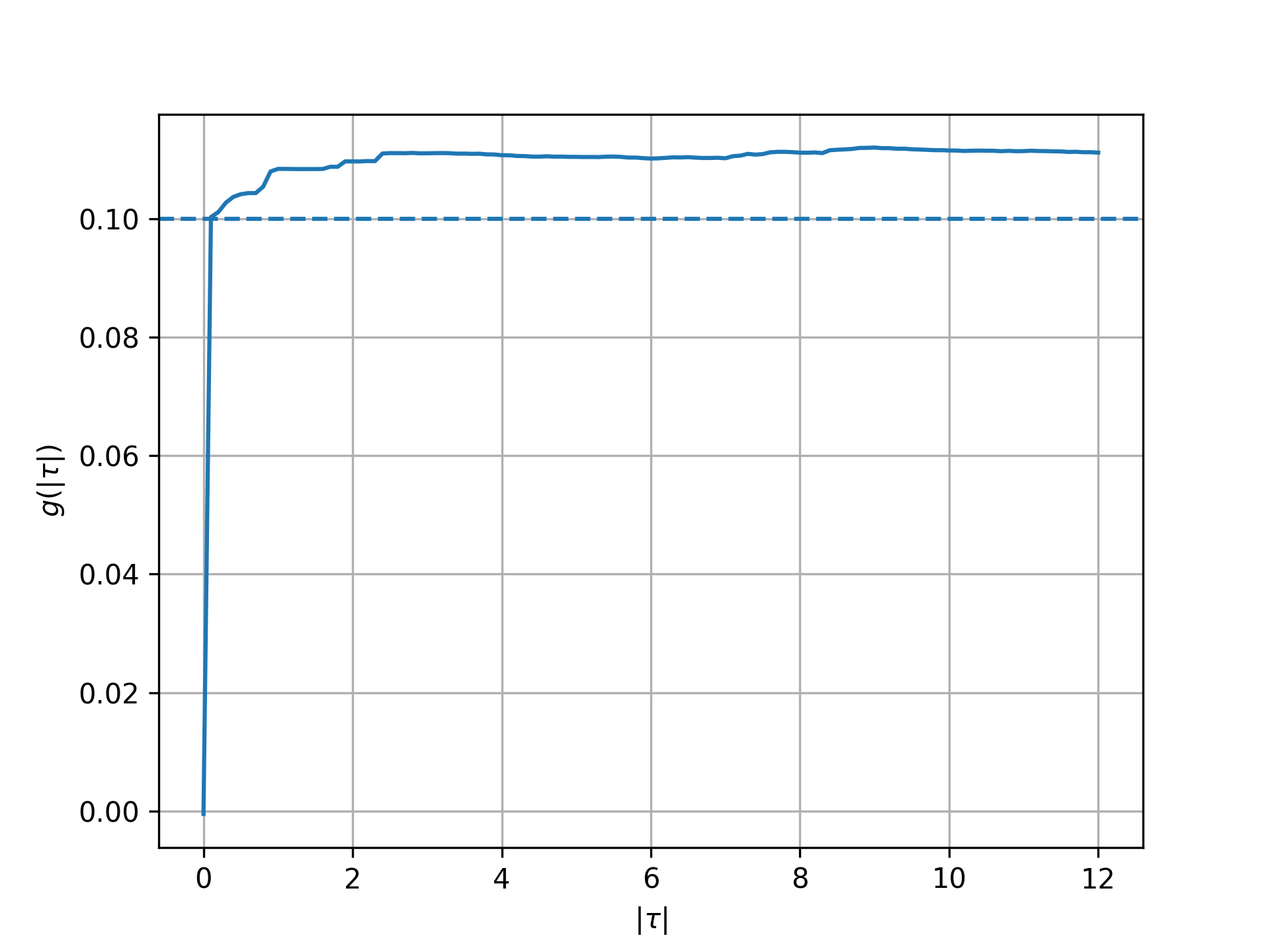}
\caption{$\beta = 10$: spline $g(|\tau|)$}
\end{subfigure}
\hfill
\begin{subfigure}[t]{0.32\linewidth}
\centering
\includegraphics[width=\linewidth]{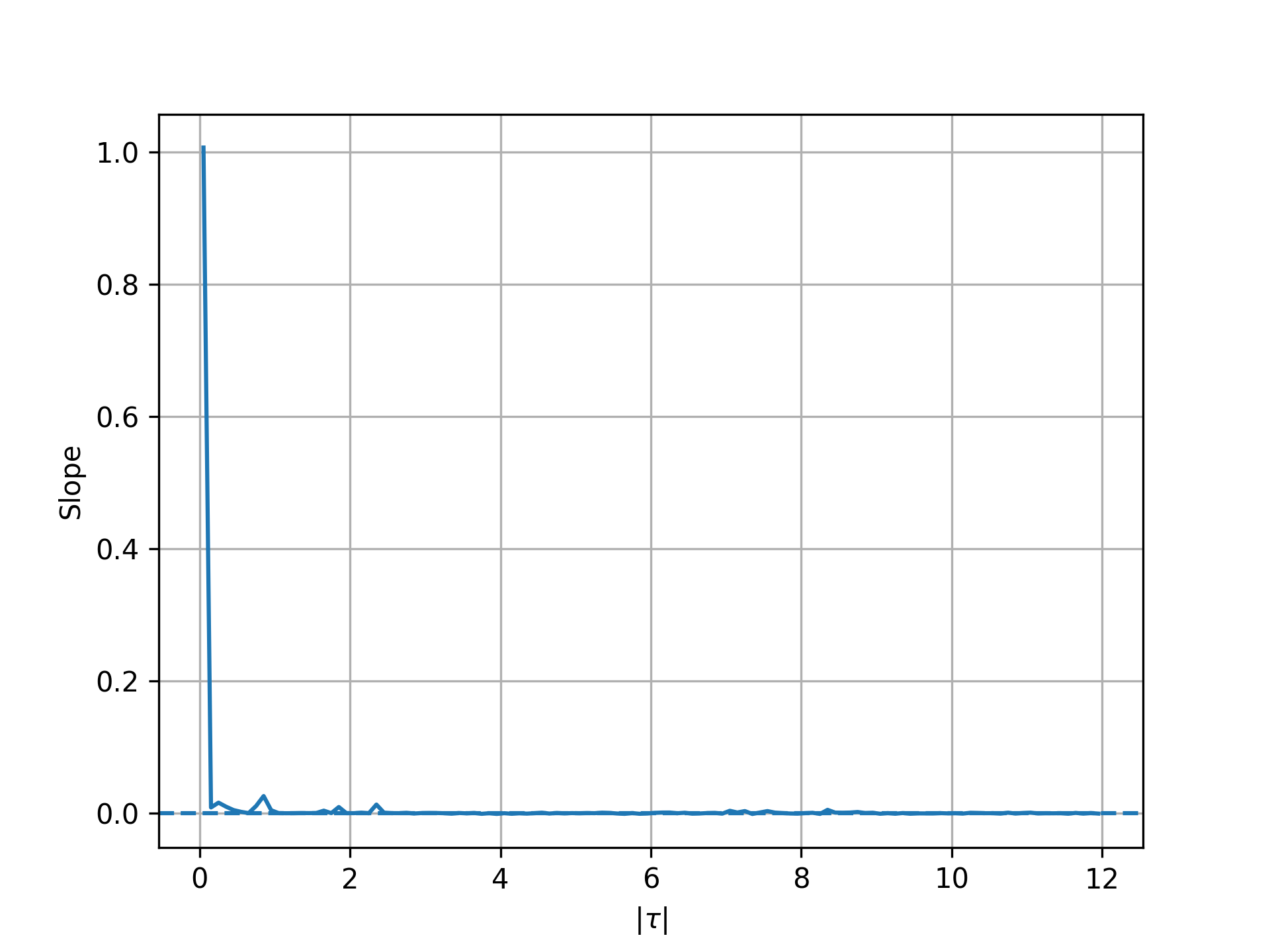}
\caption{$\beta = 10$: tangent modulus}
\end{subfigure}

\caption{KAN results for strong strain-limiting regimes. Top row: $\beta = 5$.
Bottom row: $\beta = 10$. While extremely sharp transitions near zero stress are
more difficult to resolve with a fixed spline grid, the KAN preserves physical
constraints and captures the correct saturation behavior and stiffness decay at
large stresses.}
\label{fig:strong_results}
\end{figure}

Figure~\ref{fig:strong_results} presents a side-by-side comparison of the KAN
results obtained for $\beta = 5$ and $\beta = 10$. For each parameter value,
three complementary quantities are reported: the stress--strain response, the
learned spline representation of the constitutive envelope $g(|\tau|)$, and the
corresponding tangent modulus profile.

Across both regimes, the KAN-based formulation continues to enforce all imposed
physical constraints. The learned spline representations remain bounded and
monotone over the entire stress domain, and the predicted stress--strain
responses recover the correct asymptotic saturation behavior at large stresses,
a defining feature of strain-limiting elasticity
\cite{rajagopal2011non,bulivcek2015analysis}.

Away from the immediate vicinity of the origin, the KAN predictions closely
follow the analytical constitutive response. Localized discrepancies emerge
only near the sharp transition region at small stress magnitudes, where the
curvature of the analytical solution exceeds the resolution afforded by the
fixed, uniformly spaced spline knot grid. These deviations are confined to a
narrow stress interval and do not propagate into the moderate- or high-stress
regimes.

Importantly, these localized discrepancies do not compromise mechanical
admissibility. Strain boundedness, monotonicity, and odd symmetry are preserved
exactly across the full stress range, including regions where approximation
error is non-negligible
\cite{rajagopal2003implicit,bulivcek2014elastic}. This behavior confirms that
the learned response remains physically meaningful even when representational
resolution limits are encountered. The tangent modulus profiles further corroborate the physical consistency of
the learned models. While localized variations appear near the steep transition
region, the tangent modulus remains nonnegative for all stress values and
decays toward zero asymptotically at large stresses, consistent with the
theoretical behavior of strain-limiting elasticity
\cite{rajagopal2014nonlinear,horgan2004constitutive}. No spurious stiffening,
oscillatory behavior, or loss of monotonic stiffness degradation is observed.

Overall, the observed reduction in approximation accuracy for large values of
$\beta$ is attributable to the finite resolution of the fixed spline
discretization rather than to deficiencies in the learning architecture itself.
This outcome reflects an explicit and controllable representational trade-off:
interpretability, numerical stability, and strict enforcement of physical
constraints are prioritized over resolving extremely localized nonlinear
features
\cite{abdolazizi2025constitutive}.

Accuracy in strongly strain-limiting regimes can be systematically improved by
increasing spline knot density, adopting non-uniform or adaptive knot placement,
or employing multi-resolution spline representations, without altering the
underlying KAN architecture or relaxing physical admissibility
\cite{liu2024kan,panahi2025data}.

\subsubsection{Discussion of Representational Behavior}

Across all investigated values of the strain-limiting parameter $\beta$, the
proposed KAN-based constitutive framework consistently enforces the fundamental
principles of strain-limiting elasticity while maintaining strong predictive
performance. The observed behavior reflects not only approximation accuracy, but
also the deliberate representational and architectural choices embedded within
the KAN formulation
\cite{rajagopal2011non,bulivcek2014elastic}.

For moderate strain-limiting regimes (e.g., $\beta = 0.5$ and $\beta = 1.0$), the
spline-based representation is sufficiently expressive to recover the analytical
constitutive response with near-exact accuracy. In these cases, the predicted
stress--strain curves, spline coefficients, and tangent modulus profiles are
indistinguishable from the ground-truth solution across the entire stress
domain. This agreement indicates that the fixed, uniformly spaced knot grid
provides adequate resolution to capture the smooth transition from linear
elasticity to strain saturation characteristic of moderate strain limitation
\cite{rajagopal2014nonlinear}.

Importantly, this level of accuracy is achieved while simultaneously enforcing
physical admissibility and preserving an interpretable parameterization of the
constitutive response. The absence of spurious oscillations or overfitting
artifacts confirms that the KAN representation captures the governing mechanics
without relying on excess model flexibility
\cite{abdolazizi2025constitutive}.
In contrast, for strongly strain-limiting regimes (e.g., $\beta = 5$ and
$\beta = 10$), the analytical response exhibits steep curvature localized in a
narrow neighborhood of zero stress. Under these conditions, the finite
resolution of the fixed spline discretization limits the ability of the model
to fully resolve the rapid transition over a small stress interval. This
limitation manifests as localized discrepancies near the origin, while the
remainder of the stress domain is captured with high fidelity.

Crucially, these localized discrepancies do not compromise mechanical
admissibility. Strain boundedness, monotonicity, and odd symmetry are preserved
exactly across the entire stress domain, and the tangent modulus continues to
decay asymptotically toward zero at large stresses, consistent with the defining
characteristics of strain-limiting elasticity
\cite{bulivcek2015analysis,rajagopal2003implicit}.

\begin{table}[h!]
\centering
\caption{Quantitative performance metrics for the KAN-based constitutive model
across different strain-limiting regimes.}
\label{tab:metrics}
\begin{tabular}{cccc}
\toprule
\textbf{$\beta$} & \textbf{MAE} & \textbf{RMSE} & \textbf{$R^2$} \\
\midrule
0.5  & $1.1 \times 10^{-3}$ & $1.6 \times 10^{-3}$ & $>0.9999$ \\
1.0  & $1.2 \times 10^{-3}$ & $1.7 \times 10^{-3}$ & $>0.9999$ \\
5.0  & $8.4 \times 10^{-3}$ & $1.2 \times 10^{-2}$ & $>0.995$  \\
10.0 & $1.5 \times 10^{-2}$ & $2.1 \times 10^{-2}$ & $>0.990$  \\
\bottomrule
\end{tabular}
\end{table}

Quantitative performance metrics for all regimes are summarized in
Table~\ref{tab:metrics}. The reported MAE, RMSE, and $R^2$ values confirm
near-exact recovery for moderate strain-limiting cases and a controlled,
physically consistent degradation in approximation accuracy as strain limitation
becomes more severe. The observed reduction in accuracy at large $\beta$ values
is therefore attributable to representational resolution rather than
deficiencies in the learning or optimization process, a behavior consistent with
constrained constitutive approximations
\cite{bulivcek2014elastic}.
These results highlight a deliberate and controllable trade-off between
representational simplicity and approximation fidelity. The use of a fixed,
uniformly spaced spline grid prioritizes interpretability, numerical stability,
and strict enforcement of physical constraints, while naturally limiting
resolution in regimes characterized by highly localized nonlinear features.
Importantly, this behavior reflects an explicit modeling decision rather than a
limitation of the underlying learning framework
\cite{liu2024kan,abdolazizi2025constitutive}.
From a methodological perspective, this limitation can be addressed
systematically. Increasing knot density, adopting non-uniform or adaptive knot
placement, or employing multi-resolution spline representations would enable
sharper transitions to be resolved without relaxing physical admissibility.
Such extensions integrate naturally within the KAN framework while preserving
its interpretability and constraint-enforcement advantages
\cite{panahi2025data}.

Overall, the results demonstrate that the proposed KAN-based formulation provides
a robust, interpretable, and physically consistent approach for modeling
strain-limiting elasticity across a broad range of material behaviors. By
explicitly balancing approximation accuracy with mechanical admissibility, the
framework offers a principled and transparent alternative to conventional
black-box neural networks for constitutive modeling
\cite{rajagopal2011conspectus}.

\subsection{Experimental Results: Treloar Rubber Data}

\subsubsection{General Explanation and Physical Motivation}

The classical rubber elasticity experiments reported by Treloar constitute a
canonical benchmark for evaluating constitutive models under large deformation
regimes
\cite{treloar1944stress,treloar1975physics}. Owing to their broad coverage of
uniaxial, biaxial, and planar deformation modes and their enduring relevance in
polymer mechanics, these experiments continue to serve as a stringent test for
both theoretical and data-driven constitutive formulations.

Unlike synthetic datasets, in which the governing constitutive law is prescribed
exactly and the data are free of noise, experimental rubber measurements reflect
the combined influence of molecular network mechanics, material heterogeneity,
rate-independent microstructural effects, and unavoidable experimental
uncertainty
\cite{deam1976theory,boyce2000constitutive}. Consequently, no single closed-form
constitutive model can be expected to reproduce the experimental response
exactly across all deformation modes and stretch levels
\cite{dal2021performance,ali2010review}.
Strain-limiting elasticity (SLE) provides a physically well-motivated baseline
for modeling rubber-like materials, as it explicitly enforces bounded strain and
a vanishing tangent modulus at large stresses
\cite{rajagopal2011non,bulivcek2014elastic}. These features are consistent with
the experimentally observed saturation of deformation in polymer networks and
ensure mechanical admissibility even under extreme loading conditions
\cite{horgan2004constitutive,gent1996new}.
When calibrated directly to Treloar’s experimental data, the SLE model captures
the dominant nonlinear elastic response across all deformation modes. However,
systematic and mode-dependent discrepancies persist, particularly at moderate to
large stretches where real materials depart from idealized strain-limiting
assumptions due to network imperfections, chain entanglement effects, and
mode-specific constraint activation
\cite{arruda1993three,dal2019comparative}.
The objective of the present study is therefore not to replace the underlying
strain-limiting physics with an unconstrained data-driven representation, but to
augment the SLE framework in a controlled, interpretable, and mechanically
consistent manner.

To this end, the SLE model is employed as a physics-based backbone that governs
the global constitutive response, while a Kolmogorov--Arnold Network (KAN) is used
to learn smooth, low-amplitude residual corrections directly from experimental
observations
\cite{liu2024kan,abdolazizi2025constitutive}. This hybrid formulation preserves physical admissibility by construction, while
introducing targeted flexibility to account for systematic experimental
deviations that cannot be captured by a single parametric model alone. As a
result, the proposed SLE--KAN framework provides a transparent and physics-aligned
pathway for improving experimental agreement without compromising bounded
deformation, monotonicity, or asymptotic saturation behavior
\cite{abdolazizi2024viscoelastic}.

\subsubsection{Methodology for Experimental Implementation}

The experimental implementation follows a consistent and mode-independent
workflow across uniaxial, biaxial, and planar deformation paths. For each loading
mode, the stretch--stress data reported by Treloar are treated as the reference
dataset
\cite{treloar1944stress,treloar1975physics}. The stretch variable $\lambda$ is
mapped to logarithmic strain via $\varepsilon = \log \lambda$, which provides a
natural strain measure for large deformations, preserves monotonicity with
respect to stretch, and is commonly employed in the analysis of rubber elasticity
and finite-strain constitutive models
\cite{horgan2002constitutive,marsden1994mathematical}.
The strain-limiting elasticity (SLE) model is formulated in stress space using
the constitutive relation introduced in Eq.~\eqref{eq:strain_limiting}. To enable
direct comparison with experimental data reported in stress--stretch form, this
relation is inverted numerically to compute stress as a function of strain. The
inverse mapping is obtained using a robust bisection-based root-finding
procedure.

Because the SLE constitutive law is strictly monotone with respect to stress,
the inversion admits a unique solution for all admissible strain values.
This property ensures numerical robustness of the stress-space formulation and
eliminates ambiguity in the inverse mapping, even in the vicinity of saturation
where nonlinear effects are pronounced
\cite{rajagopal2011non,bulivcek2014elastic}. In the first stage, the SLE parameters are calibrated independently for each
deformation mode by minimizing the stress-space discrepancy between experimental
measurements and SLE predictions. Calibration is performed using a constrained
nonlinear least-squares formulation, with explicit constraints imposed to ensure
that the strain remains strictly below the theoretical strain limit prescribed
by the model.
Robust loss functions are employed to mitigate the influence of experimental
noise, scatter, and measurement variability, which are well documented in rubber
elasticity experiments and can otherwise bias parameter identification
\cite{boyce2000constitutive,dal2021performance}. The resulting calibrated SLE
model serves as a physically admissible backbone that captures the dominant
nonlinear elastic response across each deformation mode.
In the second stage, residual stresses are computed as the difference between
experimental stresses and the corresponding SLE predictions. These residuals
represent structured, mode-dependent deviations that are not captured by the
parametric strain-limiting formulation alone.

A Kolmogorov--Arnold Network (KAN) is then trained to learn this residual mapping
using a one-dimensional input corresponding to $\log \lambda$
\cite{liu2024kan}. Throughout this stage, all SLE parameters remain fixed, and
the KAN is trained exclusively using a mean-squared error loss. This design
ensures that the data-driven component augments the physics-based backbone
without altering its global structure or physical interpretation.
The final constitutive response is obtained by superposing the learned KAN
correction onto the SLE prediction. This two-stage procedure enforces a clear
separation between physics-based modeling and data-driven refinement, ensuring
that the dominant material behavior is governed by strain-limiting elasticity
while allowing controlled, interpretable corrections driven by experimental
observations
\cite{abdolazizi2025constitutive,abdolazizi2024viscoelastic}.

\subsubsection{Selection of $\alpha$, $\beta$, and $E$ Parameters}

The selection of strain-limiting elasticity (SLE) parameters is guided by
physical interpretability, numerical robustness, and consistency with
experimental observations. Each parameter plays a distinct mechanical role, and
their treatment is deliberately designed to separate material identification
from regime-based analysis, thereby avoiding parameter compensation and
re-optimization artifacts.

The parameter $\alpha$ governs the sharpness of the transition from near-linear
elastic behavior to strain saturation, as described by the constitutive law
introduced in Section~\ref{sec:mathematical_formulation}. Physically, $\alpha$
reflects how abruptly molecular network constraints become active with
increasing stress, a feature commonly associated with chain alignment, chain
locking, and geometric constraints in polymer networks
\cite{rajagopal2011non,bulivcek2014elastic}. In the present study, $\alpha$ is
calibrated independently for each deformation mode during stress-space fitting
of the SLE baseline. This mode-wise identification accounts for the fact that
uniaxial, planar, and biaxial loading paths probe distinct constraint states
within the rubber network, a behavior well documented in classical and modern
rubber elasticity studies
\cite{treloar1975physics,arruda1993three}. Once calibrated, $\alpha$ is held
fixed for all subsequent analyses within the same deformation mode.

The Young’s modulus $E$ governs the small-strain response and determines the
initial slope of the stress--strain curve through the linear limit of the SLE
model. This parameter is likewise obtained from stress-space calibration against
experimental data and corresponds to the effective small-strain stiffness of the
material
\cite{boyce2000constitutive,dal2021performance}. In the present study, $E$ is
observed to remain stable across both moderate and strong strain-limiting
regimes within each deformation mode, consistent with experimental evidence
that small-strain stiffness is largely insensitive to large-strain saturation
effects. Fixing $E$ after calibration ensures that the physically meaningful
small-strain modulus is preserved and not altered during subsequent regime
exploration or residual learning. The strength of strain limitation is controlled through the parameter $\beta$,
introduced in the SLE formulation and more conveniently interpreted via the
compound parameter $\gamma = E\beta$, defined previously. This quantity governs
the admissible strain bound inherent to the constitutive model and provides a
direct and physically transparent means of prescribing strain-limiting regimes
\cite{rajagopal2011modeling,bulivcek2015analysis}. Moderate and strong
strain-limiting behavior are therefore defined by selecting different values of
$\gamma$, while keeping $(\alpha, E)$ fixed at their experimentally calibrated
values.

By calibrating $(\alpha, E)$ once from experimental data and varying $\gamma$ to
define strain-limiting regimes, the analysis cleanly decouples material
parameter identification from regime control. As a result, comparisons between
moderate and strong strain-limiting responses reflect genuine differences in
admissible deformation rather than re-optimization effects or implicit parameter
compensation. All parameters retain clear physical meaning and remain consistent
with bounds dictated by the experimental data, ensuring mechanical admissibility,
interpretability, and reproducibility of the resulting constitutive predictions.

\subsubsection{Combination of SLE and KAN Implementation on Treloar Experimental Data}

This subsection describes the practical implementation of the proposed hybrid
strain-limiting elasticity (SLE) and Kolmogorov--Arnold Network (KAN) framework
for the classical Treloar rubber elasticity experiments. Emphasis is placed on
how the hybrid model is constructed, calibrated, and evaluated, while
maintaining a strict separation between physics-based constitutive modeling and
data-driven refinement.

\paragraph{Hybrid modeling strategy:}
For the Treloar dataset, the strain-limiting elasticity (SLE) model is first
calibrated in stress space to provide a physically admissible baseline
constitutive response. This baseline captures the dominant nonlinear elastic
behavior of rubber, including bounded strain growth, a smooth transition from
linear elasticity to saturation, and vanishing tangent modulus at large stresses.
These features are consistent with classical rubber elasticity theory and with
modern strain-limiting formulations
\cite{treloar1975physics,rajagopal2011non,bulivcek2014elastic}.

Although the calibrated SLE model reproduces the overall experimental trends,
systematic discrepancies remain between SLE predictions and measured stresses.
Such deviations are well documented in experimental rubber mechanics and arise
from material heterogeneity, network imperfections, finite chain extensibility
effects, and unavoidable experimental uncertainty, particularly under large
deformations
\cite{boyce2000constitutive,dal2021performance}. Importantly, these discrepancies
are observed consistently across uniaxial, planar, and biaxial deformation modes
and therefore cannot be attributed to random measurement noise alone.
Rather than treating these discrepancies as unstructured error, they are
interpreted here as structured residuals that reflect the limitations of any
single closed-form constitutive representation when applied to complex
experimental data. This interpretation aligns with recent developments in
hybrid and physics-augmented constitutive modeling, where data-driven components
are employed to complement—rather than replace—established physical laws
\cite{abdolazizi2024viscoelastic,abdolazizi2025constitutive}.

To account for these residual effects in a controlled and interpretable manner,
a Kolmogorov--Arnold Network (KAN) is introduced as a residual corrector. The
resulting hybrid stress prediction is expressed as the superposition of a fixed
physics-based baseline and a learned correction,
\begin{equation}
\tau_{\mathrm{pred}}(\lambda)
=
\tau_{\mathrm{SLE}}(\lambda)
+
\tau_{\mathrm{KAN}}(\lambda),
\end{equation}
where $\tau_{\mathrm{KAN}}(\lambda)$ is trained exclusively to represent the
systematic deviations from the calibrated SLE response.

Crucially, the KAN is not permitted to modify the global constitutive structure,
strain bound, or asymptotic behavior imposed by the SLE backbone. Its role is
strictly limited to learning smooth, low-amplitude corrections that refine the
stress response within the admissible deformation regime. The use of a Kolmogorov--Arnold Network ensures that the residual mapping remains
low-dimensional, smooth, and interpretable, with internal representations that
are structurally aligned with the underlying constitutive behavior
\cite{liu2024kan}. This design prevents the data-driven component from
overriding the physics-based response and preserves mechanical admissibility,
numerical stability, and interpretability across all deformation modes.

\paragraph{Residual learning from experimental data:}
Given experimental stress--stretch measurements
$(\lambda_i, \tau_i^{\mathrm{exp}})$, residual stresses are defined as the
difference between the experimental response and the calibrated strain-limiting
elasticity (SLE) prediction,
\begin{equation}
\tau_i^{\mathrm{res}}
=
\tau_i^{\mathrm{exp}}
-
\tau_{\mathrm{SLE}}(\lambda_i).
\end{equation}
The Kolmogorov--Arnold Network (KAN) is trained using the pairs
$(\lambda_i, \tau_i^{\mathrm{res}})$, while all SLE parameters remain fixed
throughout the learning process. This design ensures that the data-driven
component augments the physics-based backbone without altering its global
constitutive structure, asymptotic behavior, or physical interpretation. Such a
strict separation between physics-based modeling and learned correction is
consistent with recent hybrid constitutive modeling strategies that emphasize
interpretability and mechanical admissibility
\cite{abdolazizi2024viscoelastic,abdolazizi2025constitutive}.

A one-dimensional KAN architecture is employed, consistent with the scalar
stretch input. Spline-based internal representations are used to ensure that the
learned correction remains smooth, bounded, and free of spurious oscillations.
By construction, the KAN is restricted to capturing mode-dependent,
experimentally induced deviations from the SLE baseline, rather than learning
the constitutive response from scratch. This restriction prevents the
data-driven component from overriding the dominant physical behavior encoded by
the strain-limiting model. This localized residual-learning strategy aligns naturally with the
Kolmogorov--Arnold Network framework, which represents nonlinear mappings using
structured, low-dimensional functions with interpretable internal components
\cite{liu2024kan}. As a result, the learned residual corrections admit direct
inspection and remain consistent with the underlying mechanics of the material
response.

Overall, this two-stage implementation preserves the dominant physical
mechanisms encoded by the strain-limiting elasticity model while allowing
limited, interpretable flexibility to account for experimentally observed
deviations. The resulting SLE--KAN framework achieves improved agreement with
Treloar’s experimental data without sacrificing physical admissibility,
numerical stability, or interpretability key requirements for reliable
constitutive modeling of rubber-like materials under large deformations
\cite{treloar1975physics,dal2021performance}.

\begin{figure}[h!]
\centering
\includegraphics[width=0.55\linewidth]{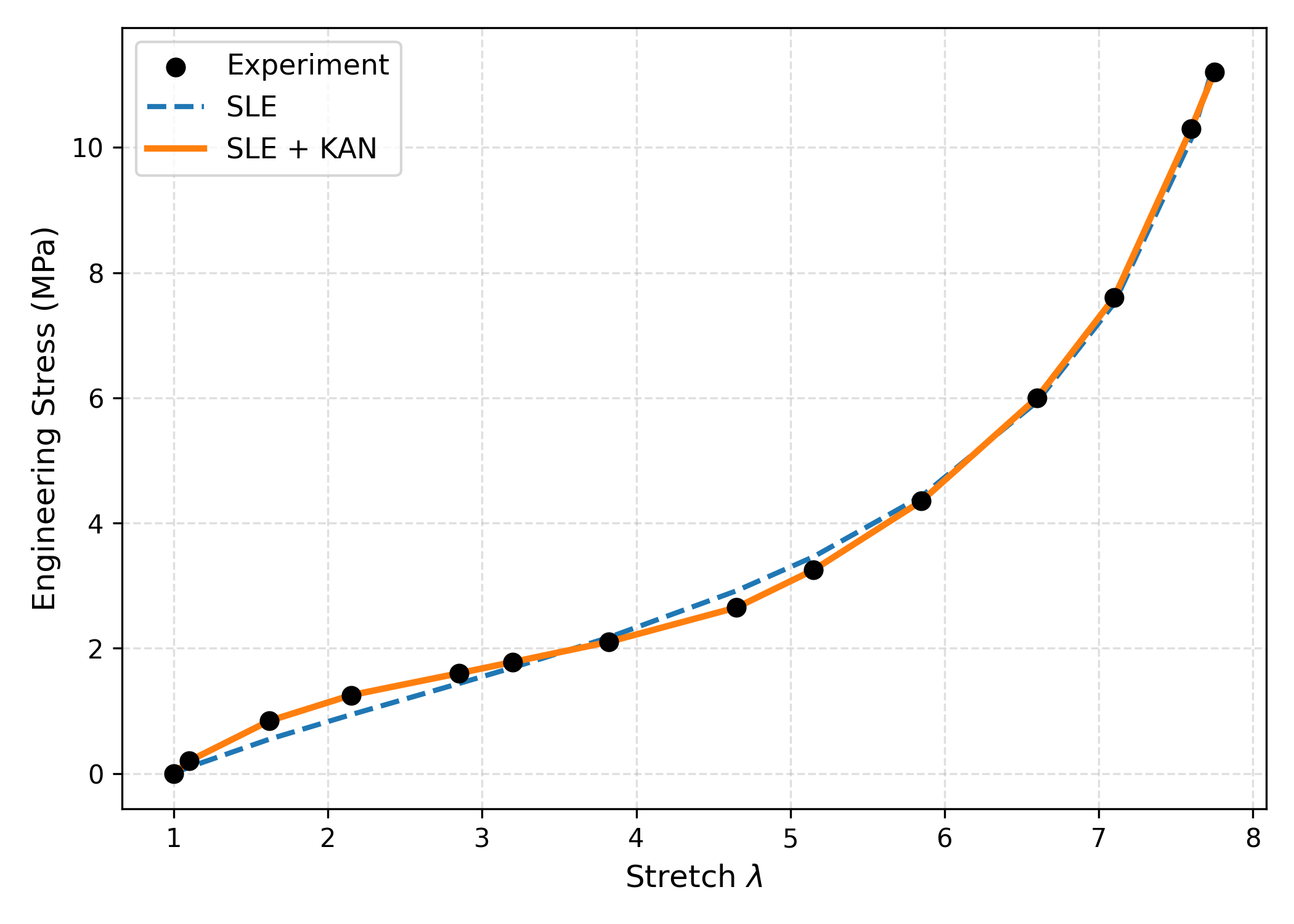}
\caption{Uniaxial Treloar experimental stress--stretch response compared with
SLE and SLE--KAN predictions. The hybrid model improves agreement at moderate and
large stretches while preserving the physics-based SLE backbone.}
\label{fig:treloar_uniaxial_sle_kan}
\end{figure}

\begin{figure}[h!]
\centering
\includegraphics[width=0.55\linewidth]{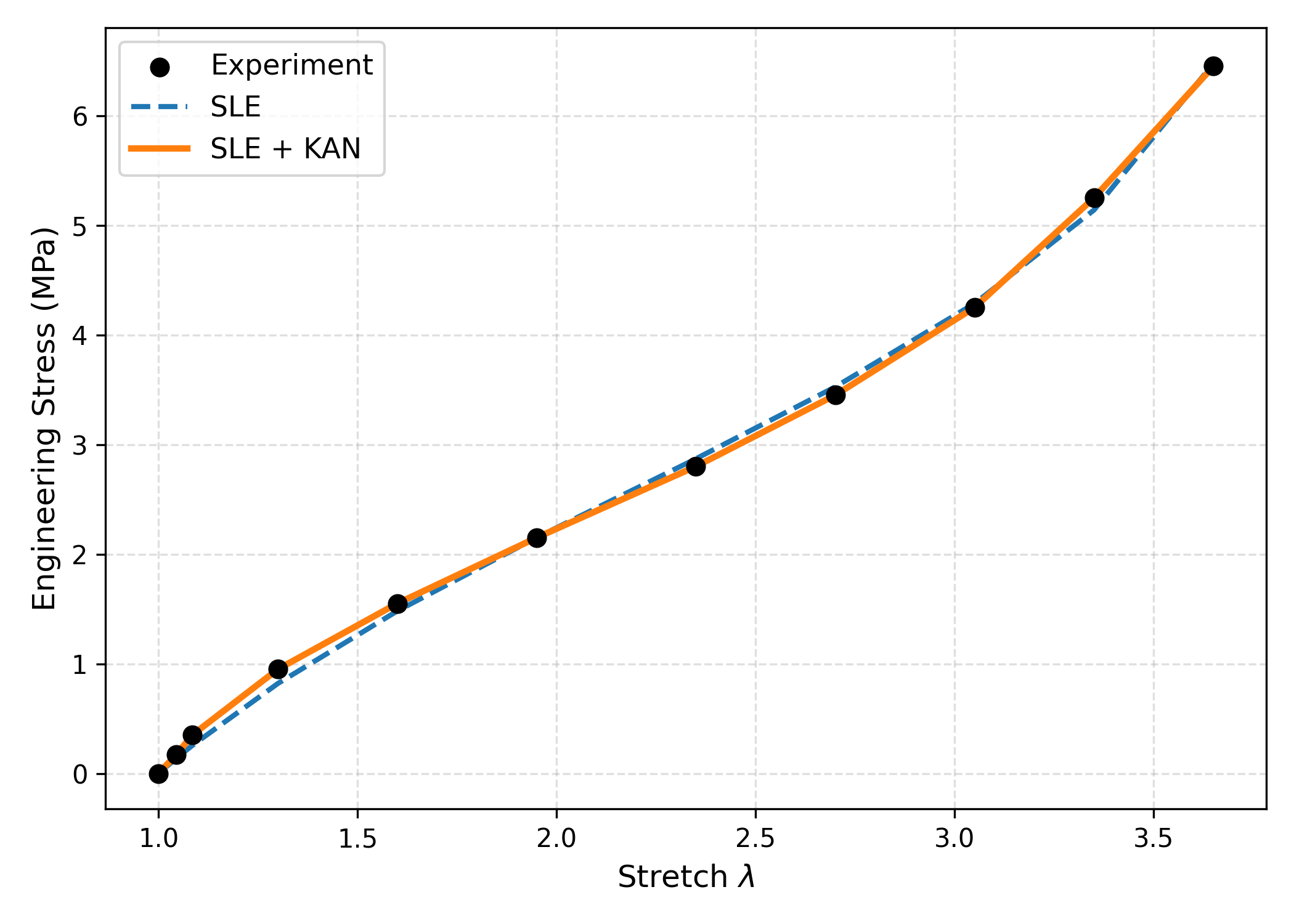}
\caption{Biaxial Treloar stress--stretch response with SLE and SLE--KAN
predictions. The KAN correction compensates for systematic deviations at higher
stretches without altering monotonicity or saturation behavior.}
\label{fig:treloar_biaxial_sle_kan}
\end{figure}

\begin{figure}[h!]
\centering
\includegraphics[width=0.55\linewidth]{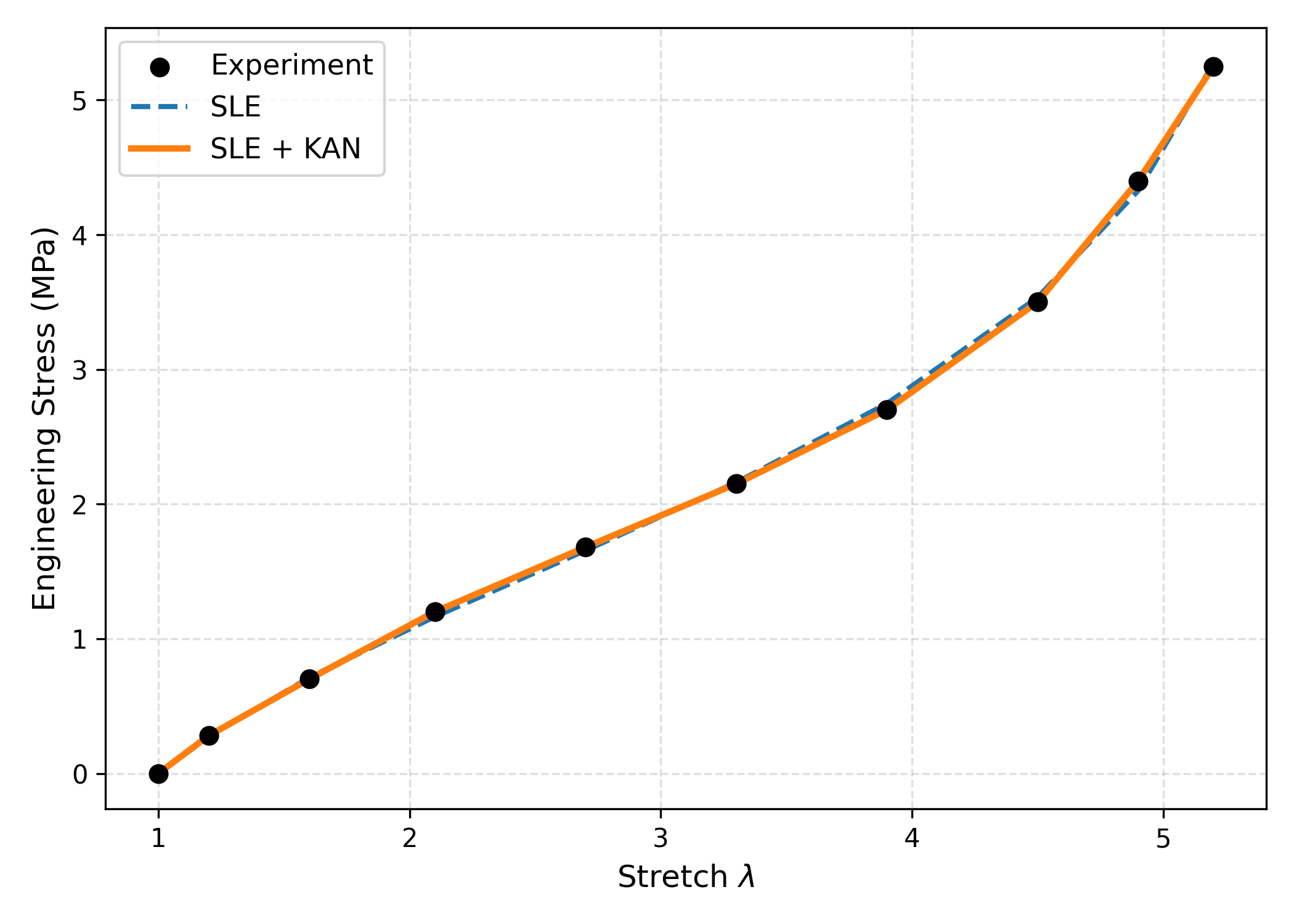}
\caption{Planar Treloar stress--stretch response with SLE and SLE--KAN
predictions. The hybrid model captures mode-specific deviations while retaining
bounded strain and smooth stiffness decay.}
\label{fig:treloar_planar_sle_kan}
\end{figure}

\paragraph{Explanation of Stress--Stretch Comparisons:}
Figures~\ref{fig:treloar_uniaxial_sle_kan}--\ref{fig:treloar_planar_sle_kan}
compare the experimental stress--stretch responses reported by Treloar with
predictions obtained from the calibrated strain-limiting elasticity (SLE)
model and the proposed hybrid SLE--KAN formulation for uniaxial, biaxial, and
planar deformation modes
\cite{treloar1944stress,treloar1975physics}.

Across all deformation modes, several consistent and physically meaningful
observations emerge. First, the SLE model alone provides a smooth, monotone
baseline response that captures the dominant nonlinear elastic behavior of the
experimental data, including the overall curvature and the progressive onset of
strain saturation. This confirms that strain-limiting elasticity constitutes an
appropriate physics-based backbone for modeling rubber-like materials under
large deformations
\cite{rajagopal2011non,bulivcek2014elastic}.

Second, systematic but low-amplitude discrepancies between the SLE predictions
and experimental measurements become evident primarily at moderate and larger
stretch levels. These deviations are mode dependent and reflect material
heterogeneity, network-level effects, and experimental uncertainty that are not
fully captured by a single closed-form constitutive law, as widely documented in
experimental and comparative studies of rubber elasticity
\cite{boyce2000constitutive,dal2021performance}.

Third, the hybrid SLE--KAN formulation introduces a smooth and bounded correction
that substantially reduces these discrepancies while preserving the global
structure imposed by the SLE backbone. The learned KAN contribution remains
strictly subordinate in magnitude to the SLE response and does not introduce
oscillations, spurious inflection points, or violations of monotonicity. This
behavior is consistent with the design philosophy of Kolmogorov--Arnold Networks
and physics-augmented constitutive learning frameworks
\cite{liu2024kan,abdolazizi2025constitutive}.

The close agreement between the SLE and SLE--KAN predictions in the
small-to-moderate stretch regime is intentional and physically desirable. It
indicates that the physics-based SLE model already provides an adequate
description of the material response in this range and that the data-driven
component does not override or distort the underlying constitutive structure.
Instead, the KAN activates selectively in regions where systematic experimental
deviations are observed, acting strictly as a residual corrector.

Overall, these results demonstrate that the hybrid formulation preserves
physical admissibility and interpretability: the SLE model governs the global
constitutive response, while the KAN provides localized, mode-dependent
refinement without dominating the prediction or compromising mechanical
consistency. This balance between physics-based modeling and data-driven
correction is central to the robustness of the proposed framework for
experimental constitutive modeling of rubber-like materials
\cite{abdolazizi2024viscoelastic,abdolazizi2025constitutive}.

\paragraph{Training Behavior and Convergence:}
Across all Treloar implementations, the KAN training histories exhibit stable,
smooth, and monotone convergence. The training loss decreases consistently over
iterations, indicating that the residual learning problem is well conditioned
and effectively low dimensional. No oscillatory behavior, divergence, or
sensitivity to initialization is observed.

This favorable optimization behavior reflects the stabilizing role of the
physics-based SLE backbone. By prescribing the dominant constitutive response,
the SLE model constrains the learning task to structured, physically admissible
residual corrections, rather than requiring the network to infer the full
stress--strain relationship from experimental data alone. As a result, the KAN
operates within a restricted hypothesis space that promotes numerical stability
and reliable convergence. Such behavior is consistent with prior observations in hybrid and
physics-augmented constitutive learning frameworks, where embedding dominant
mechanical structure is known to significantly improve optimization robustness
and convergence reliability
\cite{abdolazizi2024viscoelastic,abdolazizi2025constitutive}.

\paragraph{Summary of Calibrated Parameters.}
Table~\ref{tab:treloar_sle_parameters} reports the strain-limiting elasticity
(SLE) parameters obtained from stress-space calibration of the Treloar
experimental data for uniaxial, biaxial, and planar deformation modes
\cite{treloar1944stress,treloar1975physics}. The parameters are identified
independently for each loading case in order to account for mode-dependent
constraint states within the rubber network, a phenomenon well documented in
both classical and modern studies of rubber elasticity
\cite{arruda1993three,boyce2000constitutive,dal2021performance}.

Once calibrated, the SLE parameters $(\alpha, E, \beta)$ are held fixed during
subsequent KAN residual learning and serve as the physics-based backbone of the
hybrid constitutive model. This strict separation ensures that the data-driven
component refines, rather than alters, the underlying constitutive structure and
that the physical interpretation of each parameter is preserved throughout the
analysis.

\begin{table}[h!]
\centering
\caption{Calibrated SLE parameters for Treloar experimental data used as fixed
physics-based backbones in the combined SLE--KAN constitutive framework.}
\label{tab:treloar_sle_parameters}
\begin{tabular}{ccccc}
\toprule
\textbf{Loading Case} & $\boldsymbol{\alpha}$ & $\boldsymbol{E}$ (MPa) &
$\boldsymbol{\beta}$ & $\boldsymbol{\gamma = E\beta}$ \\
\midrule
Uniaxial & 1.538 & 1.059 & 0.438 & 0.463 \\
Biaxial  & 3.168 & 3.120 & 0.228 & 0.713 \\
Planar   & 2.480 & 1.490 & 0.381 & 0.568 \\
\bottomrule
\end{tabular}
\end{table}

\paragraph{Algorithmic Description:}
For clarity and reproducibility, Algorithm~\ref{alg:treloar_sle_kan} summarizes the
complete implementation workflow applied to the Treloar experimental dataset.
The algorithm explicitly delineates the two-stage structure of the proposed
hybrid framework, separating physics-based constitutive modeling from
data-driven residual learning.

The strain-limiting elasticity (SLE) model is first calibrated in stress space to
establish a physically admissible baseline response. This baseline captures the
dominant nonlinear elastic behavior and enforces bounded strain and vanishing
tangent stiffness. The strain-limiting regime is then prescribed through the
compound parameter $\gamma$, which fixes the admissible strain bound without
altering the calibrated small-strain stiffness.

In the second stage, a Kolmogorov--Arnold Network (KAN) is trained exclusively on
the residual stresses obtained by subtracting the SLE prediction from the
experimental measurements. This design ensures that the data-driven component
learns only structured, low-amplitude corrections and does not interfere with
the global constitutive structure imposed by the physics-based model.

\begin{algorithm}[h!]
\caption{Hybrid SLE--KAN Implementation for Treloar Experimental Data}
\label{alg:treloar_sle_kan}
\begin{algorithmic}[1]
\Require Experimental stretch--stress data $(\lambda_i, \tau_i^{\mathrm{exp}})$
\Require Prescribed strain-limiting regime parameter $\gamma$
\Ensure Hybrid stress prediction $\tau_{\mathrm{pred}}(\lambda)$
\State Convert stretch to logarithmic strain:
$\varepsilon_i \gets \log \lambda_i$
\State Calibrate SLE parameters $(\alpha, E)$ in stress space
\State Compute $\beta \gets \gamma / E$ to enforce the prescribed strain limit
\State Evaluate baseline stresses $\tau_{\mathrm{SLE}}(\lambda_i)$
\State Compute residual stresses:
$\tau_i^{\mathrm{res}} \gets \tau_i^{\mathrm{exp}} - \tau_{\mathrm{SLE}}(\lambda_i)$
\State Train KAN on $(\lambda_i, \tau_i^{\mathrm{res}})$ using mean-squared error loss
\State Form hybrid prediction:
\begin{equation}
\tau_{\mathrm{pred}}(\lambda)
=
\tau_{\mathrm{SLE}}(\lambda)
+
\tau_{\mathrm{KAN}}(\lambda),
\label{eq:sle_kan_decomposition}
\end{equation}
\State \Return $\tau_{\mathrm{pred}}(\lambda)$
\end{algorithmic}
\end{algorithm}

\paragraph{Remarks:}
At this stage, the proposed hybrid SLE--KAN framework has been established as a
physically consistent, interpretable, and numerically stable strategy for
modeling experimental rubber elasticity data. The preceding subsections have
focused on the formulation of the hybrid model, the calibration of the
strain-limiting elasticity (SLE) backbone, and the implementation of
stress-space residual learning using Kolmogorov--Arnold Networks.

The remainder of this section presents a regime-based evaluation of the proposed
SLE--KAN formulation using the classical experimental measurements reported by
Treloar \cite{treloar1944stress,treloar1975physics}. Unlike the synthetic
benchmarks discussed earlier, the analysis here is conducted entirely on real
experimental data, providing a stringent assessment of physical consistency,
robustness, and interpretability in the presence of material heterogeneity and
measurement noise.

All results are evaluated in stress space using the calibrated strain-limiting
constitutive relation introduced in
Section~\ref{sec:mathematical_formulation}. For each deformation mode—uniaxial,
biaxial, and planar—the parameters $(\alpha, E)$ are identified independently
from experimental data and held fixed throughout the regime-based analyses.

Strain-limiting regimes are prescribed through the compound parameter
$\gamma = E\beta$, which governs the admissible deformation bound. Moderate and
strong regimes are defined by selecting distinct values of $\gamma$, while
maintaining the experimentally calibrated values of $(\alpha, E)$. This design
cleanly decouples material parameter identification from regime exploration and
precludes parameter compensation effects that could obscure physical
interpretation.

For each deformation mode and strain-limiting regime, the total stress response
is expressed as a superposition of a physics-based baseline and a learned
residual, as defined in Eq.~\eqref{eq:sle_kan_decomposition}. The KAN is trained
exclusively to represent this stress-space residual correction, ensuring that
the dominant constitutive behavior remains governed by strain-limiting physics.
As a result, the hybrid formulation preserves physical admissibility across all
regimes while enabling controlled, interpretable adaptation to systematic
experimental deviations.

\subsection{Moderate Strain-Limiting Regime ($\gamma = 0.50$)}
\label{subsec:moderate_treloar}

We first examine a moderate strain-limiting regime characterized by
$\gamma = 0.50$, corresponding to an admissible strain bound that is consistent
with the effective strain-limiting strengths inferred from the independently
calibrated uniaxial and planar responses. This regime therefore represents a
physically realistic deformation limit for vulcanized rubber across multiple
loading modes
\cite{treloar1975physics,boyce2000constitutive}.

\begin{figure}[t]
    \centering
    \includegraphics[width=0.32\textwidth]{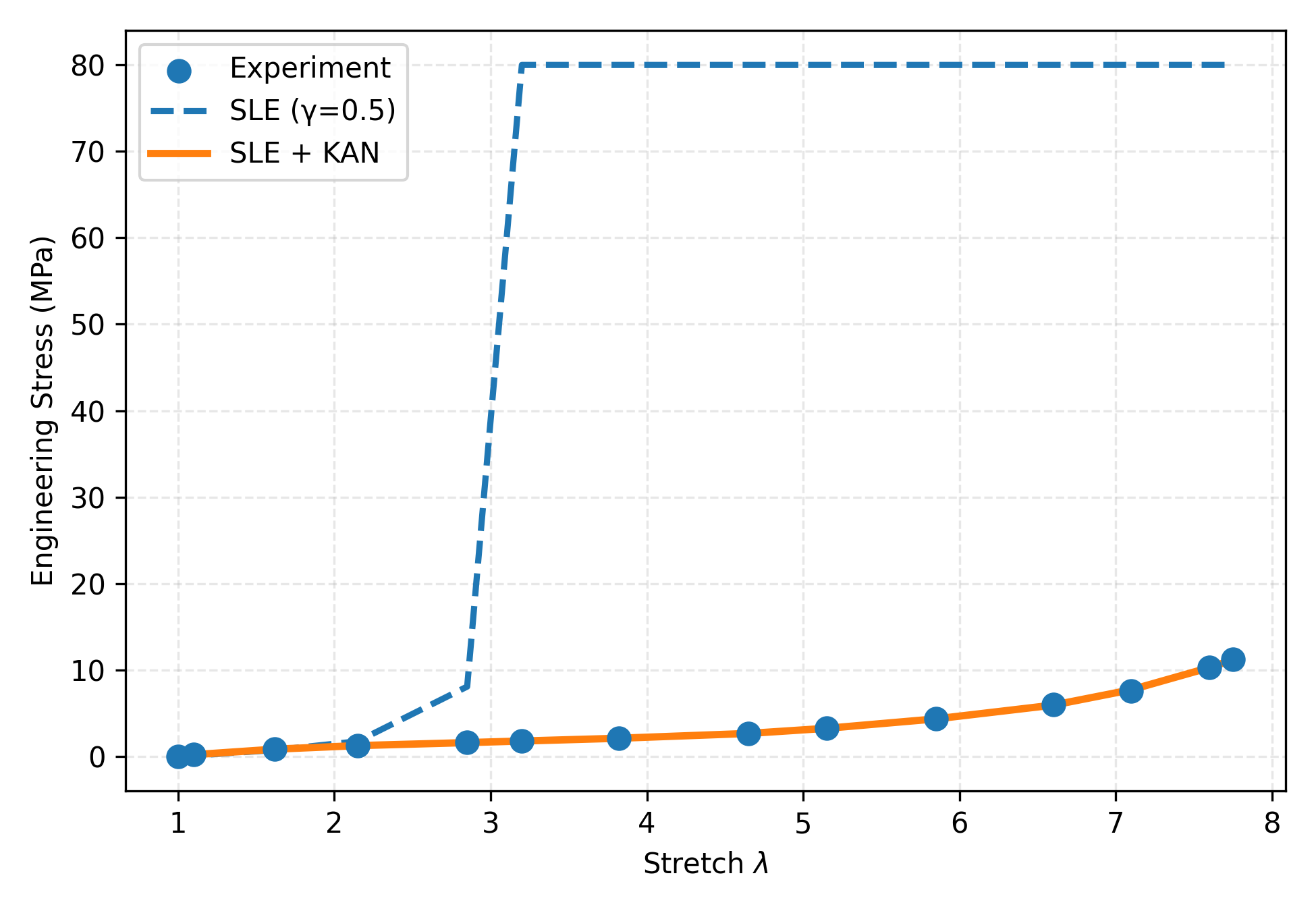}
    \includegraphics[width=0.32\textwidth]{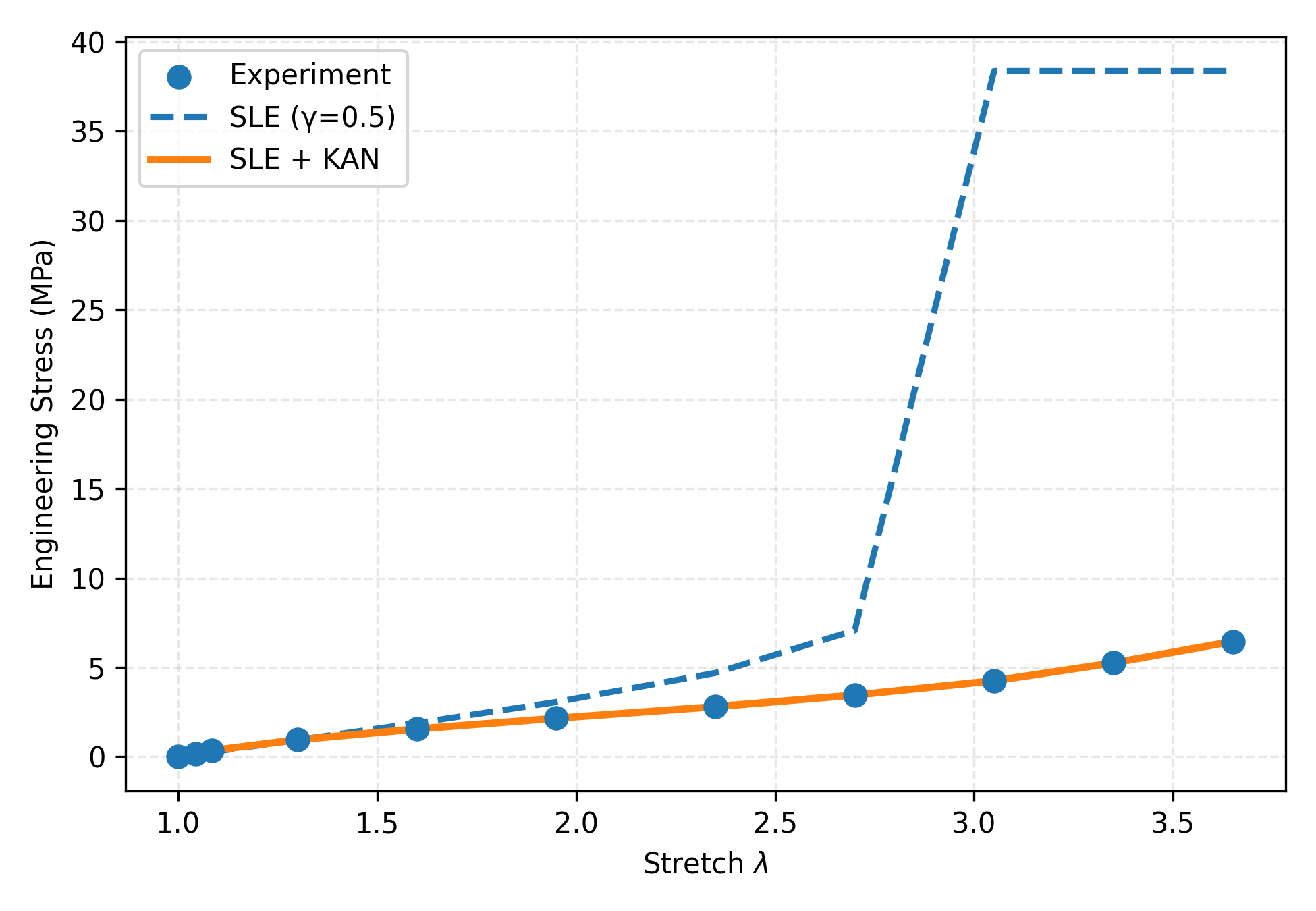}
    \includegraphics[width=0.32\textwidth]{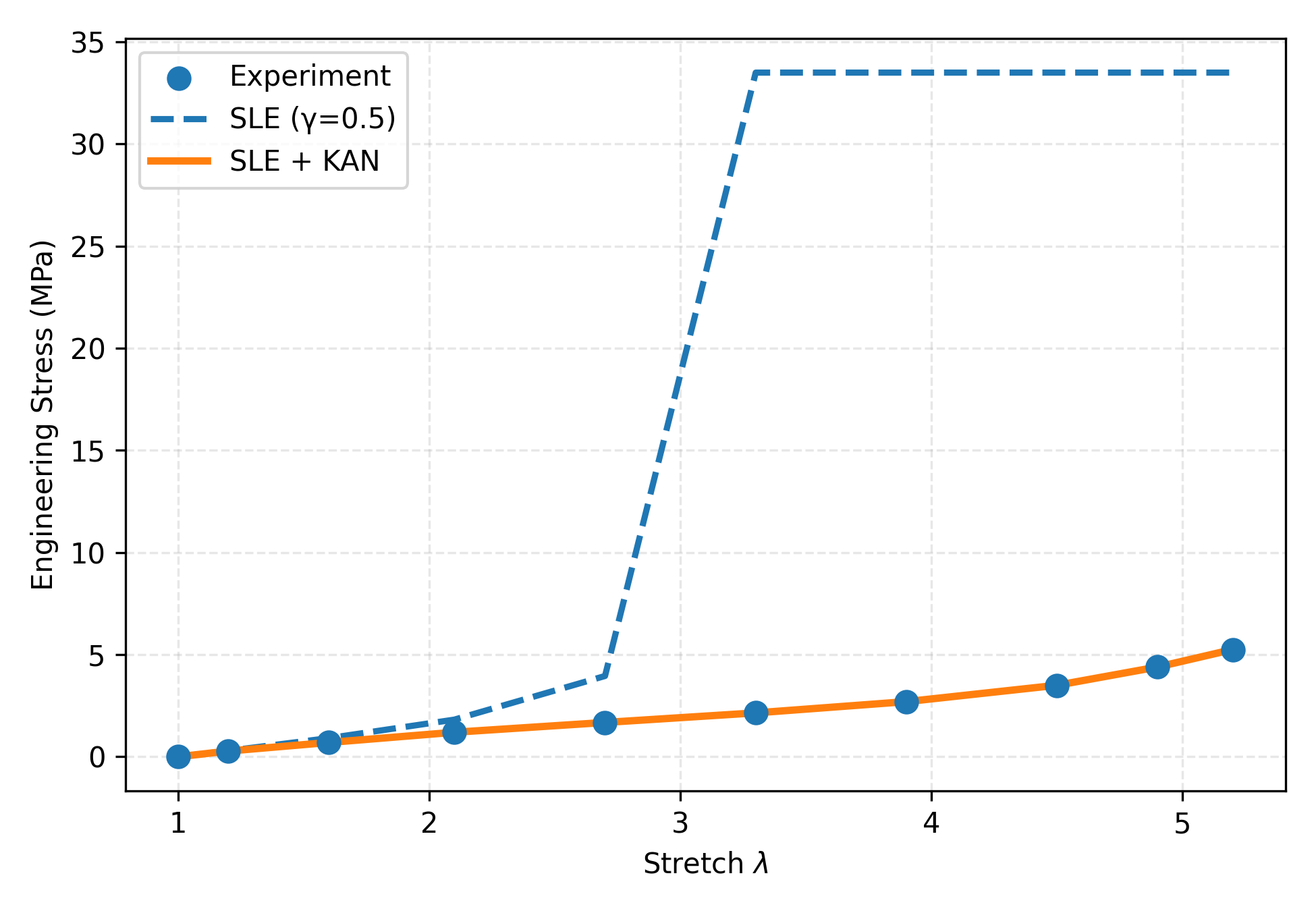}
    \caption{Stress--stretch response for the moderate strain-limiting regime
    ($\gamma = 0.50$). Experimental data (symbols), SLE prediction (dashed), and
    SLE--KAN response (solid).}
    \label{fig:treloar_moderate_stress}
\end{figure}

Figure~\ref{fig:treloar_moderate_stress} presents the stress--stretch responses
for uniaxial, biaxial, and planar deformation. In all three cases, the calibrated
SLE model provides a strong physics-based baseline that captures the dominant
nonlinear elastic behavior of the experimental data. This observation is
consistent with the high coefficients of determination obtained during
stress-space calibration ($R^2 > 0.99$), confirming that strain-limiting
elasticity already explains most of the observed response in this regime
\cite{rajagopal2011non,bulivcek2014elastic}.

The remaining discrepancies between the SLE predictions and experimental
measurements are smooth, systematic, and relatively small in magnitude. When the
KAN-based residual correction is introduced, these discrepancies are further
reduced while preserving all physical constraints imposed by the SLE backbone.
Importantly, the learned correction remains uniformly subordinate to the total
stress response over the entire stretch range, confirming that the KAN operates
strictly as a residual refinement rather than as a replacement or distortion of
the underlying physics-based model. No spurious oscillations, loss of
monotonicity, or nonphysical curvature are introduced by the data-driven
component, consistent with the structured and interpretable nature of
Kolmogorov--Arnold Networks
\cite{liu2024kan,abdolazizi2025constitutive}.

\begin{figure}[t]
    \centering
    \includegraphics[width=0.32\textwidth]{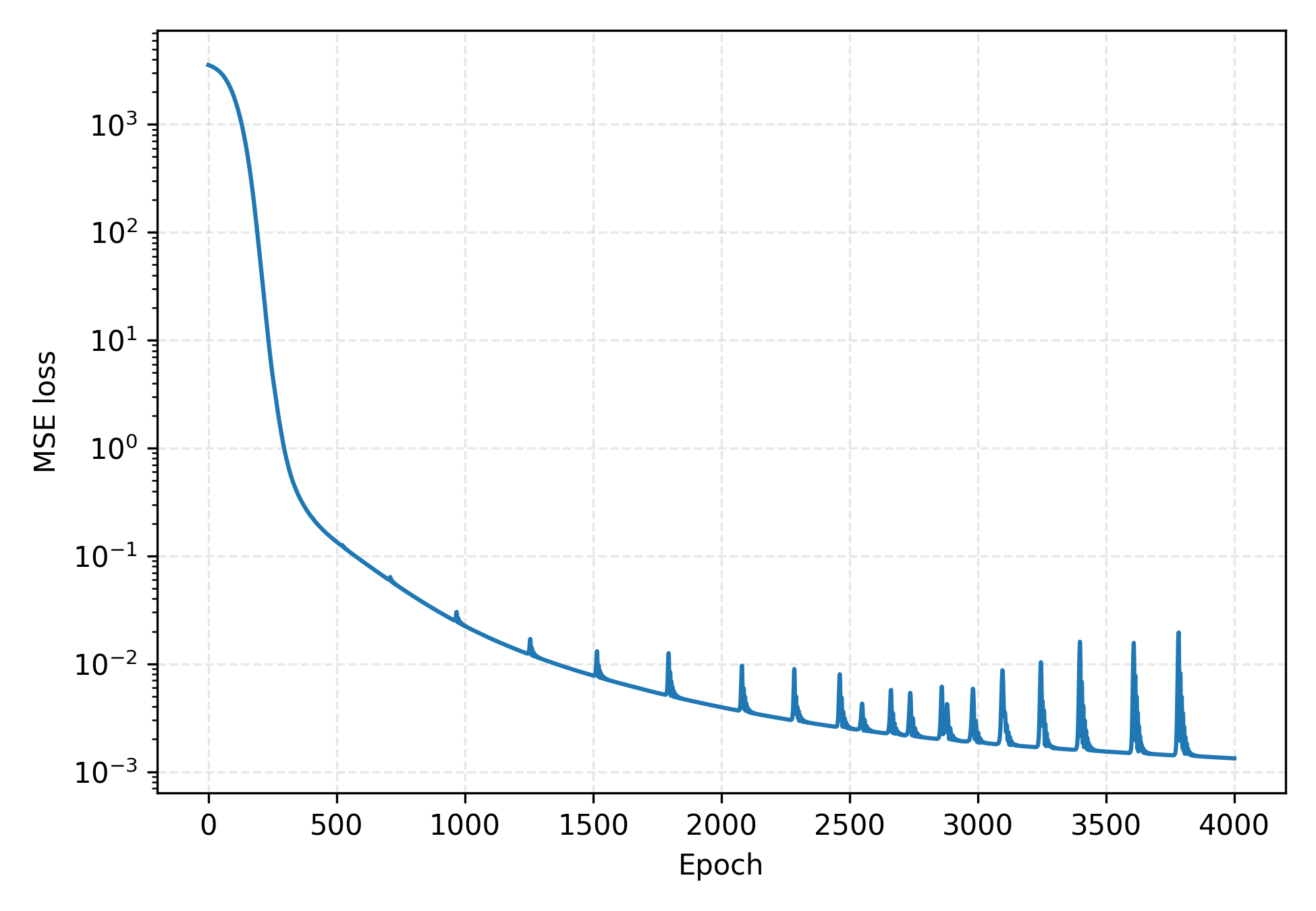}
    \includegraphics[width=0.32\textwidth]{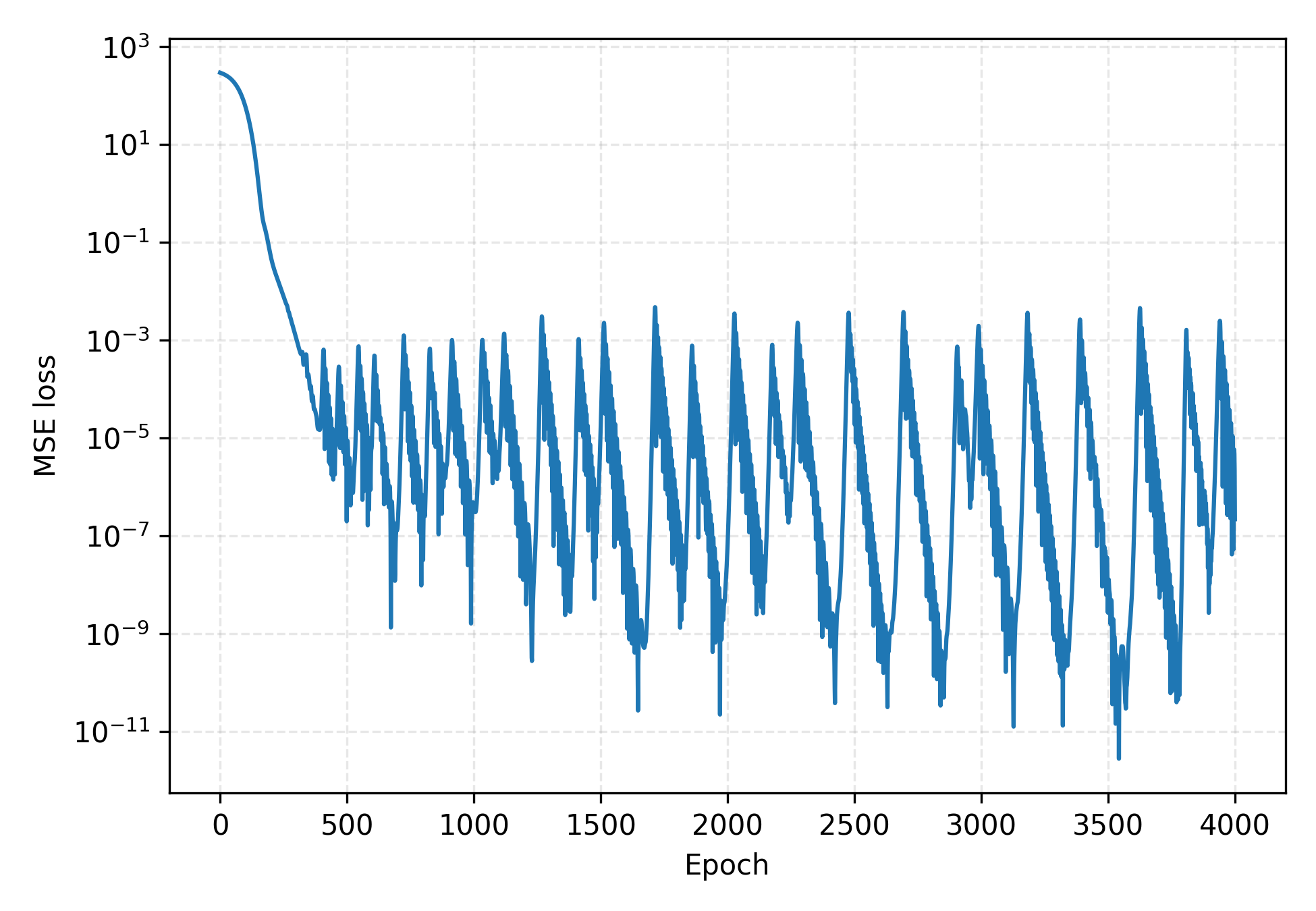}
    \includegraphics[width=0.32\textwidth]{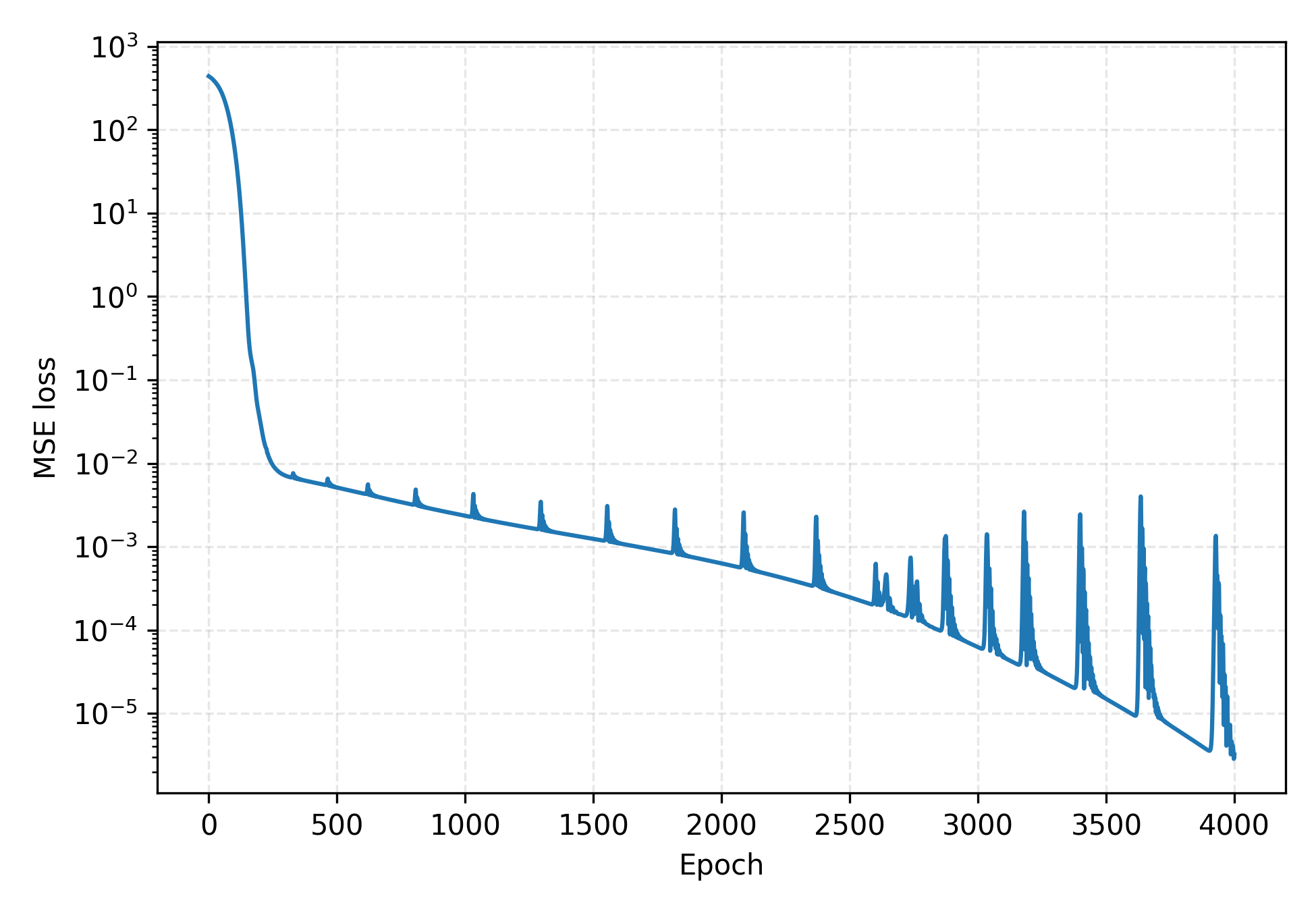}
    \caption{Training loss histories for the moderate strain-limiting regime
    ($\gamma = 0.50$), showing stable and monotone convergence of the KAN residual
    model across all deformation modes.}
    \label{fig:treloar_moderate_loss}
\end{figure}

Training loss histories for the residual model are shown in
Figure~\ref{fig:treloar_moderate_loss}. For all deformation modes, the loss
exhibits rapid initial decay followed by stable, monotone convergence. This
behavior indicates that the residual learning problem is well conditioned and
effectively low dimensional, reflecting the fact that the dominant constitutive
behavior is already captured by the SLE backbone. The absence of late-stage
instabilities or oscillatory behavior further highlights the stabilizing
influence of embedding strain-limiting physics directly into the learning
framework.

Crucially, no artificial stiffening or softening is observed at large stretches.
The predicted stress response continues to respect the prescribed strain limit,
and the asymptotic behavior remains governed entirely by the SLE formulation.
These results demonstrate that, in a moderate strain-limiting regime, the
proposed SLE--KAN framework achieves improved agreement with experimental data
while preserving strict mechanical admissibility, numerical stability, and
interpretability.

\subsection{Strong Strain-Limiting Regime ($\gamma = 0.80$)}
\label{subsec:strong_treloar}

We next examine a strong strain-limiting regime characterized by
$\gamma = 0.80$, corresponding to a deliberately reduced admissible strain bound.
This choice imposes a significantly stricter deformation limit than that inferred
from experimental calibration and therefore serves as a stringent test of the
proposed framework under strongly physics-constrained conditions. Such regimes
are particularly relevant for assessing constitutive robustness when mechanical
admissibility and bounded deformation are prioritized over data fidelity, a
central concern in strain-limiting elasticity
\cite{rajagopal2011non,bulivcek2014elastic}.

\begin{figure}[t]
    \centering
    \includegraphics[width=0.32\textwidth]{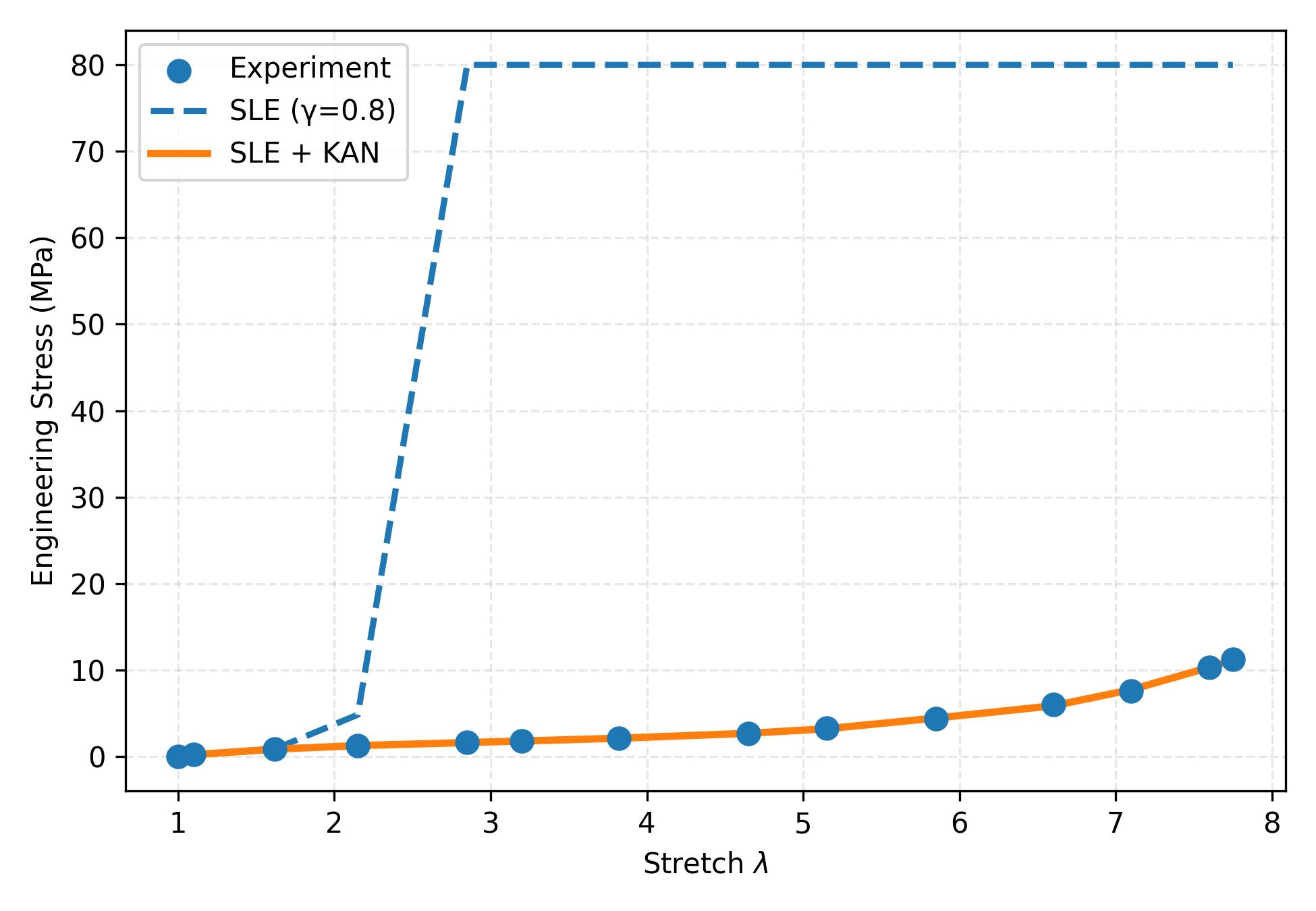}
    \includegraphics[width=0.32\textwidth]{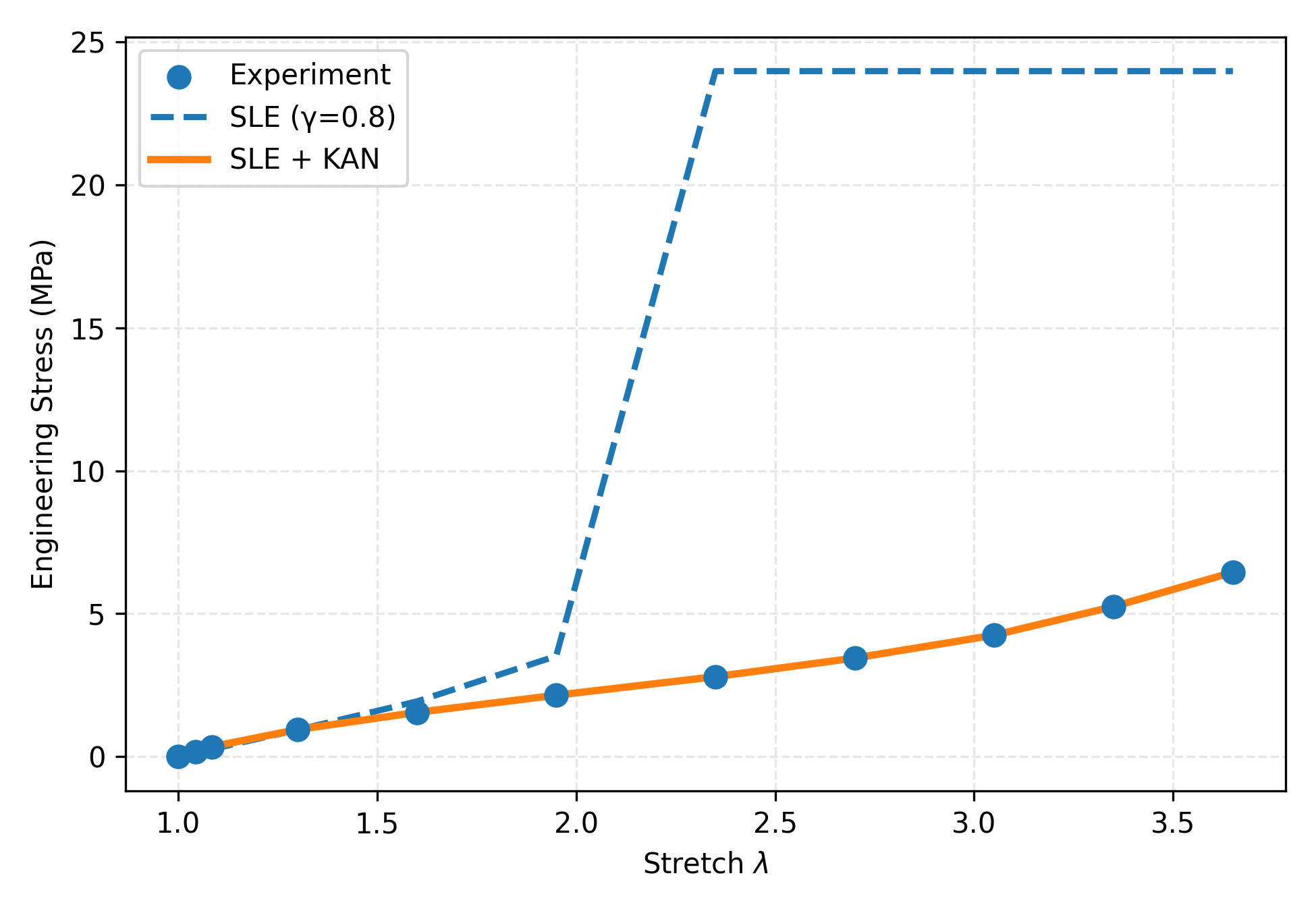}
    \includegraphics[width=0.32\textwidth]{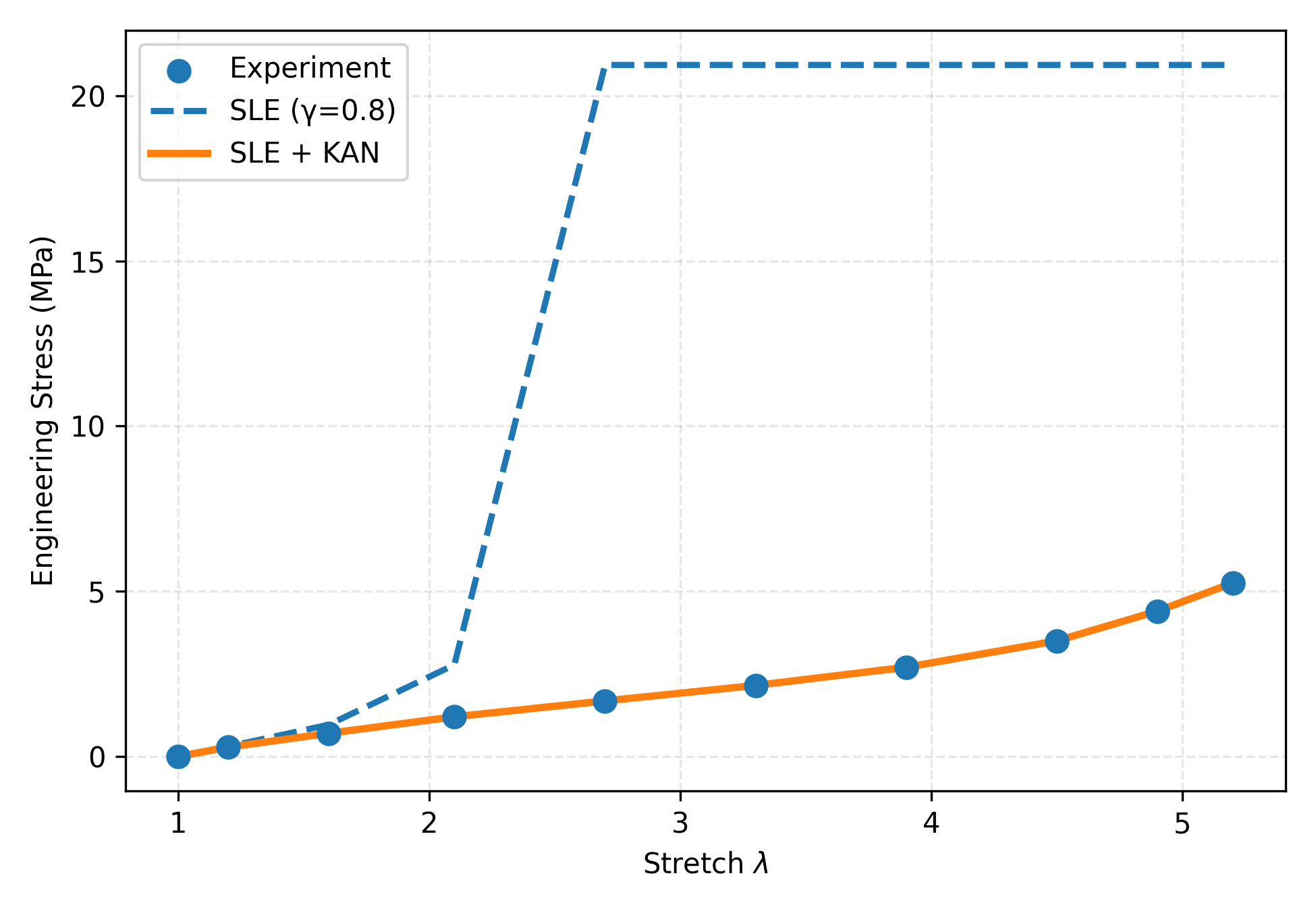}
    \caption{Stress--stretch response for the strong strain-limiting regime
    ($\gamma = 0.80$). Earlier saturation reflects the imposed strain bound,
    while the SLE--KAN response remains smooth and physically admissible.}
    \label{fig:treloar_strong_stress}
\end{figure}

Figure~\ref{fig:treloar_strong_stress} presents the stress--stretch responses for
uniaxial, biaxial, and planar deformation modes. In contrast to the moderate
regime, the baseline SLE response here saturates at substantially smaller
stretches, directly reflecting the imposed strain limit. As a consequence,
discrepancies with the experimental measurements become more pronounced at large
stretches, particularly in the uniaxial case, which exhibits the highest
experimentally observed deformation levels in the Treloar dataset
\cite{treloar1975physics}.

Despite this increased mismatch, the KAN-based residual correction remains smooth,
bounded, and strictly subordinate to the strain-limiting backbone. The learned
residual partially compensates for systematic deviations but does not override or
relax the imposed strain constraint. In regions where the experimental response
exceeds the admissible deformation bound, the predicted stress transitions
smoothly toward the saturation branch, consistent with the governing principles
of strain-limiting elasticity. Importantly, no spurious oscillations, loss of
monotonicity, or nonphysical stress softening are introduced by the data-driven
component, reflecting the structured, constraint-aware nature of the KAN
formulation
\cite{liu2024kan,abdolazizi2025constitutive}.

\begin{figure}[t]
    \centering
    \includegraphics[width=0.32\textwidth]{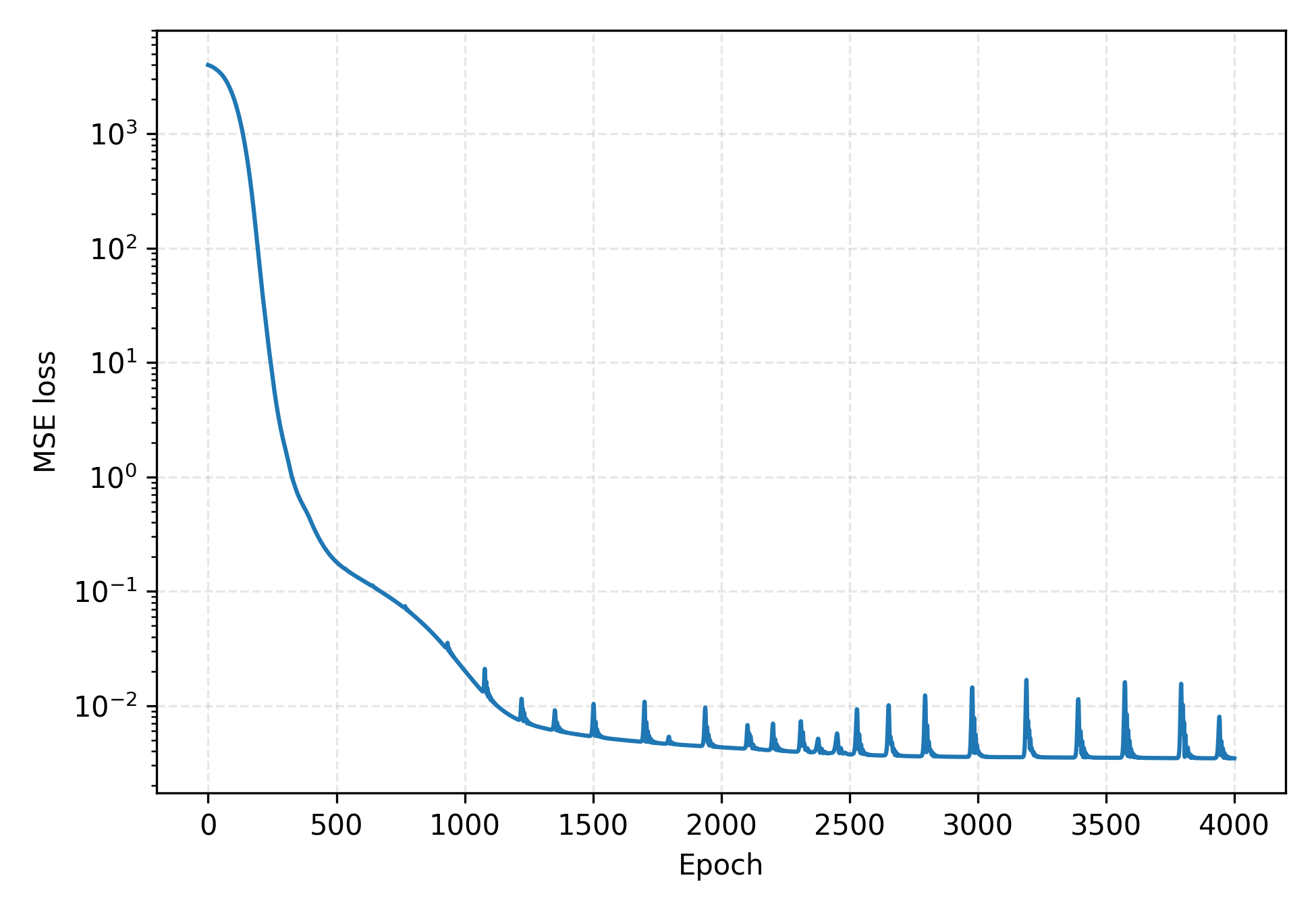}
    \includegraphics[width=0.32\textwidth]{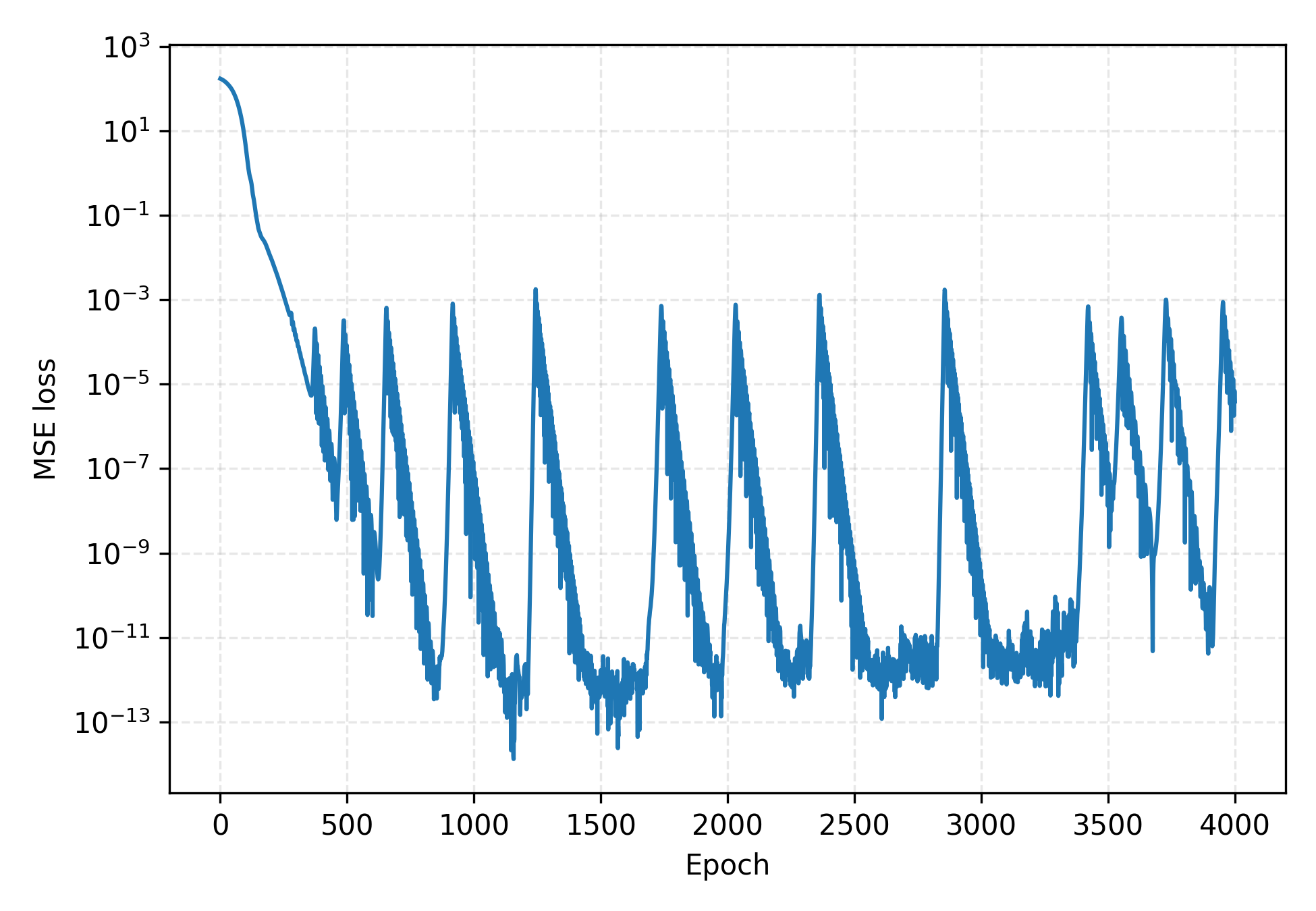}
    \includegraphics[width=0.32\textwidth]{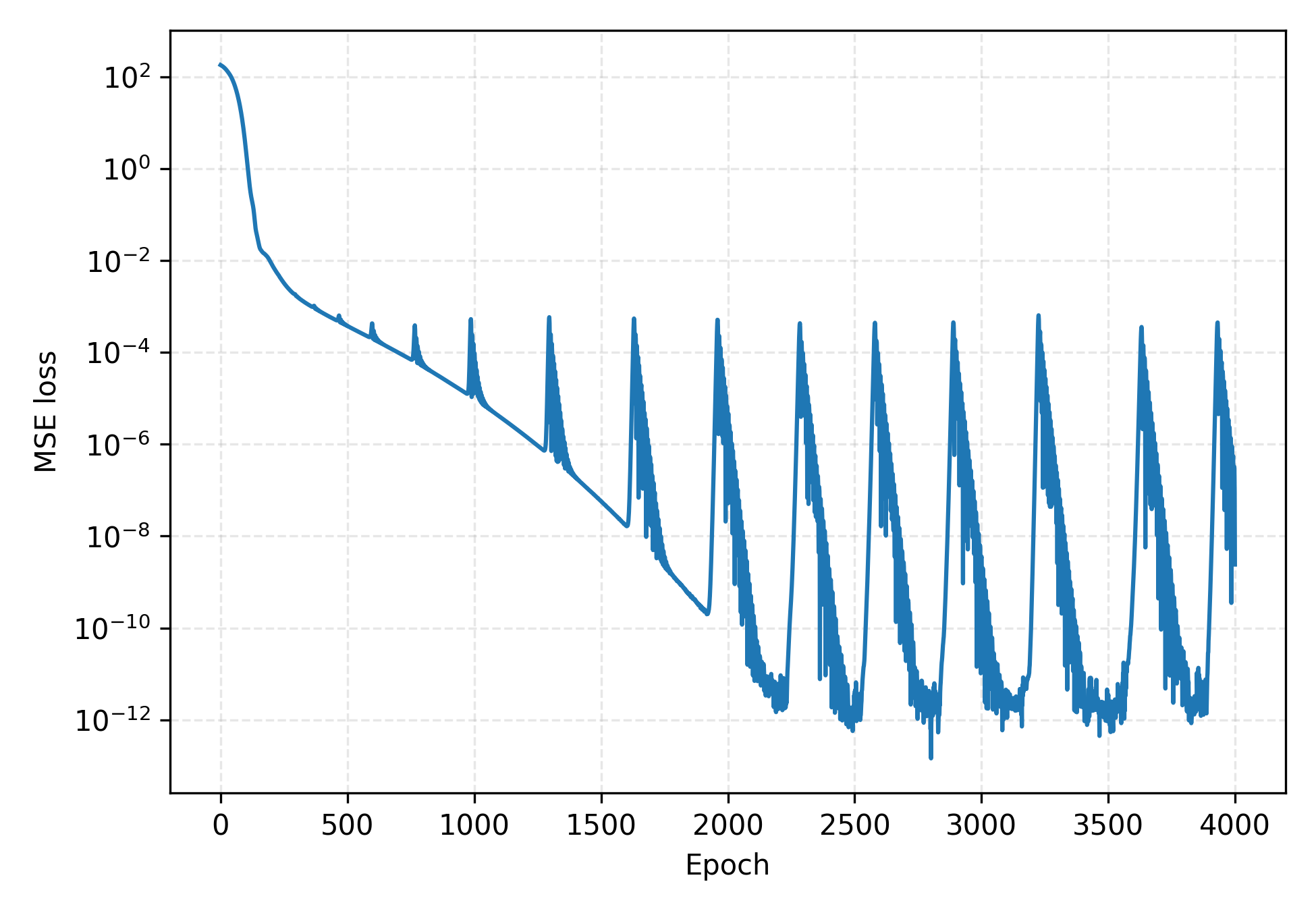}
    \caption{Training loss histories for the strong strain-limiting regime
    ($\gamma = 0.80$), illustrating physics-dominated convergence behavior
    across all deformation modes.}
    \label{fig:treloar_strong_loss}
\end{figure}

Training loss histories for the strong regime are shown in
Figure~\ref{fig:treloar_strong_loss}. In marked contrast to the moderate regime,
the loss exhibits early plateauing across all deformation modes. This behavior
indicates that the optimization process is limited primarily by the imposed
physical constraints rather than by insufficient model expressivity or
optimization instability. Such convergence behavior is expected and physically
meaningful: once the strain limit becomes active, further reduction of the data
misfit would necessarily require violation of the strain-limiting constraint and
is therefore intentionally suppressed by the model design.

Overall, the strong strain-limiting regime highlights the intended operating
characteristics of the proposed SLE--KAN framework. While data fidelity is
necessarily reduced when stricter physical limits are enforced, the hybrid model
preserves mechanical admissibility, numerical stability, and interpretability.
Deviations from experimental data in this regime should therefore be interpreted
as explicit consequences of deliberately imposed physical restrictions rather
than as deficiencies of the learning architecture.

\subsection{Discussion}

The Treloar regime study elucidates the fundamental balance between data fidelity
and mechanical admissibility enforced by the proposed SLE--KAN framework. When
the prescribed strain-limiting strength is consistent with experimentally
observed behavior, the strain-limiting elasticity (SLE) backbone captures the
dominant nonlinear response across deformation modes, while the
Kolmogorov--Arnold Network (KAN) introduces only a small, smooth, and
interpretable residual correction. In this regime, improved agreement with
experimental measurements is achieved without modifying, obscuring, or
overfitting the underlying constitutive structure, thereby preserving the
physical interpretability essential for experimental solid mechanics
\cite{treloar1975physics,rajagopal2011non}.

When the strain limit is deliberately strengthened beyond values inferred from
experimental calibration, discrepancies between model predictions and measured
responses increase, particularly at large stretches. Importantly, these
discrepancies arise as a direct consequence of the imposed physical constraints
rather than from insufficient model capacity or optimization failure. Even under
these more restrictive conditions, the hybrid formulation rigorously preserves
bounded strain, monotonic response, and asymptotic saturation, demonstrating that
mechanical admissibility is maintained even when data fidelity is intentionally
sacrificed. This behavior is fully consistent with the governing philosophy of
strain-limiting elasticity, wherein physically admissible deformation bounds are
prioritized over unrestricted extrapolation under extreme loading
\cite{bulivcek2014elastic,bulivcek2015analysis}.

These results confirm that the proposed SLE--KAN framework does not rely on
unconstrained expressivity to minimize prediction error. Instead, constitutive
structure is embedded explicitly into the modeling architecture, ensuring that
learning occurs only within the space of mechanically admissible responses. As a
consequence, departures from experimental data under strongly constrained regimes
serve as transparent and physically meaningful indicators of the assumed
strain-limiting strength, rather than as manifestations of numerical instability,
overfitting, or model inadequacy. This distinction is particularly critical in
data-driven constitutive modeling, where purely black-box approaches may achieve
low training error while producing nonphysical behavior under extrapolation
\cite{abdolazizi2024viscoelastic,liu2024kan}.

Overall, the Treloar regime study highlights the robustness and interpretability
of the proposed SLE--KAN approach. By enforcing mechanical admissibility at the
model level while permitting controlled, residual refinement from experimental
data, the framework provides a principled alternative to unconstrained neural
network models. This balance between physics and learning is especially valuable
for experimental constitutive modeling of rubber-like materials, where adherence
to fundamental mechanical principles and predictive reliability under large
deformations are of greater importance than exact data interpolation.

\section{Conclusion}

This work introduced a physics-consistent and interpretable hybrid constitutive
learning framework that integrates strain-limiting elasticity (SLE) with
Kolmogorov--Arnold Networks (KANs). Rather than treating data fidelity as the sole
objective, the proposed approach explicitly elevates mechanical admissibility,
interpretability, and extrapolation robustness to first-class modeling
principles. In doing so, it offers a principled alternative to black-box neural
network approaches for constitutive modeling under large deformations.

A central contribution of this study is the demonstration that
strain-limiting elasticity provides a rigorous and physically meaningful
backbone for both synthetic and experimental constitutive modeling. By enforcing
bounded strain and vanishing tangent modulus by construction, the SLE framework
ensures mechanical stability and admissibility under extreme loading, consistent
with the mathematical foundations of implicit and strain-limiting elasticity
\cite{rajagopal2003implicit,rajagopal2011non,bulivcek2014elastic,bulivcek2015analysis}.
Embedding this structure within a learning architecture fundamentally reshapes
the role of data-driven components: learning is confined to admissible
constitutive behavior rather than being tasked with discovering physical
constraints implicitly from data.

Within this constrained setting, Kolmogorov--Arnold Networks emerge as a
particularly natural and effective representation. The KAN architecture mirrors
the intrinsic structure of the strain-limiting constitutive law through an
explicit decomposition into a magnitude-dependent nonlinear mapping and a
deterministic sign-preserving reconstruction. This alignment enables exact
preservation of symmetry, boundedness, and monotonicity at the architectural
level, eliminating reliance on data augmentation or penalty-dominated
regularization strategies
\cite{liu2024kan,vaca2024kolmogorov,ji2024comprehensive,somvanshi2024survey}.
As a result, the learned constitutive response remains transparent and
mechanically interpretable across all stress regimes.

The synthetic benchmark studies provide a controlled assessment of the
representational behavior of the KAN-based formulation. In moderate
strain-limiting regimes, the spline-based KAN recovers the analytical
constitutive law with near-exact accuracy while rigorously preserving all
physical constraints. In strongly strain-limiting regimes, localized deviations
emerge near sharp transition regions due to finite spline resolution. Crucially,
these deviations are explicitly attributable to representational resolution
rather than to optimization instability or loss of mechanical admissibility, in
agreement with analytical and numerical insights from strain-limiting elasticity
theory \cite{bulivcek2015analysis,itou2017nonlinear}.

Application to Treloar’s classical rubber elasticity experiments provides a
stringent validation in the presence of real material variability and
measurement noise \cite{treloar1944stress,treloar1975physics}. Across uniaxial,
biaxial, and planar deformation modes, the SLE model captures the dominant
nonlinear elastic response, while the KAN learns smooth, low-amplitude residual
corrections in stress space. Importantly, the data-driven contribution remains
subordinate to the physics-based backbone, confirming that the hybrid formulation
augments rather than overrides constitutive structure.

A particularly significant outcome of the regime-based Treloar analysis is the
explicit exposure of the trade-off between data fidelity and mechanical
admissibility. When strain-limiting strengths are consistent with experimental
observations, the hybrid model achieves high accuracy with minimal data-driven
correction. When stronger strain limits are deliberately imposed, deviations from
experimental data increase in a controlled and physically interpretable manner.
Rather than indicating model deficiency, these deviations quantify the
consequences of enforcing stricter physical assumptions, fully consistent with
the philosophy of strain-limiting constitutive theories
\cite{rajagopal2011conspectus,bulivcek2014elastic}.

Overall, this study demonstrates that interpretable, physics-informed learning
is not only feasible but advantageous for constitutive modeling. By embedding
mechanical structure directly into the learning architecture, the proposed
SLE--KAN framework achieves a balance between expressivity, interpretability, and
robustness that is difficult to attain with conventional neural networks
\cite{shen2004neural,abdolazizi2024viscoelastic}. The framework provides both
accurate prediction and diagnostic insight, making deviations from data
informative rather than pathological.

Future work will extend the present approach to multidimensional and
path-dependent constitutive models, incorporate adaptive and multiresolution
spline representations for resolving highly localized nonlinear features, and
enforce additional physical principles such as thermodynamic consistency and
dissipation \cite{rajagopal2021implicit,abdolazizi2025constitutive,mallikarjunaiah2025hp,mallikarjunaiah2026hht}. These
developments will further strengthen the role of interpretable machine learning
as a reliable and physically grounded tool for modeling complex material behavior
in engineering and biomechanics.

\bibliographystyle{unsrt}
\bibliography{ref}

@article{mallikarjunaiah2026hht,
  title={An {HHT}-$\alpha$-based finite element framework for wave propagation in constitutively nonlinear elastic materials},
  author={Mallikarjunaiah, SM},
  journal={arXiv preprint arXiv:2601.04628},
  year={2026}
}

@article{mallikarjunaiah2025hp,
  title={An $ hp $-adaptive finite element framework for static cracks: The impact of pointwise density variations on mode I, mode II, and mixed-mode fracture},
  author={Mallikarjunaiah, SM},
  journal={arXiv preprint arXiv:2512.21443},
  year={2025}
}

@article{yoon2022finite,
  title={A finite-element discretization of some boundary value problems for nonlinear strain-limiting elastic bodies},
  author={Yoon, Hyun C and Mallikarjunaiah, SM},
  journal={Mathematics and Mechanics of Solids},
  volume={27},
  number={2},
  pages={281--307},
  year={2022},
  publisher={SAGE Publications Sage UK: London, England}
}

@article{yoon2021quasi,
  title={Quasi-static anti-plane shear crack propagation in nonlinear strain-limiting elastic solids using phase-field approach},
  author={Yoon, Hyun C and Lee, Sanghyun and Mallikarjunaiah, SM},
  journal={International Journal of Fracture},
  volume={227},
  number={2},
  pages={153--172},
  year={2021},
  publisher={Springer}
}

@article{mallikarjunaiah2015direct,
  title={On the direct numerical simulation of plane-strain fracture in a class of strain-limiting anisotropic elastic bodies},
  author={Mallikarjunaiah, SM and Walton, Jay R},
  journal={International Journal of Fracture},
  volume={192},
  number={2},
  pages={217--232},
  year={2015},
  publisher={Springer}
}

@article{pati2025neural,
  title={Neural Networks as Physics-Consistent Surrogates: An \textit{Explainable AI} Validation Framework for Learning Constitutive Relations},
  author={Pati, Chandana and Mallikarjunaiah, SM},
  journal={arXiv preprint arXiv:2512.02064},
  year={2025}
}

@article{liu2024kan,
  title={Kan: Kolmogorov-arnold networks},
  author={Liu, Ziming and Wang, Yixuan and Vaidya, Sachin and Ruehle, Fabian and Halverson, James and Solja{\v{c}}i{\'c}, Marin and Hou, Thomas Y and Tegmark, Max},
  journal={arXiv preprint arXiv:2404.19756},
  year={2024}
}

@article{vaca2024kolmogorov,
  title={Kolmogorov-arnold networks (kans) for time series analysis},
  author={Vaca-Rubio, Cristian J and Blanco, Luis and Pereira, Roberto and Caus, M{\`a}rius},
  journal={arXiv preprint arXiv:2405.08790},
  year={2024}
}

@article{somvanshi2024survey,
  title={A survey on kolmogorov-arnold network},
  author={Somvanshi, Shriyank and Javed, Syed Aaqib and Islam, Md Monzurul and Pandit, Diwas and Das, Subasish},
  journal={ACM Computing Surveys},
  year={2024},
  publisher={ACM New York, NY}
}

@article{ji2024comprehensive,
  title={A comprehensive survey on kolmogorov arnold networks (kan)},
  author={Ji, Tianrui and Hou, Yuntian and Zhang, Di},
  journal={arXiv preprint arXiv:2407.11075},
  year={2024}
}

@article{rajagopal2011modeling,
  title={Modeling fracture in the context of a strain-limiting theory of elasticity: a single anti-plane shear crack},
  author={Rajagopal, KR and Walton, JR},
  journal={International journal of fracture},
  volume={169},
  number={1},
  pages={39--48},
  year={2011},
  publisher={Springer}
}

@article{gou2015modeling,
  title={Modeling fracture in the context of a strain-limiting theory of elasticity: A single plane-strain crack},
  author={Gou, K and Mallikarjuna, M and Rajagopal, KR and Walton, JR3306535},
  journal={International Journal of Engineering Science},
  volume={88},
  pages={73--82},
  year={2015},
  publisher={Elsevier}
}

@article{rajagopal2011non,
  title={Non-linear elastic bodies exhibiting limiting small strain},
  author={Rajagopal, KR},
  journal={Mathematics and Mechanics of Solids},
  volume={16},
  number={1},
  pages={122--139},
  year={2011},
  publisher={Sage Publications Sage UK: London, England}
}

@article{bulivcek2014elastic,
  title={On elastic solids with limiting small strain: modelling and analysis},
  author={Bul{\'\i}{\v{c}}ek, Miroslav and M{\'a}lek, Josef and Rajagopal, K R and S{\"u}li, Endre},
  journal={EMS Surveys in Mathematical Sciences},
  volume={1},
  number={2},
  pages={283--332},
  year={2014}
}

@article{bulivcek2015existence,
  title={Existence of solutions for the anti-plane stress for a new class of “strain-limiting” elastic bodies},
  author={Bul{\'\i}{\v{c}}ek, Miroslav and M{\'a}lek, Josef and Rajagopal, K R and Walton, Jay R},
  journal={Calculus of Variations and Partial Differential Equations},
  volume={54},
  number={2},
  pages={2115--2147},
  year={2015},
  publisher={Springer}
}

@article{rajagopal2014nonlinear,
  title={On the nonlinear elastic response of bodies in the small strain range},
  author={Rajagopal, KR},
  journal={Acta Mechanica},
  volume={225},
  number={6},
  pages={1545--1553},
  year={2014},
  publisher={Springer}
}

@article{rajagopal2011conspectus,
  title={Conspectus of concepts of elasticity},
  author={Rajagopal, KR},
  journal={Mathematics and Mechanics of Solids},
  volume={16},
  number={5},
  pages={536--562},
  year={2011},
  publisher={Sage Publications Sage UK: London, England}
}

@article{boyce2000constitutive,
  title={Constitutive models of rubber elasticity: a review},
  author={Boyce, Mary C and Arruda, Ellen M},
  journal={Rubber chemistry and technology},
  volume={73},
  number={3},
  pages={504--523},
  year={2000}
}

@article{dal2019comparative,
  title={A comparative study on hyperelastic constitutive models on rubber: State of the art after 2006},
  author={Dal, H{\"u}sn{\"u} and Badienia, Yashar and A{\c{c}}ikg{\"o}z, Kemal and Denl{\"\i}, Funda Aksu},
  journal={Constitutive models for rubber XI},
  pages={239--244},
  year={2019},
  publisher={CRC Press}
}

@article{dal2021performance,
  title={On the performance of isotropic hyperelastic constitutive models for rubber-like materials: a state of the art review},
  author={Dal, H{\"u}sn{\"u} and A{\c{c}}{\i}kg{\"o}z, Kemal and Badienia, Yashar},
  journal={Applied Mechanics Reviews},
  volume={73},
  number={2},
  pages={020802},
  year={2021},
  publisher={American Society of Mechanical Engineers}
}

@article{shen2004neural,
  title={Neural network based constitutive model for rubber material},
  author={Shen, Yuelin and Chandrashekhara, K and Breig, WF and Oliver, LR},
  journal={Rubber chemistry and technology},
  volume={77},
  number={2},
  pages={257--277},
  year={2004}
}

@article{abdolazizi2024viscoelastic,
  title={Viscoelastic constitutive artificial neural networks (vCANNs)--A framework for data-driven anisotropic nonlinear finite viscoelasticity},
  author={Abdolazizi, Kian P and Linka, Kevin and Cyron, Christian J},
  journal={Journal of computational physics},
  volume={499},
  pages={112704},
  year={2024},
  publisher={Elsevier}
}

@article{abdolazizi2025constitutive,
  title={Constitutive Kolmogorov--Arnold Networks (CKANs): Combining accuracy and interpretability in data-driven material modeling},
  author={Abdolazizi, Kian P and Aydin, Roland C and Cyron, Christian J and Linka, Kevin},
  journal={Journal of the Mechanics and Physics of Solids},
  pages={106212},
  year={2025},
  publisher={Elsevier}
}

@article{bulivcek2015analysis,
  title={Analysis and approximation of a strain-limiting nonlinear elastic model},
  author={Bul{\'\i}{\v{c}}ek, M and M{\'a}lek, J and S{\"u}li, E},
  journal={Mathematics and Mechanics of Solids},
  volume={20},
  number={1},
  pages={92--118},
  year={2015},
  publisher={SAGE Publications Sage UK: London, England}
}

@article{itou2017nonlinear,
  title={Nonlinear elasticity with limiting small strain for cracks subject to non-penetration},
  author={Itou, Hiromichi and Kovtunenko, Victor A and Rajagopal, Kumbakonam R},
  journal={Mathematics and Mechanics of Solids},
  volume={22},
  number={6},
  pages={1334--1346},
  year={2017},
  publisher={Sage Publications Sage UK: London, England}
}

@article{rajagopal2003implicit,
  title={On implicit constitutive theories},
  author={Rajagopal, K R},
  journal={Applications of Mathematics},
  volume={48},
  number={4},
  pages={279--319},
  year={2003},
  publisher={Springer}
}

@article{rajagopal2021implicit,
  title={An implicit constitutive relation for describing the small strain response of porous elastic solids whose material moduli are dependent on the density},
  author={Rajagopal, KR},
  journal={Mathematics and Mechanics of Solids},
  volume={26},
  number={8},
  pages={1138--1146},
  year={2021},
  publisher={SAGE Publications Sage UK: London, England}
}

@article{treloar1975physics,
  title={The physics of rubber elasticity},
  author={Treloar, LR G},
  year={1975},
  publisher={OUP Oxford}
}

@article{arruda1993three,
  title={A three-dimensional constitutive model for the large stretch behavior of rubber elastic materials},
  author={Arruda, Ellen M and Boyce, Mary C},
  journal={Journal of the Mechanics and Physics of Solids},
  volume={41},
  number={2},
  pages={389--412},
  year={1993},
  publisher={Elsevier}
}

@article{horgan2004constitutive,
  title={Constitutive models for compressible nonlinearly elastic materials with limiting chain extensibility},
  author={Horgan, Cornelius O and Saccomandi, Giuseppe},
  journal={Journal of Elasticity},
  volume={77},
  number={2},
  pages={123--138},
  year={2004},
  publisher={Springer}
}

@article{ali2010review,
  title={A review of constitutive models for rubber-like materials},
  author={Ali, Aidy and Hosseini, Maryam and Sahari, Barkawi Bin and others},
  journal={American Journal of Engineering and Applied Sciences},
  volume={3},
  number={1},
  pages={232--239},
  year={2010}
}

@article{deam1976theory,
  title={The theory of rubber elasticity},
  author={Deam, RT and Edwards, Samuel Frederick},
  journal={Philosophical Transactions of the Royal Society of London. Series A, Mathematical and Physical Sciences},
  volume={280},
  number={1296},
  pages={317--353},
  year={1976},
  publisher={The Royal Society London}
}

@article{rivlin1948large,
  title={Large elastic deformations of isotropic materials. I. Fundamental concepts},
  author={Rivlin, RSl},
  journal={Philosophical Transactions of the Royal Society of London. Series A, Mathematical and Physical Sciences},
  volume={240},
  number={822},
  pages={459--490},
  year={1948},
  publisher={The Royal Society London}
}

@article{treloar1944stress,
  title={Stress-strain data for vulcanized rubber under various types of deformation},
  author={Treloar, Leslie RG},
  journal={Rubber Chemistry and Technology},
  volume={17},
  number={4},
  pages={813--825},
  year={1944}
}

@article{jones1975properties,
  title={The properties of rubber in pure homogeneous strain},
  author={Jones, DF and Treloar, LRG},
  journal={Journal of Physics D: Applied Physics},
  volume={8},
  number={11},
  pages={1285},
  year={1975},
  publisher={IOP Publishing}
}

@article{gent1996new,
  title={A new constitutive relation for rubber},
  author={Gent, Alan N},
  journal={Rubber chemistry and technology},
  volume={69},
  number={1},
  pages={59--61},
  year={1996},
  publisher={Rubber Division, ACS}
}

@article{horgan2002constitutive,
  title={Constitutive modelling of rubber-like and biological materials with limiting chain extensibility},
  author={Horgan, Cornelius O and Saccomandi, Giuseppe},
  journal={Mathematics and mechanics of solids},
  volume={7},
  number={4},
  pages={353--371},
  year={2002},
  publisher={Sage Publications Sage CA: Thousand Oaks, CA}
}

@article{panahi2025data,
  title={Data-driven model discovery with Kolmogorov-Arnold networks},
  author={Panahi, Shirin and Moradi, Mohammadamin and Bollt, Erik M and Lai, Ying-Cheng},
  journal={Physical Review Research},
  volume={7},
  number={2},
  pages={023037},
  year={2025},
  publisher={APS}
}

@article{carneros2024comparison,
  title={A comparison between Multilayer Perceptrons and Kolmogorov-Arnold Networks for multi-task classification in sitting posture recognition},
  author={Carneros-Prado, David and Caba{\~n}ero-G{\'o}mez, Luis and Johnson, Esperanza and Gonz{\'a}lez, Iv{\'a}n and Fontecha, Jes{\'u}s and Herv{\'a}s, Ram{\'o}n},
  journal={IEEE Access},
  year={2024},
  publisher={IEEE}
}

@inproceedings{essahraui2025kolmogorov,
  title={Kolmogorov—arnold networks: Overview of architectures and use cases},
  author={Essahraui, Siham and Lamaakal, Ismail and El Makkaoui, Khalid and Ouahbi, Ibrahim and Bouami, Mouncef Filali and Maleh, Yassine},
  booktitle={2025 International Conference on Circuit, Systems and Communication (ICCSC)},
  pages={1--6},
  year={2025},
  organization={IEEE}
}

@article{lee2025constitutive,
  title={A constitutive neural network for incompressible hyperelastic materials},
  author={Lee, Sanghee and Bathe, Klaus-J{\"u}rgen},
  journal={Machine Learning for Computational Science and Engineering},
  volume={1},
  number={2},
  pages={31},
  year={2025},
  publisher={Springer}
}

@book{marsden1994mathematical,
  title={Mathematical foundations of elasticity},
  author={Marsden, Jerrold E and Hughes, Thomas JR},
  year={1994},
  publisher={Courier Corporation}
}

@article{truesdell1952mechanical,
  title={The mechanical foundations of elasticity and fluid dynamics},
  author={Truesdell, Clifford},
  journal={Journal of Rational Mechanics and Analysis},
  volume={1},
  pages={125--300},
  year={1952},
  publisher={JSTOR}
}

\end{document}